\documentclass[11pt,reqno]{amsart}
\usepackage[a4paper,margin=1.15in]{geometry}
\usepackage{amsfonts}
\usepackage{amsmath,amssymb,amsthm,mathtools}
\usepackage{mathrsfs}
\usepackage{dsfont}
\usepackage{microtype}
\usepackage[colorlinks=true,linkcolor=blue,citecolor=magenta,urlcolor=blue]{hyperref}
\usepackage[nameinlink,noabbrev]{cleveref}

\numberwithin{equation}{section}
\newtheorem{theorem}{Theorem}[section]
\newtheorem{proposition}[theorem]{Proposition}
\newtheorem{conjecture}{Conjecture}
\newtheorem{lemma}[theorem]{Lemma}

\theoremstyle{definition}
\newtheorem{definition}[theorem]{Definition}

\newcommand{\dd}{\mathrm{d}}

\newcommand{\Ldz}{\frac{\zeta'}{\zeta}}

\newcommand{\mycomment}[1]{}

\title[Almost primes in short intervals and averages of long Dirichlet polynomials]{Almost primes in short intervals and averages of long Dirichlet polynomial approximations to $(\zeta'/\zeta)^{k}$}

\author{S. Baluyot, J.B. Conrey, A. Irfan and N.C. Snaith}
\address{School of Mathematics, University of Bristol, [Bristol, United Kingdom]}
\date{\today}

\begin{document}

\subjclass[2020]{60B20, 11M06, 11N37, 11M50}
\keywords{Riemann zeta function, random matrix theory, logarithmic derivative, Dirichlet polynomials, shifted convolution, von Mangoldt, characteristic polynomial}

\begin{abstract}
We study the asymptotic behaviour of the mean square of long truncations of the Dirichlet series for $\bigl(\zeta'/\zeta\bigr)^{k}$ on the critical line; we establish the connection between this problem and the short-interval variance of the convoluted von Mangoldt function $\Lambda^{*k}$ in short intervals and in arithmetic progressions. Under a shifted moments conjecture and the Riemann Hypothesis (RH), we prove that the short-interval variance leading coefficient is the derivative of the mean-square leading coefficient. A unitary trace-sum average governs the variance of sums of $\Lambda^{*k}$ over monic polynomials in
$\mathbb F_q[T]$; we denote its leading coefficient by $\delta_k(c)$, where the total trace degree and matrix size grow with a fixed ratio $c$. Under a weaker shifted moments conjecture and RH, we recover this $\delta_k(c)$ as the leading coefficient of the variance of short-interval sums of $\Lambda^{*k}$ over the integers. It is the counterpart of $\gamma_k(c)$, the piecewise polynomial from the variance of divisor sums, with the Dirichlet coefficients of $\zeta(s)^k$ replaced by those of $\bigl(-\zeta'/\zeta\bigr)^k$.
\end{abstract}
\maketitle
\section{Introduction}\label{sec:intro}

This paper belongs to a series of works relating almost-prime statistics to moments of logarithmic derivatives of $\zeta(s)$, motivated by function-field analogues that we have analysed using random matrix theory. The asymptotics for the unitary trace-sum average governing the function-field short-interval variance of sums of $\Lambda^{*k}$ \cite{kuperberglalin} were studied in \cite{kn:irfan26-sums}, while here we recover the predicted leading coefficients for long Dirichlet polynomial mean squares and almost-prime variances over the integers, under shifted moment conjectures and RH. Forthcoming work \cite{kn:bcis26-general-swap} examines the exact finite-$N$ moments generating these trace-sum averages through their swap contributions and their Painlev\'e $\tau$-function representations, refining conjectures for shifted moments of $\zeta'/\zeta$ and almost-prime variances in short intervals and arithmetic progressions.

\subsection{Averages of long polynomials}
The moments of the Riemann zeta function $\zeta(s)$ and its derivatives evaluated on the critical line offer insight into large values of $\zeta(s)$, its zero statistics, and even the distribution of primes. In particular, we focus on the logarithmic derivative
\begin{equation}
      \frac{\zeta'}{\zeta}(s) = -\sum_{n=1}^{\infty} \frac{\Lambda(n)}{n^{s}}, \qquad \Re{s}> 1.
 \end{equation}
Here, $\Lambda(n)$ is the von Mangoldt function defined as
 \begin{align}
 \Lambda(n) & := \begin{cases} 
   \log{p}
  &  \exists \, m \in \mathbb{N }, \, p \in \mathbb{P}: n = p^{m} \\
      0
     & \text{otherwise.}
   \end{cases}
\end{align}
This may be interpreted as a weighted counting function for prime powers. As the meromorphic continuation of $\zeta'/ \zeta(s)$ has simple poles at the zeros of $\zeta(s)$, the moments are evaluated to the right of the critical line along which the non-trivial zeros lie under RH. For a fixed $a>0$, Conjecture 2 in \cite{kn:bbbcprs19} implies there exists a function $c_k : \mathbb{R}_{+} \rightarrow \mathbb{R}$ such that
\begin{equation}
\label{eq:moments1}
\frac{(2a)^{2k-1}}{T} \int_{0}^{T}
\bigg | \frac{\zeta'}{\zeta}
\left(\tfrac{1}{2}+\tfrac{a}{\log T}+it\right)\bigg |^{2k} dt
\sim c_k(a)(\log T)^{2k},
\qquad k\in\mathbb{N}, \quad a>0,
\end{equation}
as $T\to\infty$. For fixed $k \in \mathbb{N}$, such moment asymptotics are consistent with the Conrey–Snaith \cite{kn:consna08} form of the ratios conjecture \cite{kn:cfz2} under RH. Furthermore, the same conjecture in \cite{kn:bbbcprs19} also predicts
\begin{equation}
\label{eq:c-limiting-value}
c_k(a)\sim {2k-2\choose k-1},
\end{equation}
as $a\to0^+$. The same small-$a$ limit in the dyadic setting, i.e. with \eqref{eq:moments1} averaged over $[T,2T]$, was subsequently established by Ge \cite{Ge3} under both RH and the $2k$-tuple Essential Simplicity hypothesis. Several results towards the fixed-$a$ moment problem are also known. Selberg \cite{kn:selberg43} obtained an upper bound for the corresponding second moment in the range $a \geq 10$ under RH. Farmer et al. \cite{kn:fargonll13} investigated mixed shifted generalisations of the dyadic form of \eqref{eq:moments1}, relating them to almost primes and higher correlations of the zeros of $\zeta(s)$. Later, Rodgers \cite{kn:rodgers17} proved
\begin{equation}
\frac{1}{T}\int_T^{2T}
\left|
\frac{\zeta'}{\zeta}
\left(\frac{1}{2}+\frac{a}{\log T}+it\right)
\right|^{2k}dt
\ll_{a,k}(\log T)^{2k},
\end{equation}
for every fixed $k\in\mathbb{N}$ and $a>0$ under RH alone. Thus, under RH, Rodgers obtains the conjectured order of magnitude as an upper bound, but not the asymptotic constant appearing in \eqref{eq:moments1}. A standard method for studying such moments is to replace the relevant Dirichlet series by a long Dirichlet polynomial and to study its mean values in $t$-aspect. This viewpoint was sharpened by Bettin and Conrey \cite{bettinconrey} who obtained precise asymptotics for mean squares of long Dirichlet polynomials, in their case for the mean square for $\zeta(s)^{k}$ with coefficients $d_k(n)$, under shifted moment conjectures. Here $d_{k}(n)$ is the $k$th divisor function; one may interpret $d_k(n)$ as the number of ways to write a positive integer
$n$ as a product of $k$ positive integers. 

Explicitly, one has $d_k(n)=(1*\cdots*1)(n)$ which is a $k$-fold \textit{Dirichlet convolution} of the series coefficient for $\zeta(s)$. For two arithmetic functions $f$ and $g$, one defines their Dirichlet convolution $f * g$ by
\begin{equation*}
    (f*g)(n) = \sum_{d | n} f(d) g\biggl(\frac{n}{d}\biggr),
\end{equation*}
where the sum extends over all positive divisors $d$ of $n$. Note that we use $*$ to denote a Dirichlet convolution throughout. 

In studying the moments
\begin{equation}
 \label{eq:moments}
    \frac{1}{T} \int_{0}^{T} \bigg | \frac{\zeta'}{\zeta} (\tfrac{1}{2} + \kappa + it)\bigg |^{2k} dt, 
\end{equation}
for $\kappa = \frac{a}{\log{T}}$ with $a>0$ fixed, we consider the long Dirichlet polynomials of $\bigl(-\zeta'/\zeta(s)\bigr)^{k}$ whose coefficients are the $k$-fold convolution of $\Lambda(n)$ accordingly. For $k\in \mathbb{N}$, we have
\begin{equation}
\Lambda^{*k} := \underbrace{\Lambda * \cdots * \Lambda \mathstrut}_{k \text{ times}} \, ,
\end{equation}
and the Dirichlet series
\begin{equation}\label{eq:dirichlet-series}
\Bigl(-\frac{\zeta'}{\zeta}(s)\Bigr)^{k} = \sum_{n=1}^\infty \frac{\Lambda^{*k}(n)}{n^s},
\qquad (\Re s>1).
\end{equation}
Let $N\geq 1$ be a real parameter, and fix $s\in\mathbb{C}$. 
For any $\sigma>\max(1,\Re s)$, Perron's formula gives
\begin{equation}\label{eq:perron-lambda-k}
\sum_{n\leq N}\hspace{-0.03 in}^{'}\frac{\Lambda^{*k}(n)}{n^s}
=
\frac{1}{2\pi i}
\int_{(\sigma)}
\left(-\frac{\zeta'}{\zeta}(w)\right)^k
\frac{N^{w-s}}{w-s}\,dw .
\end{equation}
The prime on the sum denotes the endpoint term $n=N$ being counted with weight
$1/2$ if $N$ is an integer and the contour $(\sigma)$ denotes integration along the vertical line
$\Re w=\sigma$.
By shifting the $w$-contour to cross the pole at $w=1$, we obtain a
main residue term given by
\begin{equation}
    \text{Res}_{w=1} \bigg ( \Big ( -\frac{\zeta'}{\zeta}(w)\Big )^{k} \frac{N^{w-s}}{w-s} \bigg ).
\end{equation}

Motivated by \eqref{eq:moments1}, we set $s= \tfrac{1}{2} + it$ and examine the mean square of the truncated series after subtracting the residue term 
\begin{align}
\label{Ik-truncated-moments}
    \mathcal{I}_{k}(T,N) &:= \frac{1}{T} \int_{0}^{T}  \bigg | \sum_{1 \leq n \leq N} \frac{\Lambda^{*k}(n)}{n^{1/2 + it} }- \text{Res}_{w=1}\biggl( \Big ( -\frac{\zeta'}{\zeta}\bigl( w \bigr) \Big )^{k} \frac{N^{w-\frac{1}{2} - it}}{w-\frac{1}{2}  -it}\biggr) \bigg |^{2} dt,
\end{align}
where $N = T^{\alpha}$ for fixed $\alpha \in \mathbb{R}_{+}$. In contrast to \eqref{eq:moments}, $\mathcal{I}_{k}(T,T^{\alpha})$ consists of a finite Dirichlet polynomial with an explicit residue correction, and can therefore be defined on the critical line. By analogy with  \eqref{eq:moments1}, for any fixed $\alpha > 0$, we expect $\mathcal{I}_{k}(T,T^{\alpha})$ to grow at the scale of $(\log{T})^{2k}$ as $T \rightarrow \infty$ \cite{kn:bbbcprs19}. This leads us to introduce
\begin{equation}
\label{eq:Mdef}
   \mathcal{M}_{k}(\alpha) := \lim_{T \rightarrow \infty} \frac{\mathcal{I}_{k}(T, T^{\alpha})}{(\log{T})^{2k}}, \qquad k \in \mathbb{N}, \alpha > 0,
\end{equation}
whenever the limit exists. Under Conjecture~\ref{eq:main-conj}, we will prove that $\mathcal{M}_k(\alpha)$ is a continuous piecewise polynomial function that satisfies similar smoothness properties to the leading coefficients for the mean square of the truncated series of $\zeta(s)^{k}$ \cite{bettinconrey}.

As an aside for readers familiar with the associated literature, the mean-square problem for the case $k=1$ is related to the pair correlation of the Riemann zeros. Pair correlation describes differences between zero ordinates on the scale of the mean spacing, including differences between non-consecutive zeros. Montgomery’s form factor $F(u,T)$ \cite{kn:mont73} represents these differences by a normalised, weighted Fourier sum. Under RH, Goldston, Gonek and Montgomery \cite{kn:ggm01} related the shifted second moment of $\zeta'/\zeta(s)$ near the critical line to an integral transform of $F(u,T)$. Goldston and Montgomery \cite{goldstonmont87} proved, also under RH, that the strong pair-correlation conjecture is equivalent to suitably uniform asymptotics for the variance of primes in short intervals. Farmer \cite{kn:farmer95} extended this perspective beyond pair correlation, proving an asymptotic for a third-order mixed moment of $\zeta'/\zeta(s)$ under RH and a three-level GUE hypothesis. Carneiro et al.\ \cite{kn:carneiro_ea22} later obtained improved RH-conditional bounds for the second moment of $\zeta'/\zeta(s)$, integrals of $F(u,T)$, and weighted prime variance in short intervals. Chirre \cite{kn:chirre22} also extended the Goldston--Gonek--Montgomery equivalence to second moments of the derivatives $(\zeta'/\zeta)^{(j)}$, with $j\geq1$, under RH.

Due to the cumulative dependence on $\alpha$ in the kernel $T^{\alpha z}/z$, 
\begin{equation}
\frac{T^{\alpha z}-1}{z}=(\log T)\int_0^\alpha T^{uz}\,du,
\end{equation}
we may take at a heuristic level that
\begin{equation}
\mathcal M_1(\alpha)=\int_0^\alpha f(u)\,du,
\end{equation}
where $f(u)$ denotes the limiting behaviour of $F(u, T)$ as $T \rightarrow \infty$, whenever this limit exists. A rigorous justification of this relation requires separate treatment of the contribution near $u=0$. Under RH, one has $f(u)=u$ for $0<u < 1$, while the strong pair-correlation conjecture predicts $f(u)=1$ for $u\ge 1$ \cite{kn:mont73}; hence
\begin{equation}
\label{m1alpha}
\mathcal M_1(\alpha)=
\begin{cases}
\frac{\alpha^2}{2}, & 0<\alpha < 1 \\[1mm]
\alpha-\tfrac12, & \alpha\ge 1,
\end{cases}
\end{equation}
under both assumptions. No explicit description appears to be available in the literature when $k \geq 2$. To obtain such a form for $\mathcal{M}_{k}(\alpha)$, we introduce a conjecture which mirrors the structure of a known theorem from random matrix theory on averages of logarithmic derivatives of characteristic polynomials (stated in \cite{kn:bbbcprs19}, equation 43) and we also derive an unsymmetrised version of this expression, given in equation \eqref{eq:unsymmetrised-form} in Section \ref{sec:Mk-derivation}. Based on the observed structure in these results, we conjecture the following form for shifted moments of the logarithmic derivative of $\zeta(s)$.
\begin{conjecture} 
\label{eq:main-conj} Let $k \in \mathbb{N}$, and let $Y = Y(T) \geq 1$ be a parameter. For $T > 1$, $\Re z_1 > 0$ and a fixed $0 < \xi < 1$, one has
\begin{align}
&\frac{1}{T}\int_0^T
\left(\frac{\zeta'}{\zeta}\bigl(\tfrac12+it+z_1\bigr)\right)^k
\left(\frac{\zeta'}{\zeta}\bigl(\tfrac12-it-z_2\bigr)\right)^k dt \\ \nonumber =& \frac{1}{k!^2 (2\pi i)^{2k} }  \int_{\mathcal C_u^k}\int_{\mathcal C_v^k} T^{\frac 12 (\sum_{j=1}^{k}u_j - \sum_{\ell = 1}^{k} v_\ell)} \prod_{j,\ell =1}^{k} (v_{\ell} - u_j - z)\prod_{1 \leq i < j \leq k}(u_i - u_j)^{2} \hspace{-1mm} \prod_{1 \leq i < j \leq k}(v_i - v_j)^{2} \\  \nonumber & \quad \qquad \times  \biggl( \sum_{j=1}^{k}\frac{1}{u_j}
   \;+\;
   \sum_{\ell=1}^{k}\frac{1}{v_\ell -z}
   \;-\;\frac{\log{T}}{2} \biggr)^{k} 
    \biggl( \, \sum_{j=1}^{k}\frac{1}{u_j + z}
   \;+\;
   \sum_{\ell=1}^{k}\frac{1}{v_\ell}
   \;+\;\frac{\log{T}}{2} \biggr)^{k} z^{-k^2}\\ \nonumber
& \quad \qquad \times \, \prod_{j=1}^{k} u_j^{-k} du_j \prod_{\ell = 1}^{k} v_{\ell}^{-k} dv_\ell  + O_{k, c, \, \xi}((\log{T})^{2k-1+\xi}(1+\min(|z|,1)\log{T})^{-\theta_{k,c}})
\end{align}
where $z = z_1 - z_2$ for $|\Im z|, |\Im z_2| \le Y$ and $0< -\Re{z_2}, \Re z \le \frac{3}{\log T}$. Here  $\mathcal C_u = \mathcal C_u(z)$ and $\mathcal C_v = \mathcal C_v(z)$ are positively oriented piecewise smooth contours, with $\mathcal C_u$ encircling once $0$ and $-z$, and $\mathcal C_v$ encircling once $0$ and $z$, chosen to avoid the poles of the integrand. Also $\theta_{k,c} := k-1$ for $k\geq 3$ and $\theta_{k,c} := 1+c$ for $k=1,2$ for some fixed $c>0$.
For $k=1,2$, we take
$0<\xi<c \leq 1$. 
\end{conjecture}
As detailed in the proof of Theorem~\ref{eq:rmt-theorem}, one can also differentiate the contour integral form of the shifted ratios conjecture (stated in \cite{kn:cfz2}, equation 6.22) to obtain the above form. In fact, this approach suggests a stronger error term of $O_{k,\varepsilon}(T^{-\frac{1}{2} + \varepsilon}(\log{T})^{2k})$; we do not employ the sharper bound here, since our subsequent argument requires only an error term sufficiently small to isolate the leading coefficient of the moment polynomial. Accordingly, Conjecture~\ref{eq:main-conj} is formulated with the weakest error term needed for that purpose.
 
\subsection{Averages of von Mangoldt convolutions} The asymptotic behaviour of the variance of sums of Dirichlet coefficients is also the subject of extensive study \cite{kn:harsou17, kn:lester16}, in part due to its connection to moment problems for $\zeta(s)$ and other $L$-functions. In the proof of Theorem~2 in \cite{bettinconrey}, Bettin and Conrey analyse the variance of the divisor function in short intervals through its relationship with moments of $\zeta(s)^k$. In this paper, we consider the analogous problem for the $k$-fold von Mangoldt convolution $\Lambda^{*k}(n)$ in relation to moments of $(-\zeta'/\zeta(s))^k$. For $x, H > 0$ where $x \in [X, 2X]$, we consider the sum
\begin{equation}\label{eq:psi}
 \Psi_k(x;H):=\sum_{x<n\le x+H}\Lambda^{*k}(n), 
\end{equation}
for which we are interested in the short-interval regime $H = o(X)$. As before, we define the variance by subtracting the main contribution from the pole of $(-\zeta'/\zeta (s))^k$ at $s=1$. Let
\begin{align}\label{eq:def-Delta}
\Delta_{k}(x;H)
:=
\Psi_k(x;H)
-
\text{Res}_{s=1}
\left[
\left(-\frac{\zeta'}{\zeta}(s)\right)^k \frac{(x+H)^s - x^s}{s}
\right].
\end{align}
We consider
\begin{equation}
\label{def:variance}
    \mathcal{V}_{k}(X,H) = \frac{1}{XH} \int_{X}^{2X} | \Delta_{k}(x;H)|^{2} \, dx
\end{equation}
for $H=X^{1-1/\alpha}$ with a real parameter $\alpha > 0$. When $k=1$, the mean square is a normalised dyadic form of the classical variance of primes in short intervals studied by Goldston-Montgomery \cite{goldstonmont87}. For general $k \in \mathbb{N}$, it should be viewed as the short-interval variance attached to the coefficients of 
$(-\zeta'/\zeta(s))^k$, and hence as a higher-order weighted almost-prime statistic. Note that almost-prime is understood in the generalised sense of Farmer et al. \cite{kn:fargonll13}, as an integer $n$ expressible as a product of $k$ \textit{prime powers}. For $\alpha > 0$ and $k \in \mathbb{N}$, let
\begin{equation}
\label{eq:Dvariance}
    \mathcal{D}_{k}(\alpha) := \lim_{X \rightarrow \infty} \frac{\mathcal{V}_{k}(X,X^{1-1/\alpha})}{(\log{X^{1/\alpha}})^{2k-1}}.
\end{equation}
Bettin and Conrey \cite{bettinconrey} showed that, in the case of $d_k(n)$, the analogous leading variance term is obtained by differentiating the leading moment polynomial, i.e their equivalent $\mathcal{M}_k(\alpha)$. We prove that such a derivative relationship holds between the leading coefficient from the moments of $(-\zeta'/\zeta (s))^k$ and the short-interval variance of $\Lambda^{*k}(n)$, as stated in Theorem~\ref{conj:derivative}.
\subsection{Averages of von Mangoldt convolutions in arithmetic progressions}
We now examine this problem in arithmetic progressions. For $q \in \mathbb{N}, X > 0$, define
\begin{align}
\label{def:W_k}
    \mathcal{W}_{k}(X; q) := \frac{1}{X} \sum_{\substack{1 \leq a \leq q \\ (a,q) = 1}} \bigg | \sum_{\substack{n \leq X \\ n \equiv a \, (\text{mod} \, {q})}} \Lambda^{*k}(n) - \frac{1}{\varphi(q)} \text{Res}_{s=1} \bigg ( \Big ( -\frac{\zeta_{q}'}{\zeta_{q}}(s)\Big )^{k} \frac{X^{s}}{s} \bigg ) \bigg |^2,
\end{align}
where $\zeta_q (s) := \sum_{n \geq 1, (n, q) = 1} n^{-s}$ and $\varphi$ is Euler's $\varphi$-function. $\mathcal{W}_{k}(X; q)$ represents a weighted higher-order almost-prime variance. To recover the leading coefficient for $\mathcal{W}_k(X,q)$, we define
\begin{align}
\label{eq:arithemtic-variance}
\mathcal{Q}_{k}(\alpha) :=
\lim_{q\to\infty}
\frac{\mathcal{W}_{k}(X;q)}{(\log q)^{2k-1}},
\qquad \alpha>0,
\end{align}
whenever the limit exists, where $q$ tends to infinity through the primes and $X=X(q)$ satisfies $\log X/\log q\to\alpha$. We restrict to prime moduli $q = p$. When $k=1$, the inner sum is $\psi(X;q,a)$, the Chebyshev function in an arithmetic progression, and its variance was studied by Friedlander–Goldston \cite{friedlander}. They obtained a lower bound for $k=1$, whereas the higher-$k$ case remains largely unexplored. For $k > 1$, Knafo \cite{kn:knafo06} proved lower bounds for the variance associated with the generalised von Mangoldt function $\Lambda_{k} = \mu * (\log)^k$, where $\mu$ is the Möbius function, rather than for the $k$-fold convolution $\Lambda^{*k}(n)$ that is considered here. 

For a Dirichlet character $\chi$ modulo $q$, let $L(s,\chi)=\sum_{n=1}^{\infty}\chi(n)n^{-s}$ for $\Re s>1$, with the same notation used for its meromorphic continuation. In the proof of our arithmetic-progression theorem, we assume the Generalised Riemann Hypothesis (GRH) for these Dirichlet $L$-functions. To determine $\mathcal{Q}_{k}(\alpha)$, we require an analogue of Conjecture \ref{eq:main-conj} which is an average over $L$-functions instead of a $t$-aspect average.
\begin{conjecture}
\label{eq:ratios-arithmetic-corrected}
Let $q$ be prime, $k \in \mathbb{N}$ and $Y = Y(q) \geq 1$ be a parameter. For $\Re z_1 > 0$ and a fixed $0 < \xi < 1$, as $q \rightarrow \infty$ one has
\begin{align}
&\frac{1}{\varphi(q)}\sum_{\chi\bmod q}
\left(\frac{L'}{L}\left(\tfrac12+z_1,\chi\right)\right)^k
\left(\frac{L'}{L}\left(\tfrac12-z_2,\overline\chi\right)\right)^k
\nonumber \\
&\quad =
\frac{1}{k!^2(2\pi i)^{2k}}
\int_{\mathcal C_u^k}\int_{\mathcal C_v^k}
q^{\frac12\left(\sum_{j=1}^k u_j-\sum_{\ell=1}^k v_\ell\right)}
\prod_{1\leq i<j\leq k}(u_i-u_j)^2
\prod_{1\leq i<j\leq k}(v_i-v_j)^2
\nonumber \\
&\qquad\times
\left(
\sum_{j=1}^k \frac{1}{u_j}
+
\sum_{\ell=1}^k \frac{1}{v_\ell-z}
-
\frac{\log q}{2}
\right)^k
\left(
\sum_{j=1}^k \frac{1}{u_j+z}
+
\sum_{\ell=1}^k \frac{1}{v_\ell}
+
\frac{\log q}{2}
\right)^k
\prod_{j,\ell=1}^k (v_\ell-u_j-z)
\nonumber \\
&\qquad\times
z^{-k^2}  \prod_{j=1}^k u_j^{-k}\,du_j
\prod_{\ell=1}^k v_\ell^{-k}\,dv_\ell
\nonumber +
O_{k,\xi}\left(
(\log q)^{2k-1+\xi}
\left(1+\min(|z|,1)\log q\right)^{-\theta_{k,c}}
\right),
\end{align}
where $z = z_1 - z_2$ for $|\Im z|, |\Im z_2| \le Y$ and $0< -\Re{z_2}, \Re z \le \frac{3}{\log q}$. Here $\mathcal C_u = \mathcal C_u(z)$ and $\mathcal C_v = \mathcal C_v(z)$ are positively oriented piecewise smooth contours, with $\mathcal C_u$ encircling once $0$ and $-z$, and $\mathcal C_v$ encircling once $0$ and $z$, chosen to avoid the poles of the integrand. The error exponent $\theta_{k,c}$ is as defined in Conjecture~\ref{eq:main-conj}.
\end{conjecture}
\subsection{Random Matrix Theory}
A central bridge in the study of the Riemann zeta function is the analogy with characteristic polynomials of random unitary matrices, which provides a framework for predicting the moments of $\zeta(s)$ and its derivatives. The function-field setting gives a precise realisation of this analogy, translating variance problems for arithmetic functions into matrix integrals over classical compact groups. 

Over the polynomial ring $\mathbb{F}_q[t]$, Keating, Rodgers, Roditty-Gershon and Rudnick \cite{kr3} studied the variance of short-interval sums of $d_k$, i.e. the divisor function on monic polynomials. They showed that the leading coefficient in its large-$q$ limit is given by a unitary matrix integral involving products of secular coefficients \cite[Theorem~1.2]{kr3}. Using similar methods, Kuperberg and Lal\'{i}n \cite{kuperberglalin} prove that the matrix integral arising from the variance of short-interval sums of $\Lambda^{*k}$ is a unitary matrix integral of trace sums \cite[Theorem~1.3]{kuperberglalin}. Such matrix integrals are closely related to moments of characteristic polynomials and their derivatives, as we make explicit below for the $\Lambda^{*k}$ case.

Let $U \in U(N)$ be a unitary matrix. The eigenvalues of $U$ lie on the unit circle and can be written as $e^{i\theta_1},\dots,e^{i\theta_N}$ for $\theta_1,\dots,\theta_N \in [0, 2\pi)$. The characteristic polynomial of $U$ is defined as
\begin{equation}
\Lambda_U(z) := \det(I - zU^{*})
= \prod_{n=1}^N \bigl(1 - z e^{-i\theta_n}\bigr),
\end{equation}
with $U^{*}$ denoting the conjugate transpose of $U$. For $|z|<1$, the logarithmic derivative of $\Lambda_U$ admits the expansion
\begin{equation}
  \frac{\Lambda'_U}{\Lambda_U}(z)
  = - \sum_{n=1}^N \frac{e^{-i\theta_n}}{1 - z e^{-i\theta_n}}
  = -\frac{1}{z}\sum_{j=1}^{\infty} z^{j}\,\mathrm{Tr}(U^{*j}),
\end{equation}
where $\mathrm{Tr}(U^{*j}) := \sum_{n=1}^{N} e^{-ij \theta_n}$ is the trace of $U^{*j}$. Thus, the moments of the logarithmic derivative of $\Lambda_{U}(z)$ over the unitary group serve as a generating function for the mean square of trace sums. For $|x|,|y|<1$, we will show that
\begin{align}
\label{line:moments-haar-invariance}
F_k(xy; N)
&:=
\int_{U(N)}
\left(-x\frac{\Lambda'_U}{\Lambda_U}(x)\right)^k
\left(-y\frac{\Lambda'_{U^*}}{\Lambda_{U^*}}(y)\right)^k\,dU
\\
\label{line:haar-invariance}
&=
\sum_{n=k}^{\infty}(xy)^n I_k(n;N),
\end{align}
where
\begin{align}
\label{eq:Ik}
    I_{k}(n;N)
    = \int_{U(N)} \Biggl| \sum_{\substack{i_1 + \cdots + i_k = n \\ i_{1}, \ldots, i_{k} \geq 1}}\text{Tr}(U^{i_{1}}) \cdots \text{Tr}(U^{i_{k}}) \, \Biggr|^{2}  dU.
\end{align}
Here, $dU$ denotes the normalised Haar measure. The matrix integral that Kuperberg and Lal\'in obtain through their function-field analysis is precisely in the form of \eqref{eq:Ik} with matrix size $N = n -  h -2$ where $n$ is the total trace degree and $h$ parametrises the short interval. In the asymptotic problem, we fix a real number $\alpha>0$ and take
\begin{equation*}
        n=\lfloor \alpha N\rfloor .
\end{equation*}
The floor $\lfloor \cdot \rfloor$ is required as the total degree $n$ must be an integer; it affects only lower-order terms in the asymptotic expansion. We define the local trace-average density by
\begin{equation}
\label{eq:delta-k-trace-def}
        \delta_k(\alpha)
        :=
        \lim_{N\to\infty}
        \frac{I_k(\lfloor \alpha N\rfloor;N)}{N^{2k-1}},
        \qquad \alpha >0,
\end{equation}
whenever the limit exists. Note that expanding the $k$th powers in \eqref{line:moments-haar-invariance} produces ordered $k$--tuples
$(i_1,\dots,i_k)$ and $(j_1,\dots,j_k)$ with $i_a,j_b\ge 1$.  By Haar invariance under $ U \rightarrow e^{i \varphi} U$ for $ \varphi \in[0,2 \pi)$, the integral vanishes unless $\sum_{a=1}^k i_a=\sum_{b=1}^k j_b=:n$ giving \eqref{line:haar-invariance}.
We define the cumulative, or cut-off, trace average 
\begin{equation}
\label{I-tilde}
    \tilde{I_{k}}(n;N) := \int_{U(N)} \Biggl| \sum_{\substack{j_{1} + \cdots + j_{k} \leq n \\ j_{1}, \ldots, j_{k} \geq 1}} \text{Tr}(U^{j_{1}}) \cdots \text{Tr}(U^{j_{k}}) \Biggr|^{2}dU.
\end{equation}
 We may relate these matrix integrals through
\begin{equation}
\label{eq:discrete-relationship}
\tilde I_{k}(n;N)=\sum_{m=k}^{n}I_k(m,N).
\end{equation}
Therefore
\begin{equation}
\sum_{n\geq k}\tilde I_{k}(n;N)t^n
=
\frac{1}{1-t}\sum_{m\geq k}I_k(m,N)t^m
=
\frac{F_k(t;N)}{1-t},
\end{equation}
and consequently,
\begin{equation}
\tilde{I_{k}}(n;N)
=
\frac{1}{2\pi i}
\oint_{|t|=r}
\frac{F_k(t;N)}{(1-t)t^{n+1}}\,dt,
\end{equation}
for any $0<r<1$. Equivalently, we may set $t=e^{-s/N}$ where $s=\sigma+i\theta$; on this contour, $r=e^{-\sigma/N}$ and the circle $|t|=r$ corresponds to a vertical segment
\begin{equation*}
-\pi N\leq \theta\leq \pi N .
\end{equation*}
In this case, the Cauchy coefficient extraction becomes a discrete Laplace inversion
\begin{equation}
\tilde I_{k}(n;N)
=
\frac{1}{2\pi i N}
\int_{\sigma-i\pi N}^{\sigma+i\pi N}
e^{n s/ N}\,
\frac{F_k(e^{-s/N};N)}{1-e^{-s/N}}\,ds.
\end{equation}
The substitution $t=e^{-s/N}$ zooms into the neighbourhood of $t=1$,
since $e^{-s/N}=1-s/N+O(N^{-2})$ for bounded $s$. We have a finite vertical segment as $e^{-s/N}$ is
$2\pi iN$-periodic in $s$; one such period maps exactly onto the original Cauchy circle. This is analogous to the use of Perron's formula \eqref{eq:perron-lambda-k} near the half-line: the corresponding change of variables puts the generating function in the scale where the cutoff parameter is visible. Therefore we study the large-\(N\) asymptotics of
$\tilde I_k(\lfloor \alpha N\rfloor,N)$, the normalised limit for which is the matrix analogue of $\mathcal M_k(\alpha)$. Let
\begin{equation}
\label{eq:rmt-mkalpha}
M_{k}(\alpha) := \lim_{N \rightarrow \infty} \frac{\tilde{I_{k}}(\lfloor \alpha N \rfloor,N)}{N^{2k}}, \qquad k \in \mathbb{N}, \alpha > 0.
\end{equation}
Due to \eqref{eq:discrete-relationship}, we can write
\begin{equation}
\label{eq:Itilde-riemann-sum}
        \frac{\tilde I_k(\lfloor \alpha N\rfloor,N)}{N^{2k}}
        =
        \frac1N
        \sum_{m=k}^{\lfloor \alpha N\rfloor}
        \frac{I_k(m,N)}{N^{2k-1}}.
\end{equation}
By \cite[Proposition~5.2]{kn:irfan26-sums}, passage to the limit in this Riemann sum gives
\begin{equation}
\label{eq:Mk-delta-integral}
        M_k(\alpha)=\int_0^\alpha \delta_k(c)\,dc.
\end{equation}
As $\delta_k$ is continuous on $(0,\infty)$, it follows that $M_k$ is continuously differentiable there, with
\begin{equation}
\label{eq:derivative-relationship-intro}
M_k'(\alpha) = \delta_k(\alpha),
\qquad \alpha>0.
\end{equation}
To determine $M_k(\alpha)$, we first study the range in which the total degree $n$ does not exceed the matrix size $N$; for \eqref{eq:Ik}, the total degree is $n=i_1+\cdots+i_k$. Simple formulae obtained from symmetric-function orthogonality apply when this degree is at most the matrix size $N$. This phenomenon is observed in averages of products of secular coefficients, as well as the present traces setting. When $n \leq N$, the Diaconis--Gamburd formula \cite{kn:diagam04} gives these evaluations for products of secular coefficients, while the Diaconis--Shahshahani formula \cite{kn:diasha} gives the evaluations for trace products. Beyond this range, one must retain the finite-$N$ restrictions: the character sums include only partitions with at most $N$ parts, so some partitions of $n$ are excluded when $n>N$. Thus, under the scaling $n=\lfloor\alpha N\rfloor$, the behaviour of $M_k(\alpha)$ for $\alpha\leq1$ is understood, while $\alpha>1$ requires an analysis of the restricted character sums as in \cite{kn:irfan26-traces}. Generally, one expects $M_k(\alpha)$ to be piecewise polynomial, with candidate breakpoints occurring at integer values of $\alpha$.
\section{Results}
We now present the main results of this paper. We begin by introducing a polynomial that appears in the leading-order moment asymptotics for both $\zeta'/\zeta$, under Conjecture~\ref{eq:main-conj}, and the logarithmic derivatives of characteristic polynomials over $U(N)$.  RH is assumed throughout the results stated, and the generalised Riemann Hypothesis (GRH) is assumed for Theorem~\ref{conj:arithmetic-variance}.

\begin{definition}
Fix $\tau \in (0,\frac{1}{2})$ and let $\alpha \in \mathbb{R}_{+}$. For an integer $0 \leq r \leq k$, let
\begin{align} \label{eq:prk}
    P_{r,k}(&\alpha) := \frac{1}{(2\pi i)^{2k}}  \sum_{m=0}^{k} (-1)^{k-m}{k \choose m}\frac{1}{(k + m)!} \\ & \qquad \, \nonumber \times \int_{c_{j}(r)}  \int_{\hat{c}_{\ell}(r)}  \biggl( \sum_{j=1}^{k}x_j-\alpha \biggr)^{k +m} \Delta(x)^{2} \Delta(y)^{2} \prod^{k}_{i,j=1}(x_i - y_j - 1)   \, \prod_{j=1}^{k}x_{j}^{-k}  \, \prod_{\ell=1}^{k}y_{\ell}^{-k} 
 \\ & \nonumber
 \qquad \qquad \qquad \quad \times  \, \biggl(\sum_{j=1}^{k}\frac{1}{x_j}
   \;+\; 
   \sum_{\ell=1}^{k}\frac{1}{y_\ell + 1} \biggr)^{m} 
  \biggl(\sum_{j=1}^{k}\frac{1}{x_j  -1}
   \;+\;
   \sum_{\ell=1}^{k}\frac{1}{y_\ell  } \biggr)^{k} 
  \, \prod_{j=1}^{k} dx_j\; \prod_{\ell =1}^{k} dy_\ell,
\end{align}
where $\Delta(u) = \prod_{i < j} (u_j - u_i)$, $c_{j}(r)$ are counter-clockwise circles of radius $\tau$ centered at 1 for $j=1, \ldots, r$ and at 0 otherwise; similarly, $\hat{c}_{\ell}(r)$ are counter-clockwise circles of radius $\tau$ centered at $-1$ for $\ell = 1, \ldots, r$ and 0 otherwise. 
For $r > k$, let $P_{r,k}(\alpha) := 0$. Note that $P_{r,k}(\alpha)$ has a zero of order at least $2r^2$ at $\alpha=r$, as proven in Proposition~\ref{cor:zero-order}. 
\end{definition}

The value of $P_{r,k}(\alpha)$ is independent of the choice of $\tau\in(0,\frac12)$, by contour deformation.
Throughout the paper, $c_j(r)$ and $\hat c_\ell(r)$ denote these fixed contours, with their dependence on $\tau$ suppressed. Firstly, we arrive at an explicit form for  $\mathcal{M}_{k}(\alpha)$ as defined in \eqref{eq:Mdef} in terms of $P_{r,k}(\alpha)$ when $k \geq 1$. For $k=1$, this matches the predicted form in \eqref{m1alpha}.

\begin{theorem}\label{mkalpha} Let $k \in \mathbb{N}$ and $\alpha \in \mathbb{R}_{+}$. Assume Conjecture~\ref{eq:main-conj} for $Y \ll T^{\gamma}$ with $\gamma = \max{(1,\alpha/2)}+\varepsilon$ for some fixed $\varepsilon > 0$. Then one has
\begin{align*}
    \mathcal{M}_{k}(\alpha) = \frac{1}{k!^{2}}{\sum_{ 0 \leq r < \alpha} } {k \choose r}^{2} P_{r,k}(\alpha),
\end{align*}
with $P_{r,k}(\alpha)$ as defined in \eqref{eq:prk} and in particular, $\mathcal{M}_{k}(\alpha)$ is continuously differentiable $2\ell^2 - 1$ times at $\alpha = \ell$ for $\ell = 1, \ldots, k$.
\end{theorem}

The polynomial range for the parameter $Y$, which bounds the imaginary components $\Im (z_1-z_2)$ and $\Im z_2$,  is needed for the sharp Perron truncation
in \eqref{eq:trunc-from-lemma312}.  The following theorem identifies
the variance coefficient $\mathcal D_k(\alpha)$, defined in \eqref{eq:Dvariance}, with $\delta_k(\alpha)$ under polylogarithmic
uniformity in $Y$. The derivative relation with $\mathcal M_k(\alpha)$
requires the stronger polynomial range for $Y$ assumed in
Theorem~\ref{mkalpha}.

\begin{theorem}
\label{conj:derivative}
Let $k \in \mathbb{N}$ and $\alpha \in \mathbb{R}_{+}$. Assume Conjecture~\ref{eq:main-conj} for $Y \ll (\log T)^{\eta}$ for some fixed $\eta>0$. Then one has
\begin{equation*}
    \mathcal{D}_{k}(\alpha) = \delta_k(\alpha).
\end{equation*}
Under the stronger hypothesis of Theorem~\ref{mkalpha}, we have $\mathcal{D}_{k}(\alpha) = \mathcal{M}_{k}'(\alpha)$.
\end{theorem}
We also prove an identity between $\mathcal{Q}_{k}(\alpha)$, i.e. the leading coefficient for the variance of von Mangoldt convolutions in arithmetic progressions defined in \eqref{def:W_k}, and the short-interval variance $\delta_{k}(\alpha)$ under Conjecture \ref{eq:ratios-arithmetic-corrected}.
\begin{theorem} 
\label{conj:arithmetic-variance}
Let $k \geq 4$ be an integer and $\alpha \in \mathbb{R}_{+}$. Assume Conjecture ~\ref{eq:ratios-arithmetic-corrected} for $Y \ll q^{\alpha/2 +\varepsilon}$ for some $\varepsilon > 0$. Then one has
\begin{equation*}
    \mathcal{Q}_{k}(\alpha) = \delta_{k}(\alpha).
\end{equation*}
\end{theorem}
A much stronger decay of the error term in Conjecture~\ref{eq:ratios-arithmetic-corrected} would extend the result to $k=1,2,3$; this is detailed in the proof of the above theorem. To verify that the same leading-order asymptotic arises in the case of the matrix integral in \eqref{I-tilde}, we prove the following result for $M_{k}(\alpha)$ which is defined in \eqref{eq:rmt-mkalpha}.
\begin{theorem}
\label{eq:rmt-theorem}
Let $k \in \mathbb{N}$ and $\alpha>0$. Then
\begin{equation*}
    M_{k}(\alpha) = \frac{1}{k!^{2}} \sum_{0\leq r < \alpha} {k \choose r}^2 P_{r,k}(\alpha),
\end{equation*}
with $P_{r,k}(\alpha)$ as defined in \eqref{eq:prk}.
\end{theorem} 
\noindent As established in \eqref{eq:derivative-relationship-intro}, $M_k(\alpha)$ is the antiderivative of $\delta_k(\alpha)$. Using Theorem~\ref{eq:rmt-theorem}, we may
replace the trace-sum matrix integral definition of $\delta_k(\alpha)$ by an explicit analytic
coefficient, allowing both $\mathcal D_k(\alpha)$ and $\mathcal Q_k(\alpha)$ to
be identified with the same function. In both proofs, we obtain explicit expressions that can be directly related to
\begin{equation}
\label{def:delta-derivative-Prk}
\delta_{k}(\alpha) =  \frac{1}{k!^{2}} \sum_{0\leq r < \alpha} {k \choose r}^2 P'_{r,k}(\alpha).
\end{equation}
\noindent To prove both Theorem ~\ref{mkalpha} and Theorem ~\ref{eq:rmt-theorem}, we will use the following three lemmas that provide the contour identities and degree estimates needed. Lemma~\ref{lem:two-centre} decomposes the multiple contour integral arising from Conjecture~\ref{eq:main-conj} according to the choice of pole centre for each variable, and shows that only terms with equal numbers of $u$- and $v$-variables at their respective shifted poles can contribute. Lemma~\ref{lem:P-to-Q-form} converts the symmetric shifted-moment expression into an unsymmetrised form, allowing direct comparison with the exact formula we use in the proof of Theorem~\ref{eq:rmt-theorem}. Lemma~\ref{lem:two-centre-degree-bound} shows that only finitely many Taylor coefficients of the holomorphic factor contribute to the residue, and provides the degree bounds needed to evaluate and estimate the resulting outer contour integrals.

\begin{lemma}
\label{lem:two-centre}
Fix $k \in \mathbb{Z}_{+}$, $L\in\mathbb C$ and $\sigma\in\mathbb C\setminus\{0\}$. Let $0<\tau<\frac12$, and choose $\mu>(1+\tau)|\sigma|$. Assume that $\Phi_L$ is holomorphic in a neighbourhood of the closed polydisc $\{u\in\mathbb C^k:\ |u_j|\le \mu,\ 1\le j\le k\}$ and is symmetric in
$u_1,\ldots,u_k$. Define
\begin{equation*}
A(u, v; \sigma) :=\sum_{j=1}^k\frac1{u_j}+\sum_{\ell=1}^k\frac1{v_\ell+\sigma} - L,\qquad
B(u, v; \sigma) :=\sum_{j=1}^k\frac1{u_j-\sigma}+\sum_{\ell=1}^k\frac1{v_\ell}.
\end{equation*}
Then we have
\begin{align*}
\int_{\substack{|u_j|=\mu\\|v_\ell|=\mu}}
\Phi_{L}(u_1, \ldots, u_k) \, \Delta(u)^2&\Delta(v)^2 \,
\prod_{j,\ell=1}^k (u_j-v_\ell-\sigma)
A(u, v; \sigma)^k B(u, v; \sigma)^k \prod_{j=1}^k \frac{du_j}{u_{j}^{k}}  \prod_{\ell=1}^k \frac{dv_\ell }{v_{\ell}^{k}}  \nonumber
 \\
&\qquad=
\sum_{r=0}^k\binom{k}{r}^2\, 
\sum_{a=0}^k
(-1)^{k-a}\binom{k}{a}\,
L^{k-a} \sigma^{k^2-(a+k)}\,
J_{L,r}^{(a)}(\sigma), 
\end{align*}
with
\begin{align*}
&J_{L,r}^{(a)}(\sigma) =
\int_{c_{j}(r)}  \int_{\hat{c}_{\ell}(r)}  
\Phi_{L}(\sigma x_1, \ldots, \sigma x_k) \,
\Delta(x)^2\Delta(y)^2\,
\prod_{j,\ell=1}^k (x_j-y_\ell-1)\,
\\
& \qquad \qquad \qquad \, \times
\Biggl(\sum_{j=1}^k\frac1{x_j}+\sum_{\ell=1}^k\frac1{y_\ell+1}\Biggr)^{a}
\Biggl(\sum_{j=1}^k\frac1{x_j-1}+\sum_{\ell=1}^k\frac1{y_\ell}\Biggr)^{k}
\prod_{j=1}^k x_j^{-k} dx_j\prod_{\ell=1}^k y_\ell^{-k} dy_\ell, \nonumber
\end{align*}
where $c_{j}(r)$ are counter-clockwise circles of radius $\tau$ centered at 1 for $j=1,\ldots, r$ and 0 otherwise; similarly, $\hat{c}_{\ell}(r)$ are counter-clockwise circles of radius $\tau$ centered at $-1$ for $\ell = 1, \ldots, r$ and 0 otherwise. 
\end{lemma}

\begin{lemma}
\label{lem:P-to-Q-form}
Let $k\geq 1$ and $z \in \mathbb C\setminus\{0\}$. Write
$U=\sum_{j=1}^k u_j$ and $V=\sum_{\ell=1}^k v_\ell$, and efine
\begin{equation*}
A_z(u,v)
=
\sum_{j=1}^k \frac{1}{u_j}
+
\sum_{\ell=1}^k \frac{1}{v_\ell-z},
\qquad
B_z(u,v)
=
\sum_{j=1}^k \frac{1}{u_j+z}
+
\sum_{\ell=1}^k \frac{1}{v_\ell}.
\end{equation*}
For finite sets $P,Q\subset\mathbb C$, define the residue
\begin{equation*}
\mathcal R_{P,Q}[F]
:=
\frac{1}{(2\pi i)^{2k}}
\int_{\Gamma_P^k}\int_{\Gamma_Q^k}
F(u,v)\,
\prod_{j=1}^k du_j
\prod_{\ell=1}^k dv_\ell,
\end{equation*}
where $\Gamma_P$ and $\Gamma_Q$ are positively oriented contours enclosing, respectively, the points of $P$ and $Q$, and no other singularities of $F$. Then, for all $s,t\in\mathbb C$,
\begin{align*}
\mathcal R_{\{0,-z\},\{0,z\}}
\Bigl[
&e^{sU+tV}
W_z(u,v)
\bigl(A_z(u,v)-s\bigr)^k
\bigl(B_z(u,v)-t\bigr)^k \Bigr]
\\ &=
\mathcal R_{\{0,-z\},\{0,z\}}
\left[
e^{(s-t)U}
W_z(u,v)
\bigl(A_z(u,v)-(s-t)\bigr)^k
B_z(u,v)^k
\right],
\end{align*}
where $W_z(u,v)=\Delta(u)^2\Delta(v)^2
\prod_{j,\ell=1}^k (u_j-v_\ell+z)
\prod_{j=1}^k u_j^{-k}
\prod_{\ell=1}^k v_\ell^{-k}.$
\end{lemma}
\begin{lemma}
\label{lem:two-centre-degree-bound}
Let $J_{L,r}^{(a)}(\sigma)$ be as defined in
Lemma~\ref{lem:two-centre}. Expand
$\Phi_L(\sigma x_1,\ldots,\sigma x_k)$ at the centres of the $x$-contours:
\begin{equation}
\label{eq:Phi-two-centre-Taylor}
\Phi_L(\sigma x_1,\ldots,\sigma x_k)
=
\sum_{\beta\in\mathbb Z_{\ge0}^{r}}
\sum_{\gamma\in\mathbb Z_{\ge0}^{k-r}}
\phi_{\beta,\gamma}^{(r)}(L,\sigma)
\prod_{j=1}^{r}(x_j-1)^{\beta_j}
\prod_{j=r+1}^{k}x_j^{\gamma_{j-r}} .
\end{equation}
If $a<r$, then $J_{L,r}^{(a)}(\sigma)=0$ and if $a\ge r$, then $J_{L,r}^{(a)}(\sigma)$ is unchanged if the above
Taylor expansion is truncated to the terms satisfying
\begin{equation*}
|\beta|\le k-r,
\qquad
|\gamma|\le a-r.
\end{equation*}
\end{lemma}
Here the empty multi-index convention is used for $\beta$ and $\gamma$. If we have no shifted variables, i.e. $r=0$, or no unshifted variables meaning that $r=k$, then the corresponding product from \eqref{eq:Phi-two-centre-Taylor} is interpreted as 1.  The two propositions that follow supply additional estimates and regularity properties. Under Conjecture~\ref{eq:main-conj}, Proposition~\ref{eq:final-corollary} gives the shifted-moment decay estimate used to control the tails in the short-interval variance argument.

\begin{proposition}
\label{eq:final-corollary}
Assume Conjecture~\ref{eq:main-conj}. Under the hypotheses
of Conjecture~\ref{eq:main-conj}, with $0<\xi<1$ if $k\ge3$ and
$0<\xi<c$ if $k=1,2$, one has
\begin{align}
\label{eq:first-corollary-statement}
&\frac{1}{T}\int_0^T
\left(\frac{\zeta'}{\zeta}\bigl(\tfrac12+it+z_1\bigr)\right)^k
\left(\frac{\zeta'}{\zeta}\bigl(\tfrac12-it-z_2\bigr)\right)^k dt
\\
&\qquad\qquad\qquad\qquad\qquad\qquad\qquad
\ll_{k,c,\xi}
(\log T)^{2k}
\bigl(1+\min(|z|,1)\log T\bigr)^{-2}
\left(1+\frac{1}{|z|\log T}\right)^{2k} \nonumber
\end{align}
where $z=z_1-z_2$. In particular,
\begin{align}
\label{eq:second-corollary-statement}
\frac{1}{T}\int_0^T
\Bigl(\frac{\zeta'}{\zeta}\bigl(\tfrac12+it+z_1\bigr)\Bigr)^k
\Bigl(\frac{\zeta'}{\zeta}\bigl(\tfrac12-it-z_2\bigr)\Bigr)^k dt
\ll
(\log T)^{2k}
\bigl(1+\min(|z|,1)\log T\bigr)^{-2}
\end{align}
uniformly on any subrange on which
$|z|\log T\ge c_0>0$. The implied constants may depend on $k$, $c$, $\xi$, and, in the final estimate, also on the lower bound $c_0$.
\end{proposition}
Throughout, the notation $A\ll B$ for complex-valued quantities means $|A|\le C |B|$ with $C>0$. We have a second proposition that determines the minimum zero-order for $P_{r,k}(\alpha)$ at $\alpha=r$ with $r$ an integer, enabling us to bound the differentiability of $\mathcal{M}_k(\alpha)$ at $\alpha =r$.
\begin{proposition}
\label{cor:zero-order}
For an integer $0\le r\le k$, $P_{r,k}(\alpha)$ has a zero of order at least $2r^2$ at $\alpha=r$.
\end{proposition}
\section{Proof of results}
\subsection{Proof  of Lemma~\ref{lem:two-centre}}
Fix $0<\tau<1/2$, $\sigma\neq 0$ and $\mu>(1+\tau)|\sigma|$. This ensures $\mu>|\sigma|$, so the contours $|u_j|=\mu$ and
$|v_\ell|=\mu$ enclose the explicit poles at $u_j=0,\sigma$ and
$v_\ell=0,-\sigma$.  We assume that $\Phi_L(u)$ is holomorphic in a neighbourhood of the closed
polydisc $\{|u_j|\le \mu,\ 1\le j\le k\}$. We consider a $2k$-fold contour integral of the form
\begin{align}
G_{L}(\sigma) \, := \, \int_{\substack{|u_j|=\mu \\|v_\ell|=\mu}}
&\Phi_{L}(u) \, \Delta(u)^2 \Delta(v)^2\,
\prod_{j,\ell=1}^k (u_j-v_\ell-\sigma) \prod_{j=1}^k u_j^{-k} \prod_{\ell=1}^k v_\ell^{-k} \\
& \biggl( \sum_{j=1}^k\frac1{u_j}+\sum_{\ell=1}^k\frac1{v_\ell+\sigma} - L \biggr)^k  \biggl( \sum_{j=1}^k\frac1{u_j-\sigma}+\sum_{\ell=1}^k\frac1{v_\ell} \biggr)^k  \prod_{j=1}^k du_j \prod_{\ell=1}^k dv_\ell. \nonumber
\end{align}
Let $\rho := \tau |\sigma|$. Each contour integral can be decomposed into a sum of two small circles of radius $\rho$. Every $u$-integral splits into an integral enclosing $0$ and another enclosing $\sigma$ while each $v$-integral splits into integrals enclosing either $0$ or $-\sigma$. Note that the condition $\mu>(1+\tau)|\sigma|$ guarantees that the small circles are contained within the original contours. Since $2\rho=2\tau|\sigma|<|\sigma|$, the small circles are pairwise disjoint and do not pass through any pole of the integrand. Also, as $\Phi_L(u)$ is holomorphic in the integration variables, no additional poles are
introduced. Therefore, by the residue theorem, we can write
\begin{align}
    G_{L}(\sigma) = \sum_{r=0}^{k} \sum_{s=0}^{k} {k \choose r} {k \choose s} G_{L,r,s}(\sigma)
\end{align}
where $r$ of the $u$-integrals are centered at $\sigma$, $s$ of the $v$-integrals are centered at $-\sigma$ and the remaining $(k-r)$ $u$-integrals and $(k-s)$ $v$-integrals are around the origin. Note that this decomposition requires that $\Phi_{L}$ is symmetric in the $u$-variables. Define
\begin{align}
\label{eq:N-contour-dependence}
G_{L,r,s}(\sigma) &= \int_{C_j^\sigma(r)} \int_{\hat C_\ell^\sigma(s)}
\Phi_{L}(u) \, \prod_{j,\ell = 1}^{k}(u_j-v_\ell  -\sigma) \prod_{i \neq j} (u_i - u_j)\prod_{i \neq j} (v_i - v_j)
\\
\nonumber &\qquad \qquad \quad\times
   A(u,v;\sigma)^{k} \, B(u,v;\sigma)^{k} \prod_j u_j  ^{-k} \prod_\ell v_\ell^{-k}
\prod_{j=1}^{k}du_j \, \prod_{\ell=1}^{k} dv_\ell,
\end{align}
with
\begin{equation}
A(u,v;\sigma) = \sum_{j=1}^{k}\frac{1}{u_j}
   \;+\;
   \sum_{\ell=1}^{k}\frac{1}{v_\ell + \sigma} - L, \qquad B(u,v;\sigma) = \sum_{\ell=1}^{k}\frac{1}{v_\ell}
   \;+\;
   \sum_{j=1}^{k}\frac{1}{u_j- \sigma},
\end{equation}
for $u=(u_1, \ldots, u_k)$ and $v=(v_1, \ldots, v_k)$. The $C_j^\sigma(r)$ are counter-clockwise circles of radius $\rho$, centered at $\sigma$ for $1\le j\le r$ and centered at $0$ for $r<j\le k$. Similarly, the $\hat C_\ell^\sigma(s)$ are counter-clockwise circles of radius
$\rho$, centered at $-\sigma$ for $1\le \ell\le s$ and centered at $0$ for $s<\ell\le k$. In order to extract all powers of $\sigma$, we expand
\begin{align}
    A(u,v;\sigma)^{k} =    \sum_{a=0}^{k}(-1)^{k-a}{k \choose a} L^{k-a} \biggl(\sum_{j=1}^{k}\frac{1}{u_j}
   \;+\;
   \sum_{\ell=1}^{k}\frac{1}{v_\ell + \sigma} \biggr)^{a}.
\end{align} 
 We perform a change of variables, $u_i \mapsto \sigma u_i$ and $v_\ell \mapsto \sigma v_\ell$, and then count the $\sigma$ contributions from each factor in \eqref{eq:N-contour-dependence}. For the factor $\prod_{i \neq j}(u_i - u_j)$, and its equivalent for the $v$-variables, we obtain an overall factor $\sigma^{2k(k-1)}$; from the mixed variable product, we extract $\sigma^{k^2}$. Each variable has an order-$k$ pole at zero; therefore the poles offer a factor of $\sigma^{-2k^2}$. With the contributions from  $A(u,v;\sigma)$ and $B(u,v;\sigma)$, we obtain an overall power
\begin{align}
    \sigma^{-(a+k) + k^2}.
\end{align}
After rescaling,
the contours have radius $\tau$; the new circles $c_{j}(r)$ are centered at $1$ for $j=1,\ldots, r$ and $0$ otherwise. The $\hat{c}_{\ell}(s)$ are centered at $-1$ for $\ell = 1, \ldots, s$ and $0$ otherwise. Thus
\begin{align}
    G_{L,r,s}(\sigma)
&= \sum^{k}_{a=0}(-1)^{k-a} {k \choose a}  L^{k-a} \, \sigma^{k^2-(a+k)} \, J_{L,r,s}^{(a)}(\sigma)
\end{align}
where 
\begin{align}
\label{eq:inner-termJ}
J_{L,r,s}^{(a)}(\sigma) &= \int_{c_{j}(r)} \int_{\hat{c}_{\ell}(s)}
\Phi_{L}(\sigma u ) 
\\
\nonumber &\quad\times
 \prod_{j,\ell}(u_j-v_\ell  -1) \prod_{i \neq j} (u_i - u_j)\prod_{i \neq j} (v_i - v_j) \prod_j u_j  ^{-k} \prod_\ell v_\ell^{-k}
\\ \nonumber &\quad\times \biggl(\sum_{j=1}^{k}\frac{1}{u_j}
   \;+\;
   \sum_{\ell=1}^{k}\frac{1}{v_\ell + 1} \biggr)^{a}
  \biggl(\sum_{\ell=1}^{k}\frac{1}{v_\ell}
   \;+\;
   \sum_{j=1}^{k}\frac{1}{u_j- 1} \biggr)^{k} \prod_{j=1}^{k}du_j \, \prod_{\ell=1}^{k} dv_\ell.
\end{align}
We now show that the contributions with $r\neq s$ vanish. Assume $r>s$ and fix $0\le a \le k$. As $0<\tau<1/2$, these contours enclose only their prescribed centers and no other point among $0,\pm1$.
We define 
\begin{equation}
\label{eq:def-S-and-T}
S(u,v) := \sum_{j=1}^{k}\frac{1}{u_j}
   \;+\;
   \sum_{\ell=1}^{k}\frac{1}{v_\ell + 1}, \qquad T(u,v) := \sum_{\ell=1}^{k}\frac{1}{v_\ell}
   \;+\;
   \sum_{j=1}^{k}\frac{1}{u_j- 1}.
\end{equation}
Thus the two bracketed terms in $J_{L,r,s}^{(a)}(\sigma)$ are $S(u,v)^a$ and
$T(u,v)^k$.
\paragraph{Shifted pole-count.} 
\label{subsubsec:shifted-pole-count}
For $j\leq r$, the only possible source of a pole at $u_j=1$ is $T(u,v)^k$, since all other factors are holomorphic there. Similarly, for $\ell\leq s$, the only possible source of a pole at $v_\ell=-1$ is $S(u,v)^a$. If $s>a$, then the iterated residue is immediately $0$, as there are not enough terms from $S(u,v)^a$ to supply a pole at each $v_1,\dots,v_s$. We may therefore assume $s\leq a$. Throughout, we refer to $u_1,\dots,u_r$ as the shifted $u$-variables, since their contours are centered at $1$, and to $v_{s+1},\dots,v_k$ as the unshifted $v$-variables, since their contours are centered at $0$. 
\\
\\
\noindent Fix a monomial term in the multinomial expansion of $T(u,v)^k$. Its dependence on the shifted $u$-variables and the unshifted $v$-variables can be written as
\begin{equation}
(u_1-1)^{-\alpha_1}\cdots (u_r-1)^{-\alpha_r}
v_{s+1}^{-\beta_{s+1}}\cdots v_k^{-\beta_k}
\times(\text{other factors}),
\end{equation}
where $\alpha_j\geq 0$ is the multiplicity with which $(u_j-1)^{-1}$ is selected and $\beta_\ell\geq 0$ is the multiplicity with which $v_\ell^{-1}$ is selected from $T(u,v)^k$. For the residue at $u_j=1$ to be non-zero, we must have
\begin{equation}
\alpha_j\geq 1
\qquad
(1\leq j\leq r),
\end{equation}
since otherwise the iterated residue across $u_1, \ldots, u_r$ is zero. Define
\begin{equation}
B:=\sum_{\ell=s+1}^{k}\beta_\ell, \qquad  E:=\sum_{j=1}^{r}(\alpha_j-1)\geq 0,
\qquad \text{so that} \quad
\sum_{j=1}^{r}\alpha_j=r+E.
\end{equation}
Since the chosen monomial in $T(u,v)^k$ contains exactly $k$ reciprocal factors in total, and $r+E$ of these are already used to supply the poles
\begin{equation}
(u_1-1)^{-\alpha_1}\cdots (u_r-1)^{-\alpha_r},
\end{equation}
we have the upper bound
\begin{equation}
\label{eq:B-upper}
B\leq k-r-E.
\end{equation}
Next we obtain a lower bound for \(B\) by considering the residues at the unshifted
\(v\)-variables \(v_{s+1},\dots,v_k\).  These residues are affected by the
powers of \(v_{s+1},\dots,v_k\) which remain in the numerator after the shifted
\(u\)-residues have been taken.  The relevant contribution comes from the mixed
product
\begin{equation}
\prod_{j,\ell=1}^{k}(u_j-v_\ell-1).
\end{equation}
In particular, consider the block involving the shifted $u$-variables and the unshifted $v$-variables
\begin{align}
\prod_{j=1}^{r}\prod_{\ell=s+1}^{k}(u_j-v_\ell-1)
&=
\prod_{j=1}^{r}\prod_{\ell=s+1}^{k}\bigl((u_j-1)-v_\ell\bigr).
\end{align}
For each fixed $j\leq r$, define
\begin{equation}
C_j(u_j;v_{s+1},\dots,v_k)
:=
\prod_{\ell=s+1}^{k}\bigl((u_j-1)-v_\ell\bigr).
\end{equation}
Expanding $C_j$ in powers of $u_j-1$, we obtain
\begin{equation}
\label{eq:two-centre-mixed-expansion}
C_j(u_j;v_{s+1},\dots,v_k)
=
\sum_{p=0}^{k-s}(u_j-1)^p Q_{p}(v_{s+1},\dots,v_k),
\end{equation}
where
\begin{equation}
Q_{p}(v_{s+1},\dots,v_k)
=
(-1)^{k-s-p}e_{k-s-p}(v_{s+1},\dots,v_k),
\end{equation}
and $e_{k-s-p}$ denotes the elementary symmetric polynomial of degree $k-s-p$ in $v_{s+1},\dots,v_k$.  Hence $Q_{p}$ is homogeneous of total $v$--degree
\begin{equation}
k-s-p.
\end{equation}
Equivalently, choosing $p$ copies of $u_j-1$ from the $k-s$ factors leaves $k-s-p$ copies of the unshifted $v$-variables. The pole $(u_j-1)^{-\alpha_j}$ requires, for a non-zero residue at $u_j=1$, the coefficient of $(u_j-1)^{\alpha_j-1}$ from the remaining integrand. Therefore, if the term
\begin{equation}
(u_j-1)^p Q_{p}(v_{s+1},\dots,v_k)
\end{equation}
contributes to the coefficient of $(u_j-1)^{\alpha_j-1}$, then necessarily
\begin{equation}
p\leq \alpha_j-1.
\end{equation}
Hence the total $v$--degree contributed by $Q_{p}$ in the variables $v_{s+1},\dots,v_k$ is
\begin{equation}
k-s-p
\geq
k-s-(\alpha_j-1).
\end{equation}
In other words, the $j$-th block starts with $k-s$ possible powers of the unshifted $v$-variables, and the $\alpha_j-1$ differentiations needed for the residue can eliminate at most $\alpha_j-1$ of them.  Thus the $j$-th block contributes lowest total unshifted $v$--degree at least
\begin{equation}
k-s-(\alpha_j-1).
\end{equation}
Since the total number of such differentiations is
\begin{equation}
\sum_{j=1}^{r}(\alpha_j-1)=E,
\end{equation}
after taking the $u_1,\dots,u_r$ residues the contribution of this block has lowest total $v$--degree in $v_{s+1},\dots,v_k$ at least
\begin{align}
\label{eq:deg-block}
\sum_{j=1}^{r}\bigl(k-s-(\alpha_j-1)\bigr)
&=
r(k-s)-\sum_{j=1}^{r}(\alpha_j-1)
\\
&=
r(k-s)-E.
\end{align}
This bound may be weak if $E>r(k-s)$, but it remains a valid lower bound; the actual lowest total degree in $v_{s+1},\dots,v_k$ is never negative. Moreover, the Vandermonde factor in the $v$-variables contains the subfactor involving only the unshifted variables, namely
\begin{equation}
\prod_{\substack{s<i\neq j\leq k}}(v_i-v_j).
\end{equation}
It contributes total $v$--degree $(k-s)(k-s-1)$ in the unshifted $v$-variables. Apart from the poles $v_{\ell}^{-k}$, the remaining factors have Taylor expansions in $v_{s+1},\dots,v_k$ with only non-negative powers. Hence none of these factors can decrease the lowest total $v$--degree. Therefore every term remaining in the numerator after taking the shifted residues at $u_1=\cdots=u_r=1$
has lowest total unshifted $v$--degree at least
\begin{equation}
\label{eq:deg-lower}
(k-s)(k-s-1)+r(k-s)-E.
\end{equation}
Now consider the poles at $v_{s+1},\dots,v_k=0$.  For $\ell>s$, the factor $v_\ell^{-k}$ is present in the integrand, and the chosen monomial from $T(u,v)^k$ contributes the additional pole $v_\ell^{-\beta_\ell}$, so the pole order at $v_\ell=0$ is $k+\beta_\ell$.  A non-zero residue at $v_\ell=0$ forces the numerator to contain a monomial of degree exactly $k+\beta_\ell-1$ in $v_\ell$; hence a non-zero residue at $v_{s+1}=\cdots=v_k=0$ requires a monomial of total $v$--degree
\begin{align}
\sum_{\ell=s+1}^{k}(k+\beta_\ell-1)
&=
(k-s)(k-1)+\sum_{\ell=s+1}^{k}\beta_\ell
\\
&=
(k-s)(k-1)+B.
\end{align}
Then a necessary condition for a non-zero contribution is
\begin{equation}
(k-s)(k-1)+B
\geq
(k-s)(k-s-1)+r(k-s)-E,
\end{equation}
which rearranges to
\begin{equation}
\label{eq:B-lower}
B\geq (k-s)(r-s)-E.
\end{equation}
Combining \eqref{eq:B-upper} and \eqref{eq:B-lower} gives
\begin{equation}
k-r-E\geq (k-s)(r-s)-E,
\end{equation}
and hence
\begin{equation}
k-r\geq (k-s)(r-s).
\end{equation}
Let $d:=r-s\geq 1$.  Then
\begin{align}
(k-s)(r-s)
&=
(k-r+d)d \nonumber
\\
&=
(k-r)d+d^2 \nonumber
\\
&>
k-r,
\end{align}
which is a contradiction. Hence the corresponding contribution is zero whenever $r>s$. The case $s>r$ is analogous, interchanging the roles of $u$ and $v$, and using $S(u,v)^a$ in place of $T(u,v)^k$. Thus the corresponding contribution is also zero when $s>r$, meaning the contribution vanishes whenever $r\neq s$. This means that $J_{L,r,s}^{(a)}(\sigma)=0$ and, consequently, $G_{L,r,s}(\sigma)=0$ whenever $r \neq s$. To conclude, we have
\begin{align}
G_{L}(\sigma)=
\sum_{r=0}^k\binom{k}{r}^2\, 
\sum_{a=0}^k
(-1)^{k-a}\binom{k}{a}\,
L^{k-a} \sigma^{k^2-(a+k)}\,
J_{L,r}^{(a)}(\sigma), 
\end{align}
where $J_{L,r}^{(a)}(\sigma) := J_{L,r, r}^{(a)}(\sigma)$ defined in \eqref{eq:inner-termJ} which is the claimed result.
\subsection{Proof of Lemma~\ref{lem:P-to-Q-form}}
If $P,Q\subset\mathbb C$ are finite sets and $F(u,v)$ is meromorphic, write
\begin{equation*}
\mathcal R_{P,Q}[F]
:=
\frac{1}{(2\pi i)^{2k}}
\int_{\Gamma_P^k}\int_{\Gamma_Q^k}
F(u,v)\,
\prod_{j=1}^k du_j
\prod_{\ell=1}^k dv_\ell,
\end{equation*}
where $\Gamma_P$ and $\Gamma_Q$ are positively oriented contours enclosing $P$ and $Q$, respectively, and no other poles of $F$. Equivalently, $\mathcal R_{P,Q}[F]$ is the sum of the $2k$-variable residues over all choices with $u_j$ at a point of $P$ for every $1 \leq j \leq k$, and $v_\ell$ at a point of Q for every $1 \leq \ell \leq k$. We use auxiliary variables to split the poles at $u_j=0,-z$ and $v_\ell=0,z$. Introduce two vectors of small parameters
\begin{equation*}
\mathbf{a} = (a_1,\ldots,a_k), \qquad \mathbf{b} = (b_1,\ldots,b_k),
\end{equation*}
and, for fixed $z\neq 0$, define
\begin{equation*}
\mathcal X:=\{a_1,\ldots,a_k,b_1-z,\ldots,b_k-z\},
\qquad
\mathcal X+z:=\{a_1+z,\ldots,a_k+z,b_1,\ldots,b_k\}.
\end{equation*}
We assume that the points of $\mathcal X$ are pairwise distinct, this assumption allows us to expand the $2k$-contour integral as a sum of simple local residues. Define the generating residue
\begin{align}
\label{eq:gen-fun}
\mathcal J_{s,t}(\mathbf a,\mathbf b;z)
:=
\mathcal   R_{\mathcal X,\mathcal X+z} 
\Bigg[
&e^{sU+tV}W_z(u,v) \\
&\times \prod_{r=1}^{k} \left[
e^{-s a_r}
\prod_{j=1}^{k} \left(1-\frac{a_r}{u_j}\right)^{-1}
\prod_{\ell=1}^{k} \left(1-\frac{a_r}{v_\ell-z}\right)^{-1}
\right] \notag \\
&\times \prod_{r=1}^{k} \left[
e^{-t b_r}
\prod_{j=1}^{k} \left(1-\frac{b_r}{u_j+z}\right)^{-1}
\prod_{\ell=1}^{k} \left(1-\frac{b_r}{v_\ell}\right)^{-1}
\right]
\Bigg] \notag
\end{align}
for $s, t \in \mathbb C$, where
\begin{equation}
 W_z(u,v)=\Delta(u)^2\Delta(v)^2
\prod_{j,\ell=1}^k (u_j-v_\ell+z)
\prod_{j=1}^k u_j^{-k}
\prod_{\ell=1}^k v_\ell^{-k}.   
\end{equation}
Differentiating \eqref{eq:gen-fun} once with respect to each $a_r$ and each $b_r$, and then setting $a_1 = \cdots = a_k = b_1 = \cdots = b_k = 0$, recovers the logarithmic derivative factors. With $A_z(u,v)$ and $B_z(u,v)$ defined as in Lemma~\ref{lem:P-to-Q-form}, we have
\begin{equation*}
\left.
\frac{\partial}{\partial a_r}
\left[
e^{-s a_r}
\prod_{j=1}^{k} \left(1-\frac{a_r}{u_j}\right)^{-1}
\prod_{\ell=1}^{k} \left(1-\frac{a_r}{v_\ell-z}\right)^{-1}
\right]
\right|_{a_r=0}
=
A_z(u,v)-s,
\end{equation*}
and similarly
\begin{equation*}
\left.
\frac{\partial}{\partial b_r}
\left[
e^{-t b_r}
\prod_{j=1}^{k} \left(1-\frac{b_r}{u_j+z}\right)^{-1}
\prod_{\ell=1}^{k} \left(1-\frac{b_r}{v_\ell}\right)^{-1}
\right]
\right|_{b_r=0}
=
B_z(u,v)-t.
\end{equation*}
Hence
\begin{align}
\label{eq:res-derivatives}
&\mathcal R_{\{0,-z\},\{0,z\}}
\left[
e^{sU+tV}
W_z(u,v)
\bigl(A_z(u,v)-s\bigr)^k
\bigl(B_z(u,v)-t\bigr)^k
\right] \notag \\
&\qquad \qquad \qquad \qquad =
\left.
\frac{\partial}{\partial a_1}\cdots\frac{\partial}{\partial a_k}
\frac{\partial}{\partial b_1}\cdots\frac{\partial}{\partial b_k}
\mathcal J_{s,t}(\mathbf a,\mathbf b;z)
\right|_{\mathbf a=\mathbf b=0}.
\end{align}

It remains to show that $\mathcal J_{s,t}$ depends on $s$ and $t$ only through the difference $s-t$. We examine the poles of the integrand defining $\mathcal J_{s,t}$. For each $u_j$, the factors involving $a_r$ and $b_r$ give
\begin{equation}
\label{eq:first-expansion-lemma-additive}
\prod_{r=1}^{k}
\left(1-\frac{a_r}{u_j}\right)^{-1}
\left(1-\frac{b_r}{u_j+z}\right)^{-1}
=
\frac{u_j^k(u_j+z)^k}
{\prod_{r=1}^{k}(u_j-a_r)\prod_{r=1}^{k}(u_j+z-b_r)} .
\end{equation}
As $W_z$ contains $u_j^{-k}$, the poles in $u_j$ are therefore at $u_j=a_r$ and $u_j=b_r-z$ for $1\leq r\leq k$. Similarly, for each $v_\ell$,
\begin{equation}
\label{eq:second-expansion-lemma-additive}
\prod_{r=1}^{k}
\left(1-\frac{a_r}{v_\ell-z}\right)^{-1}
\left(1-\frac{b_r}{v_\ell}\right)^{-1}
=
\frac{(v_\ell-z)^k v_\ell^k}
{\prod_{r=1}^{k}(v_\ell-z-a_r)\prod_{r=1}^{k}(v_\ell-b_r)} ,
\end{equation}
and since $W_z$ contains $v_\ell^{-k}$, the poles in $v_\ell$ are at $v_\ell=a_r+z$ and $v_\ell=b_r$ for $1\leq r\leq k$. Thus the $u$-poles are $u=x$ with $x\in\mathcal X$, while the corresponding $v$-poles are $v=x+z$ with $x\in\mathcal X$.

When expanding the $2k$-contour integral as a sum of local residues, each $u_j$ selects a pole $x\in\mathcal X$, and each $v_\ell$ selects a pole $x+z$ with $x\in\mathcal X$. These choices may have repetitions or may use the same base point $x$ for both a $u$-variable and a $v$-variable. For example, we have
\begin{equation*}
\mathcal X=\{a_1,a_2,b_1-z,b_2-z\}    
\end{equation*}
when $k=2$; one possible choice for the $u$-variables and $v$-variables is
\begin{equation*}
u_1=a_1,\qquad
u_2=b_2-z,\qquad
v_1=a_1+z,\qquad
v_2=b_2.
\end{equation*}
We will now show that such a paired choice does not, in fact, contribute to the $2k$-contour integral. More generally, a non-zero local residue is obtained only when the $u$-variables use a subset $S \subset\mathcal X$ of size \(k\), with no repetitions, and the $v$-variables use precisely the complementary subset $\mathcal X\setminus S$, shifted by $z$.

Suppose first that a candidate local residue contains both a $u$-variable at $x$ and a $v$-variable at $x+z$, for some $x\in\mathcal X$. To test whether this local residue is non-zero, introduce local coordinates $p$ and $q$ near these two poles: $p=0$ corresponds to $u=x$, and $q=0$ corresponds to $v=x+z$. We write
\begin{equation}
\label{eq:same-base}
u=x+p,\qquad v=x+z+q,
\end{equation}
for which the mixed factor
\begin{equation*}
\prod_{j,\ell=1}^{k}(u_j-v_\ell+z)
\end{equation*}
supplies
\begin{equation*}
u-v+z=p-q.
\end{equation*}
The corresponding simple pole factors from \eqref{eq:first-expansion-lemma-additive} and \eqref{eq:second-expansion-lemma-additive} give $1/(pq)$. As all other factors are holomorphic in $p$ and $q$ at $(0,0)$, we have
\begin{equation*}
\frac{p-q}{pq}
=
\frac{1}{q}-\frac{1}{p}.
\end{equation*}
This Laurent expansion has no $p^{-1}q^{-1}$ term, so the corresponding two-variable residue is zero. We also have that two $u$-variables cannot choose the same $x\in\mathcal X$. Indeed, if
\begin{equation*}
u_i=x+p,\qquad u_j=x+q,
\end{equation*}
then the local contribution from the Vandermonde and the two simple poles is
\begin{equation*}
\frac{(u_i-u_j)^2}{(u_i-x)(u_j-x)}
=
\frac{(p-q)^2}{pq}
=
\frac{p}{q}-2+\frac{q}{p}.
\end{equation*}
Again, after multiplication by factors holomorphic in $p$ and $q$, there can be no $p^{-1}q^{-1}$ coefficient arising from this local configuration. Hence this local residue is zero. The same argument applies to two $v$-variables choosing the same pole. It follows that every non-zero residue uses each $x\in\mathcal X$ exactly once: either as a $u$-pole $u=x$, or as a $v$-pole $v=x+z$, but not both. In other words, we choose a subset \(S\subset\mathcal X\) of size \(k\) for the \(u\)-variables, and choose the complementary subset \(\mathcal X\setminus S\), shifted by \(z\), for the \(v\)-variables. In such a contribution,
\begin{equation*}
U=\sum_{x\in S}x,
\qquad
V=\sum_{x\in\mathcal X\setminus S}(x+z),
\end{equation*}
as the $u$-variables are evaluated at the points of $S$ and the $v$-variables at the shifted points $x+z$ with $x\in\mathcal X\setminus S$. Since
\begin{equation*}
\sum_{x\in\mathcal X}x
=
\sum_{r=1}^{k}a_r+\sum_{r=1}^{k}(b_r-z)
=
\sum_{r=1}^{k}a_r+\sum_{r=1}^{k}b_r-kz,
\end{equation*}
and as $|\mathcal X\setminus S|=k$, we get
\begin{align}
\label{eq:V-form}
V
&=
\sum_{x\in\mathcal X\setminus S}x+kz \notag \\
&=
\left(\sum_{r=1}^{k}a_r+\sum_{r=1}^{k}b_r-kz-U\right)+kz \notag \\
&=
\sum_{r=1}^{k}a_r+\sum_{r=1}^{k}b_r-U .
\end{align}
Using the definition of $\mathcal J_{s,t}(\mathbf a,\mathbf b;z)$ in \eqref{eq:gen-fun}, consider the exponential part of each contributing local residue  
\begin{equation*}
e^{sU+tV}
e^{-s\sum_{r=1}^{k}a_r}
e^{-t\sum_{r=1}^{k}b_r}.
\end{equation*}
Substituting in the form of $V$ from \eqref{eq:V-form}, this becomes
\begin{equation*}
e^{sU+t(\sum a_r+\sum b_r-U)}
e^{-s\sum a_r}
e^{-t\sum b_r}
=
e^{(s-t)U}
e^{-(s-t)\sum a_r}.
\end{equation*}
This is exactly the exponential factor that appears in the corresponding local contribution to $\mathcal J_{s-t,0}(\mathbf a,\mathbf b;z)$. All other factors in the local residue are rational factors independent of $s$ and $t$. Hence each non-zero local residue contribution to $\mathcal J_{s,t}$ agrees with the corresponding contribution to $\mathcal J_{s-t,0}$, and summing these contributions gives
\begin{equation*}
\mathcal J_{s,t}(\mathbf a,\mathbf b;z)
=
\mathcal J_{s-t,0}(\mathbf a,\mathbf b;z).
\end{equation*}
Finally, taking derivatives of both sides as in \eqref{eq:res-derivatives} and then setting $\mathbf a=\mathbf b=0$, we obtain
\begin{align}
&\mathcal R_{\{0,-z\},\{0,z\}}
\left[
e^{sU+tV}
W_z(u,v)
\bigl(A_z(u,v)-s\bigr)^k
\bigl(B_z(u,v)-t\bigr)^k
\right] \notag \\
&\qquad\qquad =
\mathcal R_{\{0,-z\},\{0,z\}}
\left[
e^{(s-t)U}
W_z(u,v)
\bigl(A_z(u,v)-(s-t)\bigr)^k
B_z(u,v)^k
\right],
\end{align}
which proves the lemma.
\subsection{Proof of Lemma~\ref{lem:two-centre-degree-bound}}
We use similar
residue-degree bookkeeping as in the proof of
Lemma~\ref{lem:two-centre}. Write the two bracketed sums in the definition of $J_{L,r}^{(a)}(\sigma)$ as
\begin{equation}
\label{eq:consistency-check}
S(x,y):=
\sum_{j=1}^{k}\frac1{x_j}
+
\sum_{\ell=1}^{k}\frac1{y_\ell+1},
\qquad
T(x,y):=
\sum_{j=1}^{k}\frac1{x_j-1}
+
\sum_{\ell=1}^{k}\frac1{y_\ell}.
\end{equation}
First observe that, for $\ell\le r$, the only possible pole at
$y_\ell=-1$ comes from the factor $(y_\ell+1)^{-1}$ in $S(x,y)^a$.
Therefore a non-zero residue requires one factor $(y_\ell+1)^{-1}$ for each $\ell=1,\ldots,r$.  Hence, if $a<r$, then
$J_{L,r}^{(a)}(\sigma)=0$.

Assume now that $a\ge r$.  Consider a fixed monomial in the expansion of
$S(x,y)^a$.  Let $D$ be the total number of factors $x_j^{-1}$ selected
from the unshifted variables $x_{r+1},\ldots,x_k$, and let $H$ be the total
number of factors $(y_\ell+1)^{-1}$ selected from the shifted variables
$y_1,\ldots,y_r$.  Then
\begin{equation*}
D+H\le a,
\qquad
H\ge r
\end{equation*}
for a non-zero contribution. Now consider a fixed monomial in the expansion of $T(x,y)^k$.  For
$j\le r$, the only possible pole at $x_j=1$ comes from factors
$(x_j-1)^{-1}$ appearing in this monomial.  Suppose that the monomial
contains $(x_j-1)^{-\alpha_j}$, and set
\begin{equation*}
A:=\sum_{j=1}^{r}\alpha_j .
\end{equation*}
As this monomial has total degree $k$, we have $A\le k$. For a non-zero contribution, we
may assume $\alpha_j\ge1$ for every $j\le r$.  To take the residue at
$x_j=1$, the remaining factors which are holomorphic in $x_1,\ldots,x_r$ near $x_j=1$ must supply the coefficient of
$(x_j-1)^{\alpha_j-1}$.  Therefore the Taylor
expansion of $\Phi_L(\sigma x_1,\ldots,\sigma x_k)$ can contribute total
degree in $x_1-1,\ldots,x_r-1$ at most
\begin{equation*}
\label{eq:degree-bound-for-shifted-variables}
\sum_{j=1}^{r}(\alpha_j-1)=A-r\le k-r.
\end{equation*}

It remains to bound the possible degree coming from the unshifted variables
$x_{r+1},\ldots,x_k$.  For $j>r$, the factor $x_j^{-k}$ gives a pole of
order $k$ at $x_j=0$, and the chosen monomial from $S(x,y)^a$ gives an
additional pole of total order $D$ across the unshifted $x$-variables.
Hence the residues at $x_{r+1}=\cdots=x_k=0$ require the remaining holomorphic factors to supply total degree
\begin{equation*}
(k-r)(k-1)+D
\end{equation*}
in the variables $x_{r+1},\ldots,x_k$.

For the mixed product, use the same coefficient-wise expansion as in the proof
of Lemma~\ref{lem:two-centre}.  For fixed $\ell\le r$,
\begin{equation*}
\prod_{j=r+1}^{k}\bigl(x_j-(y_\ell+1)\bigr)
=
\sum_{p=0}^{k-r}
(-1)^p (y_\ell+1)^p
e_{k-r-p}(x_{r+1},\ldots,x_k).
\end{equation*}
If $\nu_\ell$ denotes the multiplicity of $(y_\ell+1)^{-1}$ selected from
$S(x,y)^a$, then taking the residue at $y_\ell=-1$ can lower the degree contributed by the mixed product in the unshifted $x$-variables by at most
$\nu_\ell-1$. Since
$H=\sum_{\ell=1}^{r}\nu_\ell$, the mixed block contributes degree at least
\begin{equation*}
r(k-r)-\sum_{\ell=1}^{r}(\nu_\ell-1)
=
r(k-r)-(H-r)
\end{equation*}
in the unshifted $x$-variables.  The Vandermonde factor involving only
$x_{r+1},\ldots,x_k$ contributes total degree $(k-r)(k-r-1)$. Thus the
mixed product and Vandermonde factor already give degree at least
\begin{equation*}
(k-r)(k-r-1)+r(k-r)-(H-r)
=
(k-r)(k-1)-(H-r)
\end{equation*}
in $x_{r+1},\ldots,x_k$. Consequently, the Taylor expansion of
$\Phi_L(\sigma x_1,\ldots,\sigma x_k)$ can contribute total degree in
$x_{r+1},\ldots,x_k$ at most
\begin{equation}
\label{eq:degree-bound-for-unshifted-variables}
\bigl((k-r)(k-1)+D\bigr)
-
\bigl((k-r)(k-1)-(H-r)\bigr)
=
D+H-r
\le a-r.
\end{equation}
We use multi-indices
\begin{equation*}
    \beta=(\beta_1,\ldots,\beta_r)\in\mathbb Z_{\ge0}^{r},
    \qquad
    \gamma=(\gamma_1,\ldots,\gamma_{k-r})\in\mathbb Z_{\ge0}^{k-r},
\end{equation*}
where $\beta$ records the powers of the shifted variables
$(x_1-1),\ldots,(x_r-1)$ and $\gamma$ records the powers of the unshifted
variables $x_{r+1},\ldots,x_k$. We write
\begin{equation}
\label{eq:beta-defs}
    |\beta|=\sum_{j=1}^{r}\beta_j,\qquad
    |\gamma|=\sum_{m=1}^{k-r}\gamma_m,\qquad
    \beta!=\prod_{j=1}^{r}\beta_j!,\qquad
    \gamma!=\prod_{m=1}^{k-r}\gamma_m!,
\end{equation}
with the usual convention that empty products equal $1$. As $\Phi_L(\sigma x_1,\ldots,\sigma x_k)$ is holomorphic in a neighbourhood of the product of the $x$-contours, its Taylor expansion about the corresponding contour centres converges absolutely and uniformly on these contours. Hence the Taylor expansion may be inserted into the contour integral term by term, and the resulting residue may be computed coefficient-wise. Together with \eqref{eq:degree-bound-for-shifted-variables} and \eqref{eq:degree-bound-for-unshifted-variables},
inserting the Taylor expansion
\begin{equation*}
\Phi_L(\sigma x_1,\ldots,\sigma x_k)
=
\sum_{\beta\in\mathbb Z_{\ge0}^{r}}
\sum_{\gamma\in\mathbb Z_{\ge0}^{k-r}}
\phi_{\beta,\gamma}^{(r)}(L,\sigma)
\prod_{j=1}^{r}(x_j-1)^{\beta_j}
\prod_{j=r+1}^{k}x_j^{\gamma_{j-r}},
\end{equation*}
into the definition of $J_{L,r}^{(a)}(\sigma)$ proves that every term with
$|\beta|>k-r$ or $|\gamma|>a-r$ gives zero. By linearity of the residue,
\begin{equation*}
J_{L,r}^{(a)}(\sigma)
=
\sum_{\substack{\beta\in\mathbb Z_{\ge0}^{r}\\ |\beta|\le k-r}}
\sum_{\substack{\gamma\in\mathbb Z_{\ge0}^{k-r}\\ |\gamma|\le a-r}}
\phi_{\beta,\gamma}^{(r)}(L,\sigma)
\, \mathcal C_{\beta,\gamma}^{(r,a)},
\end{equation*}
where the constants $\mathcal C_{\beta,\gamma}^{(r,a)}$ are independent of $L$ and $\sigma$.
\subsection{Proof of Theorem \ref{eq:rmt-theorem}}
\label{sec:Mk-derivation}
Recall that $\tilde{I_{k}}(n;N)$ denotes the truncated sum over traces of total degree $\leq n$, whereas $I_{k}(n;N)$ corresponds to a total degree of exactly $n$ as defined in \eqref{I-tilde} and \eqref{eq:Ik}. To obtain the leading asymptotic for $\tilde{I_{k}}(n;N)$, we begin by noting that $I_{k}(n;N)$ can be expressed as a coefficient extraction for the moments of logarithmic derivatives of characteristic polynomials evaluated within the unit circle. 

\noindent We have
\begin{align}
    \frac{1}{2\pi i} \int_{|z| = \eta} \biggl( z \frac{\Lambda'_{U}}{\Lambda_{U}}(z)\biggr)^{k}&\frac{1}{z^{n+1}}dz = \frac{1}{2 \pi i} \int_{|z|= \eta}  \biggl( -\sum_{j=1}^{\infty} \text{Tr}(U^{*j})z^{j}\biggr)^{k}\frac{1}{z^{n+1}}dz \\ \nonumber
    &= (-1)^{k} \sum_{m=1}^{\infty} \sum_{\substack{j_{1}, \ldots, j_{k} \geq 1 \\ j_{1} + \cdots + j_{k} = m}} \frac{1}{2 \pi i} \int_{|z|= \eta}   \text{Tr}(U^{*j_{1}}) \cdots \text{Tr}(U^{*j_{k}}) z^{m-(n+1)}dz \\ \nonumber
    &= (-1)^{k} \sum_{\substack{j_{1}, \ldots, j_{k} \geq 1 \\ j_{1} + \cdots + j_{k} = n}} \text{Tr}(U^{*j_{1}}) \cdots \text{Tr}(U^{*j_{k}}), 
\end{align}
where $|\eta| < 1$. One may write
\begin{align*}
    I_{k}(n;N) &:= \frac{1}{(2\pi i)^{2}} \int_{U(N)} \int_{|z_{1}| = \eta}\int_{|z_{2}| = \eta} \biggl( z_1 \frac{\Lambda'_{U}}{\Lambda_{U}}(z_1)\biggr)^{k}\biggl( z_2 \frac{\Lambda'_{U^{*}}}{\Lambda_{U^{*}}}(z_2)\biggr)^{k} \frac{1}{(z_{1}z_{2})^{n+1}}dz_{1}dz_{2} dU \\ \nonumber
    &= \sum_{\substack{j_{1}, \ldots, j_{k} \geq 1 \\ j_{1} + \cdots + j_{k} = n}}  \sum_{\substack{\ell_{1}, \ldots, \ell_{k} \geq 1 \\ \ell_{1} + \cdots + \ell_{k} = n}} \int_{U(N)} \text{Tr}(U^{j_{1}}) \cdots \text{Tr}(U^{j_{k}})    \overline{\text{Tr}(U^{\ell_{1}})} \cdots \overline{\text{Tr}(U^{\ell_{k}})}   dU \\ \nonumber 
    &= \int_{U(N)} \biggl| \sum_{\substack{j_{1}, \ldots, j_{k} \geq 1 \\ j_{1} + \cdots + j_{k} = n}} \text{Tr}(U^{j_{1}}) \cdots \text{Tr}(U^{j_{k}})  \biggr|^{2}dU.
\end{align*}
In the above, we use the phase invariance of Haar measure as shown in obtaining \eqref{eq:Ik}. Due to \eqref{line:haar-invariance}, $I_{k}(n;N)$ can be expressed as a single contour
\begin{align}
     I_{k}(n;N) = \frac{1}{(2\pi i)^{2}} \int_{|z_{1}| = \eta}\int_{|z_{2}| = \eta} \frac{F_{k}(z_{1}z_{2}; N)}{(z_{1}z_{2})^{n+1}}dz_{1}dz_{2} = \frac{1}{2\pi i} \int_{|s| = \eta^{2}} \frac{F_{k}(s; N)}{s^{n+1}}ds,
\end{align}
where 
\begin{equation} \label{eq:fn}
    F_{k}(z_{1}z_{2}; N) = \int_{U(N)}\biggl( z_1 \frac{\Lambda'_{U}}{\Lambda_{U}}(z_1)\biggr)^{k}\biggl( z_2 \frac{\Lambda'_{U^{*}}}{\Lambda_{U^{*}}}(z_2)\biggr)^{k}dU .
\end{equation}
We may also express $\tilde{I_{k}}(n;N)$ in terms of $F_{k}(s;N)$ as
\begin{equation}
\label{eq:outer-integral-for-cauchy}
    \tilde{I_{k}}(n;N) = \sum_{j=1}^{n} I_{k}(j,N) = -\frac{1}{2\pi i} \int_{|\xi| = \eta^{2}} \frac{F_{k}(\xi;N)}{\xi^{n+1}(\xi-1)}d\xi.
\end{equation}
An expression for $F_k(\xi;N)$, when $|\xi|<1$, is well-known in the literature. It was originally written using set notation (see Theorem 3.3 of \cite{kn:consna08}), but was rewritten as a multiple contour integral in equation (43) of \cite{kn:bbbcprs19}. In either case, it is derived by differentiating the Ratios Theorem (Theorem 4.1 in \cite{kn:cfz2}), which we restate in Theorem~\ref{eq:permutation-ratios}. To get an expression in a form that is useful to us, we will repeat the steps of that process here.
\begin{theorem} [\cite{kn:cfz2}, Theorem 4.1] 
\label{eq:permutation-ratios}
Let $k,N\geq 1$ be integers, and let
\begin{equation*}
\alpha_1,\ldots,\alpha_{2k},\gamma_1,\ldots,\gamma_k,\delta_1,\ldots,\delta_k\in\mathbb C,
\qquad
\Re\gamma_q>0,\quad \Re\delta_r>0
\end{equation*}
for \(1\le q,r\le k\). Then
 \begin{align}
 \label{eq:one-to-reference-for-sign}
 & \int_{U(N)}\frac{\prod_{j=1}^k\Lambda_U(e^{-\alpha_j})\prod_{l=k+1}^{2k}\Lambda_{U^*}(e^{\alpha_{l}})}{\prod_{q=1}^k\Lambda_U(e^{-\gamma_{q}})\prod_{r=1}^k\Lambda_{U^*}(e^{-\delta_{r}})}dU_N \nonumber \\
 & \quad = \sum_{\sigma \in \Xi_{k,k} } e^{N \sum_{p=1}^k (\alpha_{\sigma(p)} -
\alpha_p) }
{\prod_{p=1}^k \prod_{l=k+1}^{2k} z(\alpha_{\sigma(p)} -
\alpha_{\sigma(l)} ) \prod_{q=1}^k \prod_{r=1}^k z(\gamma_q+\delta_r)
\over
\prod_{r=1}^k \prod_{p=1}^k z(\alpha_{\sigma(p)} + \delta_r)
\prod_{q=1}^k \prod_{l=k+1}^{2k} z(\gamma_q - \alpha_{\sigma(l)}) },
\end{align}
where $z(x) := \frac{1}{1-e^{-x}}$ and $\Xi_{k,k}$ denotes the set of permutations $\sigma \in S_{2k}$ such that $1 \le \sigma(1) < \cdots <  \sigma(k) \le 2k $ and $1 \le \sigma(k+1) < \cdots <  \sigma(2k) \le 2k$.
\end{theorem}
\noindent Next we let $a_j=e^{\alpha_j}$ for $1 \leq j \leq 2k $, $c_q=e^{\gamma_q}$ and $d_r=e^{-\delta_r}$ for $1\leq q,r \leq k$. Then for $|d_r|<1$ and $|c_q|>1$, we have 
\begin{align}
\label{eq:permutation-res}
 & \int_{U(N)}\frac{\prod_{j=1}^k\Lambda_U(\frac{1}{a_j})\prod_{l=k+1}^{2k}\Lambda_{U^*}(a_l)}{\prod_{q=1}^k\Lambda_U(\frac{1}{c_q})\prod_{r=1}^k\Lambda_{U^*}(d_r)}dU_N \\
 & \nonumber \quad  =(-1)^{k} \prod_{j=1}^ka_j^{-N} \sum_{\sigma \in \Xi_{k,k} }\left(\prod_{j=1}^k a_{\sigma(j)}^{N}\right){\prod_{r=1}^k \prod_{p=1}^{k} (1-\frac{a_{\sigma(p)}}{d_r} ) \prod_{q=1}^k \prod_{l=k+1}^{2k} (a_{\sigma(l)}-c_q)
\over
\prod_{p=1}^k \prod_{l=k+1}^{2k} (a_{\sigma(p)} - a_{\sigma(l)})
\prod_{q=1}^k \prod_{r=1}^{k} (1 -\frac{c_q}{d_r}) }.
\end{align}
In the above, we obtain an overall  $(-1)^k=(-1)^{k^2}$ factor after rearranging the expressions resulting from the $z(x)$ functions in \eqref{eq:one-to-reference-for-sign}. We then apply Lemma 2.5.1 from \cite{kn:cfkrs} to express \eqref{eq:permutation-res} as
\begin{align}
\label{eq:permutation-contour}
 & \frac{\prod_{j=1}^{k}a_j^{-N}}{(k!)^2 (2\pi i)^{2k}} \, \,  \oint \cdots \oint \prod_{j=1}^k w_j^{N}\frac{\prod_{r=1}^k\prod_{p=1}^{k}(1-\frac{w_p}{d_r})\prod_{q=1}^k\prod_{l=k+1}^{2k}(w_l-c_q)}{\prod_{p=1}^k\prod_{l=k+1}^{2k}(w_p-w_l)\prod_{q=1}^k\prod_{r=1}^k(1-\frac{c_q}{d_r})}\nonumber\\
 &\qquad\qquad\qquad\qquad\qquad\qquad\qquad\qquad \qquad\qquad \, \times\frac{\prod_{1\leq i <j\leq 2k}(w_j-w_i)^2\prod_{j=1}^{2k}dw_j}{\prod_{p=1}^{2k}\prod_{j=1}^{2k}(w_p-a_j)}.
\end{align}
Note that each of the $w$-contours encloses every $a_j$ for $1 \leq j \leq 2k$. In a similar manner to \cite{kn:bbbcprs19}, we proceed by computing a series of derivatives
\begin{equation}
\left.\frac{d}{da_1}\cdots \frac{d}{da_{2k}}\int_{U(N)}\frac{\prod_{j=1}^k\Lambda_U(\frac{1}{a_j})\prod_{l=k+1}^{2k}\Lambda_{U^*}(a_l)}{\prod_{q=1}^k\Lambda_U(\frac{1}{c_q})\prod_{r=1}^k\Lambda_{U^*}(d_r)}dU_N\right|_{ {a_1,\ldots,a_k=a \atop a_{k+1},\ldots,a_{2k}=b}\atop c_q=a, d_r=b }.
\end{equation}
\noindent On the one hand, this equals
\begin{equation}
\frac{(-1)^k}{a^{2k}}\int_{U(N)} \frac{\Lambda_U'(1/a)^k\Lambda_{U^{*}}'(b)^k}{\Lambda_U(1/a)^k\Lambda_{U^{*}}(b)^k} dU_N
\end{equation}
with the $(-1)^{k}/a^{2k}$ factor from applying the chain rule. And on the other side, from \eqref{eq:permutation-contour}, it equals
\begin{align}
&\frac{d}{da_1}\cdots \frac{d}{da_{2k}} \biggl( \frac{\prod_{j=1}^{k}a_j^{-N}}{(k!)^2 (2\pi i)^{2k}}\nonumber\\
 &\qquad \qquad \qquad \times \oint \cdots \oint \prod_{j=1}^k w_j^{N}\frac{\prod_{r=1}^k\prod_{p=1}^{k}(1-\frac{w_p}{d_r})\prod_{q=1}^k\prod_{l=k+1}^{2k}(w_l-c_q)}{\prod_{p=1}^k\prod_{l=k+1}^{2k}(w_p-w_l)\prod_{q=1}^k\prod_{r=1}^k(1-\frac{c_q}{d_r})}\nonumber\\
 &\left.\quad\qquad\qquad\qquad\qquad\qquad\qquad\quad \times\frac{\prod_{1\leq i <j\leq 2k}(w_j-w_i)^2\prod_{j=1}^{2k}dw_j}{\prod_{p=1}^{2k}\prod_{j=1}^{2k}(w_p-a_j)}\right|_{ {a_1, \ldots,a_k=a \atop a_{k+1},\ldots,a_{2k}=b}\atop c_q=a, d_r=b } \biggr) .
\end{align}
For $\varepsilon \in \{0,1\}$, we have
\begin{equation}
\label{eq:derivative-rule}
\frac{d}{da}\left( \frac{a^{-\varepsilon N}}{\prod_{j=1}^{2k}(w_j-a)}\right)=\frac{a^{-\varepsilon N}}{\prod_{j=1}^{2k}(w_j-a)} \left( \frac{-\varepsilon N}{a} +\sum_{j=1}^{2k} \frac{1}{w_j-a}\right).\end{equation}
For each of the first $k$ variables $a_1,\ldots,a_k$, the
differentiated factor in the contour integral corresponds to the $\varepsilon = 1$ case in \eqref{eq:derivative-rule} and the remaining variables are the $\varepsilon = 0$ case. Therefore 
\begin{align}
\label{eq:log-der}
\frac{(-1)^k}{a^{2k}}\int_{U(N)} \frac{\Lambda_U'(1/a)^k\Lambda_{U^{*}}'(b)^k}{\Lambda_U(1/a)^k\Lambda_{U^{*}}(b)^k} &dU_N \nonumber \\
= \frac{1}{(k!)^2 (2\pi i )^{2k}}\oint\cdots \oint &\prod_{j=1}^k  w_j^{N}\left[ \frac{a^{ -N}}{\prod_{j=1}^{2k}(w_j-a)} \left( -\frac{N}{a} +\sum_{j=1}^{2k} \frac{1}{w_j-a}\right)\right]^k\nonumber \\
& \times \left[ \frac{1}{\prod_{j=1}^{2k}(w_j-b)} \left( \sum_{j=1}^{2k} \frac{1}{w_j-b}\right)\right]^k \prod_{1\leq i <j \leq 2k} (w_j-w_i)^2\nonumber \\
&\times \frac{ \prod_{j=1}^k(1-\tfrac{w_j}{b})^k \prod_{l=k+1}^{2k}(w_l-a)^k}{\left(\prod_{j=1}^k \prod_{l=k+1}^{2k}(w_j-w_l)\right) \;(1-\tfrac{a}{b})^{k^2}} \prod_{j=1}^{2k} dw_j.
\end{align}
The variables \(c_q\) and \(d_r\) are held fixed
during the differentiation and are only afterwards specialised to \(c_q=a\) and
\(d_r=b\). This means we inherit $|a|>1$ and $|b| < 1$ from the conditions $|c_q| > 1$ and $|d_r| < 1$, which ensures that the logarithmic derivatives are evaluated away from the poles on the unit circle. Proceeding from \eqref{eq:log-der}, we multiply both the numerator and denominator by $b^{k^2}$ to turn the $(1-w_j/b)$ factors into $b-w_j$ and the $(1-a/b)$  factors into $(b-a)$. We have
\begin{align}
\label{eq:unsymmetrised-form}
&\frac{1}{a^{2k}}\int_{U(N)} \frac{\Lambda_U'(1/a)^k\Lambda_{U^{*}}'(b)^k}{\Lambda_U(1/a)^k\Lambda_{U^{*}}(b)^k} dU_N\nonumber \\
&= \frac{a^{-kN}}{(b-a)^{k^2}(k!)^2 (2\pi i)^{2k}}\oint\cdots \oint \frac{\prod_{j=1}^k w_j^{N} \prod_{j=1}^k(w_j-b)^k \;\prod_{l=k+1}^{2k} (w_l-a)^k}{\left(\prod_{j=1}^k \prod_{l=k+1}^{2k}(w_j-w_l)\right)\prod_{j=1}^{2k}(w_j-a)^k(w_j-b)^k } \nonumber \\
&\quad \qquad \qquad \times \left(\sum_{j=1}^{2k} \frac{1}{w_j-a}-\frac{N}{a}\right)^k\; \left( \sum_{j=1}^{2k}\frac{1}{w_j-b}\right)^k \prod_{1\leq i <j \leq 2k} (w_j-w_i)^2\prod_{j=1}^{2k} dw_j.
\end{align}
Note that the $w$-contours are centered at the origin and have a radius greater than $\max(|a|,|b|) = |a|$. Next we change variables so that $u_j=w_j-a$ for $j=1,\ldots,k$ and $v_j=w_{k+j}-b$ for $j=1,\ldots k$.  The new contours of integration of the $u$ variables are circles centered at $-a$ and the contours of the $v$ variables are centered at $-b$; the radius remains unchanged. We have
\begin{align}
&\frac{1}{a^{2k}}\int_{U(N)} \frac{\Lambda_U'(1/a)^k\Lambda_{U^{*}}'(b)^k}{\Lambda_U(1/a)^k\Lambda_{U^{*}}(b)^k} dU_N  \\
&=\frac{a^{-kN}}{(b-a)^{k^2}(k!)^2 (2\pi i)^{2k}}\oint\cdots \oint \frac{\prod_{j=1}^k (u_j+a)^{N}}{\left(\prod_{j=1}^k \prod_{l=1}^{k}(u_j+a-v_l-b)\right)}\frac{1}{\prod_{j=1}^k u_j^kv_j^k}\nonumber \\
&\qquad\qquad\qquad\qquad\qquad \times \left(\sum_{j=1}^{k} \frac{1}{u_j}+\sum_{j=1}^k \frac{1}{v_j+b-a}-\frac{N}{a}\right)^k\; \left( \sum_{j=1}^{k}\frac{1}{v_j}+\sum_{j=1}^k \frac{1}{u_j+a-b} \right)^k\nonumber \\
&\qquad\qquad\qquad\qquad\qquad \times \Delta(u)^2 \Delta(v)^2 \prod_{j=1}^k\prod_{l=1}^k (v_l+b-u_j-a)^2 \prod_{j=1}^{k} du_jdv_j, \nonumber
\end{align}
where $\Delta(u) = \prod_{i < j} (u_j - u_i)$. A single power of the $(v_{\ell} + b - u_j -a)^2$ factors cancels against the $(u_j+a-v_l-b)$ poles which introduces a net sign of $(-1)^k$. After this cancellation, there are no remaining mixed poles of the form
\(u_j+a-v_l-b=0\). Thus the possible poles in each \(u_j\)
are at \(0\) and \(b-a\) while the possible poles in each \(v_j\) are at
\(0\) and \(a-b\). Therefore, we arrive at the simplified form 
\begin{align}
&\frac{1}{a^{2k}}\int_{U(N)} \frac{\Lambda_U'(1/a)^k\Lambda_{U^{*}}'(b)^k}{\Lambda_U(1/a)^k\Lambda_{U^{*}}(b)^k} dU_N \\
& \qquad\qquad\qquad\qquad\qquad =\frac{(-1)^k \, a^{-kN}}{(b-a)^{k^2}(k!)^2 (2\pi i)^{2k}}\oint\cdots \oint \prod_{j=1}^k (u_j+a)^{N} \Delta(u)^2 \Delta(v)^2\nonumber \\
&\qquad\qquad\qquad\qquad\qquad\quad \times \prod_{j=1}^k\prod_{l=1}^k (v_l-u_j+b-a) \left(\sum_{j=1}^{k} \frac{1}{u_j}+\sum_{j=1}^k \frac{1}{v_j+b-a}-\frac{N}{a}\right)^k \nonumber \\
&\qquad\qquad\qquad\qquad\qquad\quad \times \left( \sum_{j=1}^{k}\frac{1}{v_j}+\sum_{j=1}^k \frac{1}{u_j+a-b} \right)^k \frac{1}{\prod_{j=1}^k u_j^kv_j^k}\prod_{j=1}^{k} du_jdv_j. \nonumber
\end{align} 
We now scale every $u$ and $v$ variable by $a$ to recover a form for the moments above that depend only on the single variable $\xi = \tfrac b a$. The scaling produces an $a^{-2k}$ factor that cancels the prefactor on the left-hand side. We have
\begin{align}
&\int_{U(N)} \frac{\Lambda_U'(1/a)^k\Lambda_{U^{*}}'(b)^k}{\Lambda_U(1/a)^k\Lambda_{U^{*}}(b)^k} dU_N  \\
&\quad\quad\qquad\qquad\qquad =\frac{(-1)^k}{(\xi-1)^{k^2}(k!)^2 (2\pi i)^{2k}}\oint\cdots \oint \prod_{j=1}^k (u_j+1)^{N} \Delta(u)^2 \Delta(v)^2\nonumber \\
&\quad\qquad\qquad\qquad\quad\quad \, \times \prod_{j=1}^k\prod_{l=1}^k (v_l-u_j+\xi-1) \left(\sum_{j=1}^{k} \frac{1}{u_j}+\sum_{j=1}^k \frac{1}{v_j+\xi-1}-N\right)^k \nonumber \\
&\quad\qquad\qquad\qquad\quad\quad \, \times \left( \sum_{j=1}^{k}\frac{1}{v_j}+\sum_{j=1}^k \frac{1}{u_j+1-\xi} \right)^k \frac{1}{\prod_{j=1}^k (u_j v_j)^k}\prod_{j=1}^{k} du_jdv_j, \nonumber
\end{align}
\noindent where $|\xi|=|b/a|<1$, as $|a| > 1$ and $|b| <1$. After the previous cancellation of the mixed
poles, the only possible poles in each $u_j$ are at $0$ and $\xi-1$, and the
only possible poles in each $v_\ell$ are at $0$ and $1-\xi$. Thus, the $u$-contours are centred at $-1$ and the $v$-contours are centred at $-\xi$; the contours have a radius greater than $1$. Fix
$0<\tau<1/2$; we may then deform the contours to positively oriented circles
$|u_j|=|v_\ell|=\mu$ such that $\mu>(1+\tau)|\xi-1|$, without crossing any poles. If we now return to \eqref{eq:fn}, for $|\xi| < 1$ we can write 
\begin{align}
&F_k(\xi;N)=(-1)^k\xi^k\frac{1}{(\xi-1)^{k^2}(k!)^2 (2\pi i)^{2k}}\int_{|u_j| =\mu} \int_{|v_{\ell}| = \mu} \prod_{j=1}^k (u_j+1)^{N} \Delta(u)^2 \Delta(v)^2 \\
&\,\qquad\qquad\qquad\qquad\qquad \times \prod_{j=1}^k\prod_{l=1}^k (v_l-u_j+\xi-1) \left(\sum_{j=1}^{k} \frac{1}{u_j}+\sum_{j=1}^k \frac{1}{v_j+\xi-1}-N\right)^k \nonumber \\
&\,\qquad\qquad\qquad\qquad\qquad  \times \left( \sum_{j=1}^{k}\frac{1}{v_j}+\sum_{j=1}^k \frac{1}{u_j-(\xi-1)} \right)^k \frac{1}{\prod_{j=1}^k u_j^k v_j^k}\prod_{j=1}^{k} du_jdv_j. \nonumber
\end{align}
Since $\prod_{j,l=1}^k(v_l-u_j+\xi-1)
=
(-1)^{k^2}\prod_{j,l=1}^k(u_j-v_l-(\xi-1))
$, the sign cancels the prefactor $(-1)^k$ in $F_k(\xi;N)$. In this form, we may apply Lemma \ref{lem:two-centre} with $L = N$, $\sigma = \xi -1$ and $\Phi_{N}(u_1, \ldots, u_k) = \prod_{j=1}^{k}(1+u_j)^{N}$ to give 
\begin{align}
\label{eq:FN-two-centre}
&F_k(\xi;N)=\xi^k\frac{1}{(2\pi i)^{2k}(k!)^2 }\sum_{r=0}^k\binom{k}{r}^2\, 
\sum_{a=0}^k
(-1)^{k-a}\binom{k}{a}\, \frac{N^{k-a}}{(\xi-1)^{k+a}}
\,
J_{N,r}^{(a)}(\xi -1), 
\end{align}
with
\begin{align*}
&J_{N,r}^{(a)}(\xi -1) =
\int_{c_{j}(r)}  \int_{\hat{c}_{\ell}(r)}  
\prod_{j=1}^k (1 + (\xi - 1)x_j)^{N} \,
\Delta(x)^2\Delta(y)^2\,
\prod_{j,\ell=1}^k (x_j-y_\ell-1)\,
\\
& \qquad  \qquad \qquad \, \times
\Biggl(\sum_{j=1}^k\frac1{x_j}+\sum_{\ell=1}^k\frac1{y_\ell+1}\Biggr)^{a}
\Biggl(\sum_{j=1}^k\frac1{x_j-1}+\sum_{\ell=1}^k\frac1{y_\ell}\Biggr)^{k}
\prod_{j=1}^k x_j^{-k} dx_j\prod_{\ell=1}^k y_\ell^{-k} dy_\ell. \nonumber
\end{align*}
Here, $c_{j}(r)$ are counter-clockwise circles of radius $\tau$ centered at 1 for $j=1,\ldots, r$ and 0 otherwise; similarly, $\hat{c}_{\ell}(r)$ are counter-clockwise circles of radius $\tau$ centered at $-1$ for $\ell = 1, ..., r$ and 0 otherwise.  The \textit{shifted} variables $x_1,\ldots,x_r$ are centered at $1$ and  $y_1,\ldots,y_r$ are centered at $-1$, and the remaining unshifted variables are centered at $0$. We now isolate the exponential dependence on $\xi$ in each $r$-term. For the shifted variables $j\le r$, we expand at the centre $x_j=1$:
\begin{equation*}
1+(\xi-1)x_j
=
1+(\xi-1)\bigl(1+(x_j-1)\bigr)
=
\xi+(\xi-1)(x_j-1),
\end{equation*}
and therefore
\begin{equation*}
\label{eq:shifted-binomial-concise}
(1+(\xi-1)x_j)^N
=
\sum_{\beta_j=0}^{N}\binom N{\beta_j}
\xi^{N-\beta_j}(\xi-1)^{\beta_j}(x_j-1)^{\beta_j}.
\end{equation*} 
We use the convention $\binom{N}{\beta}$ is zero if $N<\beta$. For the unshifted variables $j=r+1,\ldots,k$, we expand at the centre $x_j=0$:
\begin{equation*}
\label{eq:unshifted-binomial-concise}
(1+(\xi-1)x_j)^N
=
\sum_{\gamma_{j-r}=0}^{N}\binom N{\gamma_{j-r}}
(\xi-1)^{\gamma_{j-r}}x_j^{\gamma_{j-r}}.
\end{equation*}
Combining the two expansions above, we can write
\begin{equation*}
\prod_{j=1}^k (1 + (\xi - 1)x_j)^{N} =
\sum_{\beta\in\mathbb Z_{\ge0}^{r}}
\sum_{\gamma\in\mathbb Z_{\ge0}^{k-r}}
\phi_{\beta,\gamma}^{(r)}(N,\xi - 1)
\prod_{j=1}^{r}(x_j-1)^{\beta_j}
\prod_{j=r+1}^{k}x_j^{\gamma_{j-r}},
\end{equation*}
where
\begin{equation}
\label{eq:phi-concise}
\phi_{\beta,\gamma}^{(r)}(N,\xi-1)
=
\left(\prod_{j=1}^{r}\binom N{\beta_j}\right)
\left(\prod_{m=1}^{k-r}\binom N{\gamma_m}\right)
\xi^{Nr-|\beta|}(\xi-1)^{|\beta|+|\gamma|}.
\end{equation}
Now refer to the definitions from \eqref{eq:beta-defs}; every non-zero term has $\beta_j\le N$ for $1\le j\le r$, so
$|\beta|\le Nr$ and the power $\xi^{Nr-|\beta|}$ is non-negative. Using Lemma~\ref{lem:two-centre-degree-bound}, we have
\begin{equation}
\label{eq:J-zero-concise}
J_{N,r}^{(a)}(\xi-1)=0\qquad (a<r),
\end{equation}
and, for $a\ge r$,
\begin{align}
\label{eq:J-expansion-concise}
J_{N,r}^{(a)}(\xi-1)
&=
\sum_{\substack{\beta\in\mathbb Z_{\ge0}^{r}\\ |\beta|\le k-r}}
\sum_{\substack{\gamma\in\mathbb Z_{\ge0}^{k-r}\\ |\gamma|\le a-r}}
\mathcal C_{\beta,\gamma}^{(r,a)}
\left(\prod_{j=1}^{r}\binom N{\beta_j}\right)
\left(\prod_{m=1}^{k-r}\binom N{\gamma_m}\right)
\xi^{Nr-|\beta|}(\xi-1)^{|\beta|+|\gamma|},
\end{align}
where
\begin{align}
\label{eq:C-beta-gamma-def}
\mathcal C_{\beta,\gamma}^{(r,a)}
&:=
\int_{c_j(r)}\int_{\hat c_\ell(r)}
\prod_{j=1}^{r}(x_j-1)^{\beta_j}
\prod_{j=r+1}^{k}x_j^{\gamma_{j-r}}
\Delta(x)^2\Delta(y)^2
\prod_{j,\ell=1}^{k}(x_j-y_\ell-1)
\\
&\qquad\qquad\times
\Biggl(\sum_{j=1}^{k}\frac1{x_j}
+
\sum_{\ell=1}^{k}\frac1{y_\ell+1}\Biggr)^a
\Biggl(\sum_{j=1}^{k}\frac1{x_j-1}
+
\sum_{\ell=1}^{k}\frac1{y_\ell}\Biggr)^k
\prod_{j=1}^{k}x_j^{-k}\,dx_j
\prod_{\ell=1}^{k}y_\ell^{-k}\,dy_\ell . \nonumber
\end{align}
We substitute \eqref{eq:FN-two-centre} and
\eqref{eq:J-expansion-concise} into
\begin{equation}
\tilde I_{k}(n;N)
=
-\frac1{2\pi i}
\int_{|\xi|=\eta^2}\frac{F_k(\xi;N)}{\xi^{n+1}(\xi-1)}\,d\xi .
\end{equation}
For each admissible pair $(\beta,\gamma)$ this reduces the outer integral to
\begin{equation}
\label{eq:K-coeff-concise}
\frac1{2\pi i}
\int_{|\xi|=\eta^2}
\xi^{Nr-n+k-1-|\beta|}
(\xi-1)^{-d}\,d\xi,
\qquad
 d:=a+k+1-|\beta|-|\gamma|.
\end{equation}
The degree bounds imply $d>0$: if $r\ge1$, then $d\ge2r+1$, and if $r=0$,
then $d\ge k+1$.  Hence, for $|\xi|<1$,
\begin{equation}
(\xi-1)^{-d}=(-1)^d(1-\xi)^{-d}
=(-1)^d\sum_{q\ge0}\binom{q+d-1}{d-1}\xi^q,
\end{equation}
and Cauchy's integral formula gives the exact identity
\begin{equation}
\label{eq:cauchy-coeff}
\frac1{2\pi i}
\int_{|\xi|=\eta^2}
\frac{\xi^{Nr-n+k-1-|\beta|}}{(\xi-1)^{d}}\,d\xi
=
(-1)^d
\binom{n-Nr-k+|\beta|+d-1}{d-1}
\mathds{1}_{\{n-Nr-k+|\beta|\ge0\}}.
\end{equation}
It remains to identify the leading coefficient for these terms. For fixed
$\beta$ and $\gamma$,
\begin{equation}
\label{eq:leading-1}
\prod_{j=1}^{r}\binom N{\beta_j}
\prod_{m=1}^{k-r}\binom N{\gamma_m}
=
\frac{N^{|\beta|+|\gamma|}}{\beta!\,\gamma!}
+O(N^{|\beta|+|\gamma|-1}).
\end{equation}
We can isolate the role of the parameter $r$. Inserting the finite Taylor
expansion \eqref{eq:J-expansion-concise} into the outer integral \eqref{eq:outer-integral-for-cauchy} and applying Cauchy's theorem, the
contribution of the fixed pair $(\beta,\gamma)$ vanishes unless
\begin{equation}
\label{eq:cauchy-selection-condition-bettin}
    n-Nr-k+|\beta|\ge0.
\end{equation}
In particular, if $r>\alpha$, then
$n-Nr=(\alpha-r)N+O(1)$ tends to $-\infty$, while $|\beta|$ is bounded in
terms of $k$; hence \eqref{eq:cauchy-selection-condition-bettin} fails for
all sufficiently large $N$.  Thus no term with $r>\alpha$ contributes.  If
$r<\alpha$, then \eqref{eq:cauchy-selection-condition-bettin} is satisfied
for every admissible $\beta$ once $N$ is large enough.  The only remaining
case is the boundary case $r=\alpha$, which can occur only when $\alpha$ is
an integer.  In that case, the selected binomial coefficient in
\eqref{eq:cauchy-coeff} is $O(1)$ as $n-N
\alpha = O(1)$. With the prefactor $N^{k-a}$ from \eqref{eq:FN-two-centre} and using \eqref{eq:leading-1}, we have at 
most $N^{k-a+|\beta|+|\gamma|}=O(N^{2k-2})$ as $|\beta| \leq k-r$ and $|\gamma| \leq a - r$ due to the degree-bounds from Lemma \ref{lem:two-centre-degree-bound}. Therefore, we will see shortly that the boundary term is lower order and for the leading term, we proceed with only $r$ satisfying $r < \alpha$.

Since $r<\alpha$, we have
\begin{equation}
\label{eq:leading-2}
\binom{n-Nr-k+|\beta|+d-1}{d-1}
=
\frac{N^{d-1}(\alpha-r)^{d-1}}{(d-1)!}+O(N^{d-2}).
\end{equation}
As $|\beta|+|\gamma|+d-1=a+k$, the product of these two asymptotics contributes $N^{a+k}$.  After
multiplication by the prefactor $N^{k-a}$ in \eqref{eq:FN-two-centre}, each surviving $(r,a)$-term is therefore of size $N^{2k}$. We now explain how the finite sum over $\beta$ and $\gamma$ recombines into
the polynomial in the final formula. Let $q:=d-1=a+k-|\beta|-|\gamma|$. After including the minus sign from the outer integral, the leading coefficient multiplying $\prod_{j=1}^{r}(x_j-1)^{\beta_j}
\prod_{j=r+1}^{k}x_j^{\gamma_{j-r}} $ is
\begin{equation}
\label{eq:coefficient-match-noH}
\frac{(-1)^q(\alpha-r)^q}{q!\,\beta!\,\gamma!}.
\end{equation}
On the other hand,
\begin{equation}
\label{eq:shifted-sum-for-multinomial-noH}
\sum_{j=1}^{k}x_j-\alpha
=
\sum_{j=1}^{r}(x_j-1)+\sum_{j=r+1}^{k}x_j-(\alpha-r).
\end{equation}
Therefore, by the multinomial theorem,
\begin{align}
\label{eq:full-multinomial-expanded}
\sum_{\substack{q,\beta,\gamma\ge0\\ q+|\beta|+|\gamma|=a+k}}
\frac{(-1)^q(\alpha-r)^q}{q!\,\beta!\,\gamma!}
\prod_{j=1}^{r}(x_j-1)^{\beta_j}
\prod_{j=r+1}^{k}x_j^{\gamma_{j-r}} = \frac{\left(\sum_{j=1}^{k}x_j-\alpha\right)^{a+k}}{(a+k)!}.
\end{align}
Due to the multinomial sum collapse above and using the leading-order forms of both \eqref{eq:leading-1} and \eqref{eq:leading-2}, the $2k$-contour integral expression for $J_{N,r}^{a}(\xi - 1)$ in \eqref{eq:J-expansion-concise} is simplified greatly. As the $\xi$-integral has been evaluated, we have
\begin{align}
\label{eq:Itilde-leading-noH}
\tilde I_{k}(n;N)
&=
\frac{N^{2k}}{(2\pi i)^{2k} \, k!^2}
\, \sum_{r=0}^{k}\mathds{1}_{\{r<\alpha\}}\binom kr^2
\sum_{a=r}^{k} \, (-1)^{k-a} \, \binom ka
\frac1{(a+k)!}
\nonumber\\
&  \times
\int_{c_j(r)}\int_{\hat c_\ell(r)}
\Biggl(\sum_{j=1}^{k}x_j-\alpha\Biggr)^{a+k}
\Biggl(\sum_{j=1}^{k}\frac1{x_j}
+
\sum_{\ell=1}^{k}\frac1{y_\ell+1}\Biggr)^a
\Biggl(\sum_{j=1}^{k}\frac1{x_j-1}
+
\sum_{\ell=1}^{k}\frac1{y_\ell}\Biggr)^k 
\nonumber\\
&  \times
\Delta(x)^2\Delta(y)^2
\prod_{j,\ell=1}^{k}(x_j-y_\ell-1)
\prod_{j=1}^{k}x_j^{-k}\,dx_j
\prod_{\ell=1}^{k}y_\ell^{-k}\,dy_\ell
+O(N^{2k-1}).
\end{align}
Finally, dividing by $N^{2k}$ and using the definition of $M_k(\alpha)$ in \eqref{eq:rmt-mkalpha} yields the stated result in the theorem.

\subsection{Proof of Theorem \ref{mkalpha}}  Let $k \geq 1$ be an integer and assume $\alpha$ is a real positive non-integer; the integer case for $\alpha$ is treated separately at the end. We assume Conjecture~\ref{eq:main-conj} throughout with the error decay exponent $\theta_{k,c}$ defined as
\begin{equation}
\label{eq:theta_kc}
\theta_{k,c}=
\begin{cases}
1+c, & k=1,2,\\
k-1, & k\ge3,
\end{cases}
\end{equation}
for some fixed $0<c\leq 1$ and, in particular, when $k=1,2$ we choose $0<\xi< c$. We recall that $\Lambda^{*k}=\Lambda*\cdots*\Lambda$, hence
\begin{equation*}
\Lambda^{*k}(n)=\sum_{n_1\cdots n_k=n}\Lambda(n_1)\cdots\Lambda(n_k). 
\end{equation*}
For any ordered factorisation $n=n_1\cdots n_k$ we have $0\le \Lambda(n_j)\le \log n_j\le\log n$, so each summand is
bounded by $(\log n)^k$. The number of ordered $k$--tuples $(n_1,\dots,n_k)$ with product $n$ is precisely $d_k(n)$,
the $k$--fold divisor function. Therefore we have a pointwise bound
\begin{equation}
\label{eq:pointwise-dk}
\lvert\Lambda^{*k}(n)\rvert\le (\log n)^k d_k(n)\qquad (n\ge 2). 
\end{equation}
To obtain a truncated Perron formula for
\begin{equation*}
\sum_{n\le \lfloor T^\alpha\rfloor}\frac{\Lambda^{*k}(n)}{n^{1/2+it}},
\end{equation*}
we apply the Lemma stated below.
\begin{lemma}[Titchmarsh \cite{kn:titchmarsh}, Lemma 3.12]
\label{lemma:titchmarsh}
Let
\begin{equation*}
f(s)=\sum_{n=1}^{\infty}\frac{a_n}{n^{s}}\qquad (\Re{s} >1),
\end{equation*}
with $a_n=O\{\psi(n)\}$, $\psi(n)$ being non-decreasing, and
\begin{equation*}
\sum_{n=1}^{\infty}\frac{|a_n|}{n^{\sigma}}
=O\!\left\{\frac{1}{(\sigma-1)^{\mu}}\right\}
\end{equation*}
as $\sigma\to 1^{+}$. For $s = \sigma + it$ and $c>0$ so that $\sigma+c>1$, if $x$ is not an integer and $L$ is the
integer nearest to $x$ then
\begin{align*}
\sum_{n<x}\frac{a_n}{n^{s}}
&=\frac{1}{2\pi i}\int_{c-iT}^{c+iT} f(s+w)\,\frac{x^{w}}{w}\,dw
+O\!\left\{\frac{x^{c}}{T(\sigma+c-1)^{\mu}}\right\} \notag\\
&\qquad\qquad\qquad\qquad\qquad
+O\!\left\{\frac{\psi(2x)\,x^{1-\sigma}\log x}{T}\right\}
+O\!\left\{\frac{\psi(L)\,x^{1-\sigma}}{T\,|x-L|}\right\}.
\end{align*}
\end{lemma}
\noindent The abscissa of absolute convergence is $\sigma_a=1$, and, by positivity together with the estimate
\begin{equation}
\label{eq:unfortunatelysuperimportant}
-\frac{\zeta'}{\zeta}(\sigma)\asymp(\sigma-1)^{-1}
\qquad (1<\sigma\le2),
\end{equation}
meaning that the two quantities are bounded above and below by positive
constant multiples of one another uniformly in this range, we have
\begin{equation}
\label{eq:bounding-series}
\sum_{n\ge 1}\frac{|\Lambda^{*k}(n)|}{n^\sigma}
= \bigg(\sum_{n\ge 1}\frac{\Lambda(n)}{n^\sigma}\bigg)^{\!k}
=\bigg(-\frac{\zeta'(\sigma)}{\zeta(\sigma)}\bigg)^{\!k}
\ll_k (\sigma-1)^{-k}.
\end{equation}
Thus, the series bound condition of Lemma~\ref{lemma:titchmarsh} holds with $\mu=k$. For the complex sequence condition, we use the non-decreasing majorant
\begin{equation}
\label{eq:bdef}
B(x):=\max_{m\le x, \, m \in \mathbb{Z}_{+}}(\log m)^k d_k(m), \qquad x \geq 1, \, x \in \mathbb{R},
\end{equation}
for which the standard bound $d_k(m)\ll_{k,\varepsilon_1} m^{\varepsilon_1}$ for any fixed $\varepsilon_1>0$ \cite{iwaniec-kowalski} yields
$B(x)\ll_{k,\varepsilon_1} (\log x)^k x^{\varepsilon_1}$. If one integrates only up to height $T$, the resulting truncation error from applying the lemma is $o(1)$ only under restrictions on $\alpha$, namely $\alpha<2$ and $\alpha(\tfrac12+\varepsilon_1)<1$. To avoid these restrictions, we choose a larger Perron height. Fix $\varepsilon>0$ and take \begin{equation} \gamma:=\max(1,\alpha/2)+\varepsilon, \qquad T':=T^\gamma. \end{equation}
Fixing $\sigma=\tfrac12$ and letting $\kappa:=\frac12+\frac{1}{\log T'}$, Lemma~\ref{lemma:titchmarsh} gives
\begin{equation}
\label{eq:applied}
\sum_{n < N}\frac{\Lambda^{*k}(n)}{n^{1/2+it}}
= \frac{1}{2\pi i}\int_{\kappa-iT'}^{\kappa+iT'}\!
\Big(-\frac{\zeta'}{\zeta}(\tfrac12+it+w)\Big)^{\!k}\!\,\frac{N^{w}}{w}\,dw
\;+\;E_{\mathrm{trunc}}(t), 
\end{equation}
with
\begin{equation}\label{eq:trunc-from-lemma312}
E_{\mathrm{trunc}}(t)
\ \ll_{k} \, \frac{N^{\,\kappa}}{T'(\sigma+\kappa-1)^{k}}
\;+\;\frac{B(2N)\,N^{1-\sigma}\log N}{T'}
\;+\;\frac{B(\lfloor N \rfloor)\,N^{1-\sigma}}{T'\,|N-\lfloor N\rfloor|}.
\end{equation}
Here we set $N=\big\lfloor T^{\alpha}\big\rfloor+\tfrac12$, which is equidistant from two consecutive integers, so we choose $L=\lfloor N \rfloor$ as the nearest integer in \eqref{eq:trunc-from-lemma312}. Writing $N=(T')^\beta$, we have
\begin{equation*}
\beta:=\frac{\log N}{\log T'} =\frac{\alpha}{\gamma}+O\!\left(\frac{1}{\log T}\right) = \frac{\alpha}{\max(1,\alpha/2)+\varepsilon} +O\!\left(\frac{1}{\log T}\right)<2
\end{equation*}
for all sufficiently large $T$. Since $\sigma+\kappa-1=1/\log T'$, the truncation error becomes
\begin{equation}\label{eq:Etrunc-Tprime}
E_{\mathrm{trunc}}(t)
\ll_{k}
(\log T')^{k}\,(T')^{\beta/2-1}
+B(2N)\,(\log N)\,(T')^{\beta/2-1}
+B(N)\,(T')^{\beta/2-1},
\end{equation}
where the factor $N^{1/\log T'}=e^\beta$ has been absorbed into the implied constant. Using $B(N)\le B(2N) \ll_{k,\varepsilon_1}(\log N)^k N^{\varepsilon_1}$ yields the explicit form
\begin{equation*}
E_{\mathrm{trunc}}(t)
\ll_{k,\varepsilon_1, \alpha, \varepsilon}
(\log T')^{k}(T')^{\beta/2-1}
+(\log T')^{k+1}(T')^{\beta/2-1+\beta\varepsilon_1}
\end{equation*}
for any fixed $\varepsilon_1 > 0$. As $\limsup_{T\to\infty}\beta<2$, we may choose
$\varepsilon_1=\varepsilon_1(\alpha,\varepsilon)>0$ so small that
\begin{equation*}
\beta\left(\frac12+\varepsilon_1\right)<1
\end{equation*}
for all sufficiently large $T$. Hence there exists
$\epsilon'=\epsilon'_{\varepsilon,\varepsilon_1}(\alpha)>0$ such that, for all sufficiently large $T$,
\begin{equation*}
\frac{\beta}{2}-1\le -\epsilon'
\qquad\text{and}\qquad
\frac{\beta}{2}-1+\beta\varepsilon_1\le -\epsilon'.
\end{equation*}
Therefore
\begin{equation*}
\sum_{n < N}\frac{\Lambda^{*k}(n)}{n^{1/2+it}}
=\frac{1}{2\pi i}\int_{\kappa-iT'}^{\kappa+iT'}
\Big(-\frac{\zeta'}{\zeta}\big(\tfrac12+it+w\big)\Big)^{k}\frac{N^{w}}{w}\,dw
+E_{\mathrm{trunc}}(t),
\end{equation*}
where
\begin{equation}\label{eq:Etrunc-decay}
E_{\mathrm{trunc}}(t)
\ll_{k,\alpha,\varepsilon,\varepsilon_1}
(\log T')^{k+1}(T')^{-\epsilon'}, \qquad 0\le t\le T.
\end{equation}
Next, we choose a contour on the conjectural scale. Fix $0<\eta<\frac32$, and set
\begin{equation*}
\sigma_0:=\frac{\eta}{\log T}.
\end{equation*}
Then, for all sufficiently large $T$, we have $0<\sigma_0<\frac12<\kappa$.
We shift the vertical line of integration to $\Re w=\sigma_0$,
closing the rectangle by the two horizontal segments at heights
$\pm T'$. The integrand has a pole at $w=\tfrac12-it$, arising from the pole of $-\zeta'/\zeta(s)$ at $s=1$. This point
lies on $\Re w=\tfrac12$ and is enclosed provided $T'>|t|$, so we
pick up its residue. Assuming RH, the poles arising from the
non-trivial zeros $\rho$ of $\zeta(s)$ occur at $w=\rho-\bigl(\tfrac12+it\bigr)$
and lie on $\Re w=0$. Since we stop at
$\Re w=\sigma_0>0$, these poles are not enclosed, and no
zero-residue terms appear. The pole at $w=0$ from the Perron kernel
and the poles arising from the trivial zeros of $\zeta(s)$ all lie
to the left of $\Re w=\sigma_0$, so they are not crossed in this
shift. By the residue theorem,
\begin{align}
\label{eq:perrons-application}
\frac{1}{2\pi i}
\int_{\kappa-iT'}^{\kappa+iT'}
\Bigl(-\frac{\zeta'}{\zeta}
 \bigl(\tfrac12+it+z\bigr)\Bigr)^k
\frac{N^z}{z}\,dz
&=
\operatorname*{Res}_{w=1}
\left[
\Bigl(-\frac{\zeta'}{\zeta}(w)\Bigr)^k
\frac{N^{w-\frac12-it}}{w-\frac12-it}
\right] \nonumber
\\
&\quad+
\frac{1}{2\pi i}
\int_{\sigma_0-iT'}^{\sigma_0+iT'}
\Bigl(-\frac{\zeta'}{\zeta}
 \bigl(\tfrac12+it+z\bigr)\Bigr)^k
\frac{N^z}{z}\,dz
+E_{\mathrm{hor}}(t).
\end{align}
To estimate the horizontal contribution, $E_{\text{hor}}(t)$, write $z=u\pm iT'$ with $\sigma_0\le u\le \kappa$. Then
\begin{equation*}
\frac12+it+z
=
\frac12+u+i(t\pm T').
\end{equation*}
Since $T'=T^\gamma$ with $\gamma>1$, we have $T/T'\to0$ and $|t\pm T'|\asymp T'$ uniformly for $0\le t\le T$. Moreover
\begin{equation*}
\frac12+\sigma_0
\le
\Re\!\left(\frac12+it+z\right)
\le
1+\frac1{\log T'}.
\end{equation*}
By Titchmarsh's formula $(14.15.2)$ in \cite{kn:titchmarsh}, under RH, and using
$u\ge \eta/\log T$, we obtain
\begin{equation}
\label{eq:titchmarsh}
\frac{\zeta'}{\zeta}\!\left(\frac12+it+z\right)
\ll_{\eta, \alpha, \varepsilon}(\log T')^2
\end{equation}
on the horizontal edges. As $|z|\asymp T'$, $N^{\Re z}\le N^\kappa$, and the horizontal edges have bounded length, we obtain
\begin{equation}
E_{\mathrm{hor}}(t)
\ll_{k,\eta, \alpha, \varepsilon}
\frac{N^\kappa(\log T')^{2k}}{T'}
\end{equation}
uniformly for $0\le t\le T$. As $N^{\kappa}=N^{\frac12+\frac1{\log T'}}=N^{1/2}N^{1/\log T'}$ and
\begin{equation*}
N^{1/\log T'}
=
\exp\!\Big(\frac{\log N}{\log T'}\Big)
=
e^\beta,
\end{equation*}
we have $N^{\kappa}=e^{\beta}N^{1/2}$. As before, the factor $e^{\beta}=O(1)$ may be absorbed into the implied constant. Therefore $N^\kappa=T^{\alpha/2+o(1)}$ and  as $T'^{-1}=T^{-\gamma}$, it follows that
\begin{equation}
\label{eq:horizontal-error-supremum}
\sup_{0 \le t \le T}|E_{\mathrm{hor}}(t)|
\ll_{k,\eta, \alpha, \varepsilon}
(\log T')^{2k}T^{\alpha/2-\gamma+o(1)}.
\end{equation}
Since 
\begin{equation}
\gamma-\frac{\alpha}{2} = \max(1-\alpha/2,0)+\varepsilon \ge \varepsilon, \end{equation}
the exponent of $T$ is negative for all sufficiently large $T$. From this point, assume the implied constants depend on the fixed parameters $k,\eta,\alpha,\varepsilon$, and $\varepsilon_1$. It is convenient to define $E(t) = E_{\text{hor}}(t) + E_{\text{trunc}}(t)$; we have
\begin{equation}
|E(t)|\ \ll \mathcal M, \qquad 
\mathcal M:= (\log T')^{2k}\,T^{\alpha/2-\gamma+o(1)} + (\log T')^{k+1}\,T'^{-\epsilon'}.
\end{equation}
Since $\log T'\asymp_{\alpha, \varepsilon} \log T$, we may write
\begin{equation*}
\mathcal M \ \ll\ (\log T)^{2k}\,T^{-\delta},
\end{equation*}
where
\begin{equation}
\delta:=\frac12\min\!\Big(\gamma-\tfrac{\alpha}{2},\ \gamma\,\epsilon'\Big)>0.
\end{equation}
Here the factor $\gamma-\alpha/2$ comes from the estimate for
$E_{\mathrm{hor}}(t)$ in \eqref{eq:horizontal-error-supremum}, while $\gamma\epsilon'$ follows from
\eqref{eq:Etrunc-decay} upon writing
$(T')^{-\epsilon'}=T^{-\gamma\epsilon'}$; the factor $\tfrac 12$ allows the $o(1)$ term to be absorbed for sufficiently large
$T$. For each fixed $\alpha>0$, this gives a positive value of $\delta$. Introduce
\begin{equation}
V(t):=\frac{1}{2\pi i}\int_{\sigma_0-iT'}^{\sigma_0+iT'}\bigg(-\frac{\zeta'}{\zeta}\big(\tfrac12+it+z\big)\bigg)^{k}\frac{N^{z}}{z}\,dz,
\qquad 
\mathcal I:=\frac{1}{T}\int_{0}^{T}|V(t)|^{2}\,dt.
\end{equation}
We want to obtain an asymptotic expression of $\mathcal I_k(T,N)$ as defined in \eqref{Ik-truncated-moments}. Hence, by using \eqref{eq:perrons-application} and the definition of $\mathcal I_k(T,N)$ in \eqref{Ik-truncated-moments}, we may write
\begin{equation*}
\mathcal I_k(T,N)
=
\frac1T\int_0^T |V(t)+E(t)|^2\,dt.
\end{equation*}
Using $|A+B|^{2}=|A|^{2}+2\Re(A\overline{B})+|B|^{2}$, we obtain
\begin{align}
\mathcal I_k(T,N)-\mathcal I
=&\ \frac{2}{T}\Re\!\int_{0}^{T}V(t)\overline{E(t)}\,dt
+\frac{1}{T}\int_{0}^{T}|E(t)|^{2}\,dt,
\end{align}
and by the Cauchy--Schwarz inequality,
\begin{align}
\left|\mathcal I_k(T,N)-\mathcal I\right|
\le&\ 
2\sqrt{\mathcal I}\,
\Big(\frac{1}{T}\int_{0}^{T}|E(t)|^{2}\,dt\Big)^{1/2}
+\frac{1}{T}\int_{0}^{T}|E(t)|^{2}\,dt .
\end{align}
Using $|E(t)|\ll \mathcal M$, we deduce
\begin{equation}
\label{eq:full-cauchy-schwarz}
\mathcal I_k(T,N)=\mathcal I+O \!\left(\sqrt{\mathcal I}\,\mathcal M+\mathcal M^{2}\right),
\end{equation}
where
\begin{align}
\label{eq:the-real-one-unfortunately}
\mathcal I
=\frac{1}{4\pi^{2}T}\int_{\sigma_0 -iT'}^{\sigma_0 +iT'}\int_{-\sigma_0 -iT'}^{-\sigma_0 +iT'}
\frac{N^{z_{1}-z_{2}}}{z_{1}\,z_{2}}
\left[\int_{0}^{T}
\biggl(\frac{\zeta'}{\zeta}\big(\tfrac12+it+z_{1}\big)\biggr)^{k} \hspace{-0.5mm}
\biggl(\frac{\zeta'}{\zeta}\big(\tfrac12-it -z_{2}\big)\biggr)^{k}dt \right] dz_{2}dz_{1}.
\end{align}
It remains to bound the error term
$\sqrt{\mathcal I}\mathcal M+\mathcal M^2$. On the vertical line $z=\sigma_0+iv$, we have
\begin{equation*}
\frac12+it+z=\frac12+\sigma_0+i(t+v),
\end{equation*}
with $|t+v|\leq T+T'\ll T'$. When $|t+v|$ is sufficiently
large, the same argument used to obtain \eqref{eq:titchmarsh} gives
\begin{equation}
\label{eq:full-range-bound-titch}
\frac{\zeta'}{\zeta}
\left(\frac12+\sigma_0+i(t+v)\right)
\ll_{\eta,\alpha,\varepsilon}(\log T')^2.
\end{equation}
When $|t+v|$ is bounded, a local partial-fraction expansion gives
\begin{equation*}
\frac{\zeta'}{\zeta}
\left(\frac12+\sigma_0+i(t+v)\right)
\ll \sigma_0^{-1}\ll_\eta\log T.
\end{equation*}
Therefore, the first estimate \eqref{eq:full-range-bound-titch} holds uniformly for
$0\leq t\leq T$ and $|v|\leq T'$.
Using $N^{\sigma_0}
=
\exp\!\left(\frac{\eta\log N}{\log T}\right)
=
e^{\alpha\eta+o(1)}
=
O_{\alpha,\eta}(1)$ gives
\begin{equation}
|V(t)|
\ll
(\log T')^{2k}
\int_{-T'}^{T'}\frac{dv}{|\sigma_0+iv|}
\ll
(\log T')^{2k}\log T.
\end{equation}
Hence $\mathcal I\ll(\log T)^{4k+2}$, and as $\mathcal M\ll(\log T)^{2k}T^{-\delta}$, it follows that $\sqrt{\mathcal I}\,\mathcal M+\mathcal M^2
=o(1)$ as $T\to\infty$. At this stage, Conjecture~\ref{eq:main-conj} may be applied to the innermost integral in $\mathcal{I}$ with $Y=2T'$. As $\sigma_0=\eta/\log T$ with $0<\eta<3/2$, on the contours
$\Re z_1=\sigma_0$ and $\Re z_2=-\sigma_0$ we have
\begin{equation*}
0<\Re z_1,\ -\Re z_2,\ \Re(z_1-z_2)\le \frac{3}{\log T}.
\end{equation*}
Applying Conjecture~\ref{eq:main-conj} gives the error term
\begin{equation*}
O\!\left(
(\log T)^{2k-1+\xi}
\bigl(1+\min(|z_1-z_2|,1)\log T\bigr)^{-\theta_{k,c}}
\right)
\end{equation*}
for the normalised innermost integral with a fixed $\xi > 0$. In the cases $k=1,2$, where $\theta_{k,c}=1+c$, we choose $0<\xi<c$; for $k\ge3$, we take any fixed $0<\xi<1$. Going forth, the implied constants may also depend on $\xi$ and $c$. Note that
\begin{equation*}
\left|N^{z_1-z_2}\right|
=
N^{2\sigma_0}
=
\exp\!\left(\frac{2\eta\log N}{\log T}\right)
\ll 1.
\end{equation*}
Thus the contribution of the conjectural error term, which we will denote by $\mathcal{I}_e$, is not enlarged by the outer factor $N^{z_1-z_2}$ from \eqref{eq:the-real-one-unfortunately}. It follows that 
\begin{equation}
\mathcal I_e
\ll
(\log T)^{2k-1+\xi}
\int_{\sigma_0-iT'}^{\sigma_0+iT'}
\int_{-\sigma_0-iT'}^{-\sigma_0+iT'}
\frac{
\bigl(1+\min(|z_1-z_2|,1)\log T\bigr)^{-\theta_{k,c}}
}{
|z_1|\,|z_2|
}\,|dz_2|\,|dz_1|.
\end{equation}
Here $|dz_j|$ denotes the positive length element along the contour,
rather than the directed complex differential $dz_j$; it appears because the triangle inequality has been used to bound the absolute value of the
contour integral. On the vertical lines, write $z_1=\sigma_0+it_1$ and $z_2=-\sigma_0+it_2$. Then
\begin{equation*}
|z_1|^{-1}\ll \frac{1}{\sigma_0+|t_1|},
\qquad
|z_2|^{-1}\ll \frac{1}{\sigma_0+|t_2|}.
\end{equation*}
We may also write
\begin{equation*}
1+\min(|z_1-z_2|,1)\log T
\asymp
1+\min(|t_1-t_2|,1)\log T
\end{equation*}
as $\sigma_0=\eta/\log T$. This gives
\begin{equation}
\mathcal I_e
\ll
(\log T)^{2k-1+\xi}
\int_{-T'}^{T'}\int_{-T'}^{T'}
\frac{
\bigl(1+\min(|t_1-t_2|,1)\log T\bigr)^{-\theta_{k,c}}
}{
(\sigma_0+|t_1|)(\sigma_0+|t_2|)
}\,dt_2\,dt_1.
\end{equation}
To estimate this double integral, we define
\begin{equation*}
f(t):=\frac{1}{\sigma_0+|t|},
\end{equation*}
and let $u=t_1-t_2$. For fixed $u$, on writing $t=t_1$, the original
range of integration is
\begin{equation*}
t\in[-T',T']\cap[u-T',u+T'].
\end{equation*}
Since the integrand is non-negative, we may enlarge this range to
$[-T',T']$, so the double integral is bounded by
\begin{equation}
\label{eq:u-integral-bound}
\int_{-2T'}^{2T'}
\bigl(1+\min(|u|,1)\log T\bigr)^{-\theta_{k,c}}
\left(\int_{-T'}^{T'}f(t)f(t-u)\,dt\right)\,du.
\end{equation}
For $|u|\leq\sigma_0$, the Cauchy--Schwarz inequality gives
\begin{align*}
\int_{-T'}^{T'}f(t)f(t-u)\,dt
&\leq
\left(\int_{\mathbb R}f(t)^2\,dt\right)^{1/2}
\left(\int_{\mathbb R}f(t-u)^2\,dt\right)^{1/2} \\
&=
\int_{\mathbb R}\frac{dt}{(\sigma_0+|t|)^2}
\ll\frac{1}{\sigma_0}.
\end{align*}
Suppose now that $|u|>\sigma_0$. We split the $t$-integral into the
regions $|t|\leq |u|/2$, $|t-u|\leq |u|/2$, and their complement.
On the first region, the reverse triangle inequality gives $|t-u|\geq |u|-|t|\geq\frac{|u|}{2}$, and therefore
\begin{equation*}
f(t-u)\ll\frac{1}{\sigma_0+|u|}.
\end{equation*}
It follows that
\begin{align*}
\int_{|t|\leq |u|/2}f(t)f(t-u)\,dt
&\ll
\frac{1}{\sigma_0+|u|}
\int_{|t|\leq |u|/2}\frac{dt}{\sigma_0+|t|} =
\frac{2}{\sigma_0+|u|}
\log\left(1+\frac{|u|}{2\sigma_0}\right) \\
&\ll
\frac{1+\log(|u|/\sigma_0)}{\sigma_0+|u|}.
\end{align*}
The region $|t-u|\leq |u|/2$ is handled identically. On the remaining region, both $|t|$ and $|t-u|$ exceed $|u|/2$.
On the part where $|t|\leq2|u|$, the integrand is
$O(|u|^{-2})$ over an interval of length $O(|u|)$, giving a
contribution $O(|u|^{-1})$. For the tails $|t|>2|u|$, we have $|t-u|\geq |t|-|u|>\frac{|t|}{2}$, 
so the integrand is $O(t^{-2})$ and contributes
\begin{equation*}
\ll\int_{|t|>2|u|}\frac{dt}{t^2}
\ll\frac{1}{|u|}
\ll\frac{1}{\sigma_0+|u|}.
\end{equation*}
Combining these estimates and using
$\sigma_0=\eta/\log T$, with $\eta$ fixed, gives
\begin{equation}
\label{eq:convolution-bound}
\int_{-T'}^{T'}f(t)f(t-u)\,dt
\ll
\frac{1+\log(2+|u|\log T)}{\sigma_0+|u|}.
\end{equation}
In order to bound \eqref{eq:u-integral-bound} explicitly using the convolution bound \eqref{eq:convolution-bound}, we first split the $u$-integral into the ranges
\begin{equation*}
|u|\le \frac{1}{\log T},\qquad
\frac{1}{\log T}<|u|\le 1,\qquad
1<|u|\le 2T'.
\end{equation*}
The first range contributes $O(1)$, while the second supplies
\begin{equation*}
\ll
\int_{1/\log T<|u|\le 1}
\bigl(|u|\log T\bigr)^{-\theta_{k,c}}
\frac{\log(2+|u|\log T)}{|u|}\,du.
\end{equation*}
By changing variables $x=|u|\log T$, this is
\begin{align*}
\int_1^{\log T}
x^{-\theta_{k,c}}\frac{\log(2+x)}{x}\,dx
&\ll
\int_1^{\log T}x^{-\theta_{k,c}}\,dx \\
&\leq
\int_1^\infty x^{-\theta_{k,c}}\,dx
=
\frac{1}{\theta_{k,c}-1}
\ll_{k,c}1,
\end{align*}
as $\theta_{k,c}>1$ and $x^{-1}\log{(2+x)} \ll 1$ for $x\geq 1$. In addition, on the third range we have
$\min(|u|,1)=1$ and $\sigma_0+|u|\asymp |u|$. Hence its contribution
is
\begin{align*}
&\ll
(\log T)^{-\theta_{k,c}}
\int_{1<|u|\leq2T'}
\frac{\log(2+|u|\log T)}{|u|}\,du
&=
2(\log T)^{-\theta_{k,c}}
\int_1^{2T'}
\frac{\log(2+u\log T)}{u}\,du.
\end{align*}
As $\log(2+u\log T)
\leq\log(2+2T'\log T)
\ll\log T'$ throughout this range, it follows that
\begin{equation*}
\ll
(\log T)^{-\theta_{k,c}}(\log T')
\int_1^{2T'}\frac{du}{u}
\ll
(\log T)^{-\theta_{k,c}}(\log T')^2
\ll
(\log T)^{2-\theta_{k,c}},
\end{equation*}
where the final bound uses $\log T'=\gamma\log T$. We may combine the estimates from these ranges to conclude that the double integral in $\mathcal I_e$ is $\ll
1+(\log T)^{2-\theta_{k,c}}$ and, consequently,
\begin{equation*}
\mathcal I_e
\ll
(\log T)^{2k-1+\xi}
\left(1+(\log T)^{2-\theta_{k,c}}\right).
\end{equation*}
If $k\ge3$, then $\theta_{k,c}=k-1\ge2$, so $\mathcal I_e
\ll (\log T)^{2k-1+\xi} = o((\log T)^{2k})$ provided $0<\xi<1$.
If $k=1,2$, then we have $\theta_{k,c}=1+c$, and hence
\begin{equation*}
\mathcal I_e
\ll
(\log T)^{2k-1+\xi}
\left(1+(\log T)^{1-c}\right)
\ll
(\log T)^{2k+\xi-c},
\end{equation*}
which is $o((\log T)^{2k})$ for $0<\xi<c$. Overall, with our earlier choice of $\xi$, we have $\mathcal I_e=o((\log T)^{2k})$ for every $k\ge1$. 

Next, we return to the leading-order term from applying Conjecture~\ref{eq:main-conj} with $\mathcal C_u = \mathcal{C}_u(z)$ and $\mathcal C_v = \mathcal C_v(z) $ being positively oriented piecewise smooth contours where $\mathcal C_u$ encircles once $0$ and $-z$, and $\mathcal C_v$ encircles once $0$ and $z$. For the above choice of $\xi$, we have
\begin{align*}
\mathcal I_k(T,N)
=& \frac{1}{4\pi^2 \, k!^2}
\int_{\sigma_0-iT'}^{\sigma_0+iT'}
\int_{-\sigma_0-iT'}^{-\sigma_0+iT'}
z^{-k^2}\frac{N^{z}}{z_1z_2} \mathcal{W}_{k}(z; T) dz_2 dz_1 +
o\bigl((\log T)^{2k}\bigr)
\end{align*}
with
\begin{align}
\label{eq:Wk-def}
& \mathcal{W}_{k}(z; T) = \frac{(-1)^{k}}{(2\pi i)^{2k}}\int_{\mathcal C_u^k}\int_{\mathcal C_v^k}
T^{\frac12\left(\sum_{j=1}^k u_j-\sum_{\ell=1}^k v_\ell\right)}
\prod_{1\le i<j\le k}(u_i-u_j)^2
\prod_{1\le i<j\le k}(v_i-v_j)^2 
\nonumber\\
&\qquad \qquad \times
\biggl(
\sum_{j=1}^{k}\frac1{u_j}
+
\sum_{\ell=1}^{k}\frac1{v_\ell-z}
-\frac{\log T}{2}
\biggr)^k
\biggl(
\sum_{j=1}^{k}\frac1{u_j+z}
+
\sum_{\ell=1}^{k}\frac1{v_\ell}
+\frac{\log T}{2}
\biggr)^k
\nonumber\\
&\qquad \qquad \times
 \prod_{j,\ell=1}^{k}(u_j-v_\ell+z) \, \prod_{j=1}^{k}u_j^{-k}\,du_j
\prod_{\ell=1}^{k}v_\ell^{-k}\,dv_\ell 
\end{align} 
where $z=z_1 - z_2$. Note that the $(-1)^{k}$ factor above is from rewriting $\prod_{j,\ell}(v_\ell - u_j - z)$ as $\prod_{j,\ell}(u_j - v_\ell + z)$. For each fixed $z$ on $\Re{z} = 2\sigma_0$, we now use Lemma~\ref{lem:P-to-Q-form} to rewrite $\mathcal{W}_{k}(z; T)$ from $\mathcal I_k(T,N)$ in a simpler and unsymmetrised form. Write $U=\sum_{j=1}^{k}u_j$, $V=\sum_{\ell=1}^{k}v_\ell$ and set $L:=\log T$ so that $T^{\frac12(U-V)}=e^{\frac L2(U-V)}$.
\mycomment{The $u,v$-variables dependent part of the integrand in \eqref{eq:before-unsym} is then
\begin{equation}
e^{\frac L2(U-V)}
W_z(u,v)
\left(A_z(u,v)-\frac L2\right)^k
\left(B_z(u,v)+\frac L2\right)^k,
\end{equation}
where $A_z$, $B_z$ and $W_z$ are as in Lemma~\ref{lem:P-to-Q-form}.}By applying the Lemma with $s=\frac L2$ and $t=-\frac L2$, we have
\mycomment{\begin{align}
&\operatorname*{Res}_{u,v}
\, T^{\frac12(U-V)}
W_z(u,v)
\left(A_z(u,v)-\frac{\log T}{2}\right)^k
\left(B_z(u,v)+\frac{\log T}{2}\right)^k
\nonumber \\
&\qquad\qquad\qquad\qquad\qquad\qquad =
\operatorname*{Res}_{u,v}
\, T^U
W_z(u,v)
\left(A_z(u,v)-\log T\right)^k
B_z(u,v)^k.
\end{align}}
\begin{align}
\mathcal{W}_{k}(z;T)
=&
\, \, \frac{(-1)^{k}}{(2\pi i)^{2k}}
 \int_{\mathcal C_u^k}\int_{\mathcal C_v^k}
T^{\sum_{j=1}^{k}u_j}
\Delta(u)^2\Delta(v)^2
\nonumber \\
&\times \left(
\sum_{j=1}^{k}\frac{1}{u_j}
+
\sum_{\ell=1}^{k}\frac{1}{v_\ell-z}
-\log T
\right)^k
\left(
\sum_{j=1}^{k}\frac{1}{u_j+z}
+
\sum_{\ell=1}^{k}\frac{1}{v_\ell}
\right)^k
\nonumber \\
&\times
\prod_{j,\ell=1}^{k}(u_j-v_\ell+z) \prod_{j=1}^{k}u_j^{-k} du_{j}
\prod_{\ell=1}^{k}v_\ell^{-k}dv_{\ell}.
\end{align} 
Thus, the symmetric pair of shifts
$(-\log T/2,+\log T/2)$ has been replaced by the unsymmetrised pair $(-\log T,0)$ at the level of the $2k$-variable residue $\mathcal{W}_{k}(z;T)$. To further decompose this residue expression, Lemma~\ref{lem:two-centre} can be applied by setting $\sigma = - z$, $L= \log{T}$ and by using the symmetric function $\Phi_{L}(u_1, \ldots, u_k) = T^{\sum_{j=1}^{k} u_j}$. This gives
\begin{align}
\label{eq:two-centres-Ik}
\mathcal I_k(T,N)
&=
\frac{1}{4 \pi^2 k!^2}
\sum_{r=0}^{k}\binom{k}{r}^{2}
\sum_{a=0}^{k}
\binom{k}{a}
L^{k-a}
\int_{\sigma_0-iT'}^{\sigma_0+iT'}
\int_{-\sigma_0-iT'}^{-\sigma_0+iT'}
\frac{N^{z}}{z_1z_2} z^{-a-k}J_{r}^{(a)}(-z;T)
dz_2 dz_1
\nonumber\\
&\qquad\qquad\qquad\qquad\qquad\qquad\qquad\qquad\qquad\qquad\qquad\qquad\quad
+
o\bigl((\log T)^{2k}\bigr)
\end{align}
where
\begin{align}
\label{def:J}
J_{r}^{(a)}(-z; T)
&= \frac{1}{(2\pi i)^{2k}}
\int_{c_j(r)}\int_{\hat c_\ell(r)}
T^{-z\sum_{j=1}^{k}x_j}
\Delta(x)^2\Delta(y)^2
\prod_{j,\ell=1}^{k}(x_j-y_\ell-1)
\nonumber\\
&\quad\times
\left(
\sum_{j=1}^{k}\frac1{x_j}
+
\sum_{\ell=1}^{k}\frac1{y_\ell+1}
\right)^a
\left(
\sum_{j=1}^{k}\frac1{x_j-1}
+
\sum_{\ell=1}^{k}\frac1{y_\ell}
\right)^k
\prod_{j=1}^{k}x_j^{-k}\,dx_j
\prod_{\ell=1}^{k}y_\ell^{-k}\,dy_\ell
\end{align}
with $c_j(r)$ and $\hat c_\ell(r)$ as defined in Lemma~\ref{lem:two-centre}. For each fixed pair $(r,a)$ in the sum above, the change of variables
$z=z_1-z_2$ gives the following collapsed form of the double Perron kernel:
\begin{align}
&-\frac{1}{(2\pi i)^2}
\int_{\sigma_0-iT'}^{\sigma_0+iT'}
\int_{-\sigma_0-iT'}^{-\sigma_0+iT'}
 z^{-a-k}J_r^{(a)}(-z;T)
\frac{N^z}{z_1z_2}\,\dd z_2\dd z_1 
\\
&\qquad \qquad\qquad\qquad\qquad\qquad\qquad =
\frac{1}{2\pi i}
\int_{2\sigma_0-2iT'}^{2\sigma_0+2iT'}
 z^{-a-k}J_r^{(a)}(-z;T)N^z K_{T'}(z)\,\dd z, \nonumber
\label{eq:collapsed-termwise}
\end{align}
where
\begin{equation}
\label{eq:perron-kernel-estimate-final}
K_{T'}(z):= -
\frac{1}{2\pi i}
\int_{-\sigma_0+i\max(-T',-T'-\Im z)}^{-\sigma_0+i\min(T',T'-\Im z)}
\frac{\dd z_2}{(z+z_2)z_2}
=
\frac1z+O\!\left(\frac{\log T}{T'}\right)
\end{equation}
uniformly for $z=2\sigma_0+iy$, $|y|\le 2T'$. We prove \eqref{eq:perron-kernel-estimate-final} by considering two
ranges; suppose first that $|y|\leq T'/2$ and close the
$z_2$-contour to the left. The pole at $z_2=-z$ is enclosed,
whereas the pole at $z_2=0$ lies to the right. Moreover,
\begin{equation*}
\operatorname*{Res}_{z_2=-z}
\frac{1}{(z+z_2)z_2}
=
-\frac1z,
\end{equation*}
so the minus sign in the definition of $K_{T'}(z)$ gives the main
term $1/z$. To estimate the two horizontal parts of the closed
contour, we write $z_2=x+iv$, with
$x\leq-\sigma_0$, and use $|y|\leq T'/2$ to get
\begin{equation*}
|z_2|\gg T'+|x|,
\qquad
|z+z_2|\gg T'+|x|.
\end{equation*}
Their total contribution is therefore
\begin{equation*}
\ll
\int_{-\infty}^{-\sigma_0}
\frac{dx}{(T'+|x|)^2}
\ll \frac1{T'},
\end{equation*}
and this gives
\begin{equation}
\label{eq:KT-small-y}
K_{T'}(z)
=
\frac1z+O\!\left(\frac1{T'}\right),
\qquad |y|\leq \frac{T'}{2}.
\end{equation}
Now we consider the case $T'/2<|y|\leq2T'$. Writing
$z_2=-\sigma_0+iv$ and applying
\eqref{eq:convolution-bound} with $u=-y$, we obtain
\begin{align*}
|K_{T'}(z)|
&\ll
\int_{-T'}^{T'}
\frac{dv}
{(\sigma_0+|v|)(\sigma_0+|v+y|)} \ll
\frac{1+\log(2+|y|\log T)}{\sigma_0+|y|}
\ll
\frac{\log T}{T'}.
\end{align*}
Here the last estimate follows from $T'/2<|y|\leq2T'$ and
$T'=T^\gamma$. In addition, $|z|\geq|y|>T'/2$, and hence
$|z|^{-1}\ll T'^{-1}$. It follows that
\begin{equation}
\label{eq:KT-large-y}
\left|K_{T'}(z)-\frac1z\right|
\leq |K_{T'}(z)|+\frac1{|z|}
\ll\frac{\log T}{T'}.
\end{equation}
Combining \eqref{eq:KT-small-y} and \eqref{eq:KT-large-y} proves
\eqref{eq:perron-kernel-estimate-final} uniformly for
$|y|\leq2T'$.

The above kernel estimate reduces the remaining task to controlling the $z$-dependence of $J_r^{(a)}(-z;T)$, and this dependence is best exposed by separating the contribution of the two contour centres. Since
\begin{equation}
\label{eq:pretty-obvious}
\sum_{j=1}^{k}x_j
=
r+
\sum_{j=1}^{r}(x_j-1)
+
\sum_{j=r+1}^{k}x_j,
\end{equation}
we have the factorisation
\begin{equation}
J_r^{(a)}(-z;T)=T^{-rz}\tilde J_r^{(a)}(z;T)
\label{eq:J-factor-revised}
\end{equation}
with
\begin{align}
\tilde J_r^{(a)}(z;T)
&:=
\frac{1}{(2\pi i)^{2k}}
\int_{c_j(r)}\int_{\hat c_\ell(r)}
\exp\left(
-z\log T
\left\{
\sum_{j=1}^{r}(x_j-1)
+
\sum_{j=r+1}^{k}x_j
\right\}
\right)
\nonumber\\
&\quad\times \left(
\sum_{j=1}^{k}\frac1{x_j}
+
\sum_{\ell=1}^{k}\frac1{y_\ell+1}
\right)^a
\left(
\sum_{j=1}^{k}\frac1{x_j-1}
+
\sum_{\ell=1}^{k}\frac1{y_\ell}
\right)^k \prod_{j,\ell=1}^{k}(x_j-y_\ell-1) 
\nonumber\\
&\quad\times
\Delta(x)^2\Delta(y)^2
\prod_{j=1}^{k}x_j^{-k}\,dx_j
\prod_{\ell=1}^{k}y_\ell^{-k}\,dy_\ell .
\label{eq:Jtilde-def-revised}
\end{align}
We now bound the $z$-dependence in $\tilde J_r^{(a)}(z;T)$. First, we expand the exponential factor at the two centres to give
\begin{align}
&\exp\left(
-z\log T
\left\{
\sum_{j=1}^{r}(x_j-1)
+
\sum_{j=r+1}^{k}x_j
\right\}
\right)
\nonumber\\
&\qquad\qquad\qquad \qquad \qquad =
\sum_{\beta\in\mathbb Z_{\ge0}^{r}}
\sum_{\gamma\in\mathbb Z_{\ge0}^{k-r}}
\frac{(-z\log T)^{|\beta|+|\gamma|}}{\beta!\gamma!}
\prod_{j=1}^{r}(x_j-1)^{\beta_j}
\prod_{j=r+1}^{k}x_j^{\gamma_{j-r}} .
\label{eq:Jtilde-expansion-revised}
\end{align}
This is precisely the Taylor expansion appearing in
Lemma~\ref{lem:two-centre-degree-bound}, with
$\sigma=-z$ and $\Phi_{\log T}(u_1,\ldots,u_k)=T^{\sum_j u_j}$,
after the scalar factor $T^{-rz}$ has been removed. By Lemma~\ref{lem:two-centre-degree-bound}, if $a<r$, then
\begin{equation}
\tilde J_r^{(a)}(z;T)=0.
\label{eq:Jtilde-zero-revised}
\end{equation}
If $a\ge r$, then the residue is unchanged after truncating
\eqref{eq:Jtilde-expansion-revised} to the terms satisfying
\begin{equation*}
|\beta|\le k-r,
\qquad
|\gamma|\le a-r.
\end{equation*}
Consequently, for $1\le r\le k$, only powers
\begin{equation*}
(z\log T)^m,
\qquad
0\le m\le (k-r)+(a-r)=k+a-2r,
\end{equation*}
can contribute to $\tilde J_r^{(a)}(z;T)$. This gives the bound
\begin{equation}
\tilde J_r^{(a)}(z;T)
\ll_k
\bigl(1+|z|\log T\bigr)^{k+a-2r},
\qquad
1\le r\le k,\quad a\ge r.
\label{eq:Jtilde-bound-r-positive-revised}
\end{equation}
For $r=0$, the shifted variable product in \eqref{eq:Jtilde-expansion-revised} is empty and Lemma~\ref{lem:two-centre-degree-bound}
leaves only the restriction $|\gamma|\le a$. Therefore
\begin{equation}
\tilde J_0^{(a)}(z;T)
\ll_k
\bigl(1+|z|\log T\bigr)^a.
\label{eq:Jtilde-bound-r-zero-revised}
\end{equation}
Thus
\begin{equation}
\tilde J_r^{(a)}(z;T)=0
\qquad (a<r),
\end{equation}
and, for $a\ge r$,
\begin{equation}
\tilde J_r^{(a)}(z;T)
=
\begin{cases}
O_k\bigl((1+|z|\log T)^a\bigr), & r=0,\\[3pt]
O_k\bigl((1+|z|\log T)^{k+a-2r}\bigr), & 1\le r\le k.
\end{cases}
\label{eq:Jtilde-final-bound-revised}
\end{equation}
\noindent We may now return to the kernel $K_{T'}$ in
\eqref{eq:perron-kernel-estimate-final}. For $a\ge r$, the
$(r,a)$-term in $\mathcal{I}_{k}(T,N)$ contains the integral
\begin{equation}
\frac{1}{2\pi i}
\int_{2\sigma_0-2iT'}^{2\sigma_0+2iT'}
(\log T)^{k-a}
\left(\frac{N}{T^r}\right)^z
z^{-a-k}
\tilde J_r^{(a)}(z;T)
K_{T'}(z)\,\dd z.
\label{eq:collapsed-term-Jtilde}
\end{equation}
The estimate \eqref{eq:perron-kernel-estimate-final} is independent
of $a$ and $r$, and gives
\begin{equation*}
\left|K_{T'}(z)-\frac1z\right|
\ll \frac{\log T}{T'}
\end{equation*}
uniformly for $z=2\sigma_0+iy$ with $|y|\leq2T'$. The terms with
$a<r$ vanish by \eqref{eq:Jtilde-zero-revised}, so it remains only
to consider $a\geq r$. On the other hand, the degree bounds for $\tilde J_r^{(a)}(z;T)$ imply that the remaining part of the
$(r,a)$-integrand has $L^1$-norm
\begin{equation}
\label{eq:bound-result}
\int_{-2T'}^{2T'}
(\log T)^{k-a}|z|^{-a-k}
\left|\tilde J_r^{(a)}(z;T)\right|\,dy
\ll_k
(\log T)^{2k-1},
\qquad z=2\sigma_0+iy,
\end{equation}
uniformly in $r,a$. We arrive at the above bound by using
\begin{equation*}
\left|
\left(\frac{N}{T^r}\right)^z
\right|
=
\exp\left(2\sigma_0(\log N-r\log T)\right)
\ll_{\alpha,\eta,k}1,
\label{eq:Nr-uniform-bound}
\end{equation*}
since $\sigma_0=\eta/\log T$ and $\log N=\alpha\log T+o(\log T)$. For \eqref{eq:bound-result}, we also use
\eqref{eq:Jtilde-final-bound-revised} and $|z|\asymp \sigma_0+|y|$ to obtain, for $r=0$,
\begin{align*}
&(\log T)^{k-a}
|z|^{-a-k}
|\tilde J_0^{(a)}(z;T)|
\nonumber\\
&\qquad \qquad \qquad\ll_k
(\log T)^{k-a}
|z|^{-a-k}
\bigl(1+|z|\log T\bigr)^a
\ll_k
(\log T)^k(\sigma_0+|y|)^{-k}.
\label{eq:kernel-error-bound-r0}
\end{align*}
For $k\ge2$, this gives
\begin{equation*}
\int_{-2T'}^{2T'}
(\log T)^k(\sigma_0+|y|)^{-k}\,\dd y
\ll_k
(\log T)^k\sigma_0^{1-k}
\ll
(\log T)^{2k-1}.
\label{eq:L1-bound-r0}
\end{equation*}
For $k=1$, direct evaluation of the residues in
\eqref{eq:Jtilde-def-revised} gives
\begin{equation*}
\tilde J_0^{(0)}(z;T)=0,
\qquad
\tilde J_0^{(1)}(z;T)=1.
\end{equation*}
Consequently, for $a=0,1$ and $z=2\sigma_0+iy$,
\begin{align*}
&\int_{-2T'}^{2T'}
(\log T)^{1-a}|z|^{-a-1}
|\tilde J_0^{(a)}(z;T)|\,\dd y \\
&\qquad\ll
\int_{-2T'}^{2T'}
\frac{\dd y}{4\sigma_0^2+y^2}
\ll \sigma_0^{-1}
\ll \log T.
\end{align*}
Thus the $r=0$ contribution satisfies
\eqref{eq:bound-result} for every $k\ge1$.

Similarly, for $1\le r\le k$ and $a\ge r$,
\begin{align*}
&(\log T)^{k-a}
|z|^{-a-k}
|\tilde J_r^{(a)}(z;T)|
\nonumber\\
&\qquad\qquad\qquad\qquad\ll_k
(\log T)^{k-a}
|z|^{-a-k}
\bigl(1+|z|\log T\bigr)^{k+a-2r}
\nonumber\\
&\qquad\qquad\qquad\qquad\ll_k
(\log T)^{2k-2r}(\sigma_0+|y|)^{-2r},
\label{eq:kernel-error-bound-rpos}
\end{align*}
and this leads to
\begin{equation*}
\int_{-2T'}^{2T'}
(\log T)^{2k-2r}(\sigma_0+|y|)^{-2r}\,\dd y
\ll_k
(\log T)^{2k-2r}\sigma_0^{1-2r}
\ll
(\log T)^{2k-1}
\label{eq:L1-bound-rpos}
\end{equation*}
for $1\le r\le k$ which proves \eqref{eq:bound-result}. Therefore the contribution of the error term $K_{T'}(z)-1/z$ to each $(r,a)$-piece is
\begin{equation*}
\ll_k
\frac{(\log T)^{2k+1}}{T'}
=
o\bigl((\log T)^{2k}\bigr),
\end{equation*}
and the same bound holds for the total error across the sums over $r$ and $a$. In particular, this means $K_{T'}(z)$ may be replaced by $1/z$ in
\eqref{eq:collapsed-term-Jtilde} with an error $o((\log T)^{2k})$. Consequently,
\begin{align}
\label{eq:after-kernel-replacement}
\mathcal I_k(T,N)
&=
\frac{1}{k!^2}
\sum_{r=0}^{k}\binom{k}{r}^{2}
\sum_{a=r}^{k}
\binom{k}{a}
\frac{1}{2\pi i}
\int_{2\sigma_0-2iT'}^{2\sigma_0+2iT'}
(\log T)^{k-a}
\left(\frac{N}{T^r}\right)^z
\\
&\qquad\qquad\qquad\qquad \qquad \qquad\qquad\qquad\qquad\times
z^{-a-k-1}
\tilde J_r^{(a)}(z;T)\,\dd z
+
o\bigl((\log T)^{2k}\bigr). \nonumber
\end{align}
Here the lower limit $a=r$ comes from the vanishing condition
\eqref{eq:Jtilde-zero-revised}. We now evaluate the remaining finite $z$-integral by shifting the vertical segment according to the sign of $\log N-r\log T$. For a fixed
pair $(r,a)$, write
\begin{equation}
\lambda_{r,T}:=\log N-r\log T .
\end{equation}
Since $N=\lfloor T^\alpha\rfloor+\tfrac12$, we have
\begin{equation}
\lambda_{r,T}
=
(\alpha-r)\log T+o(\log T).
\label{eq:lambda-rT}
\end{equation}
As $\alpha\notin\mathbb Z$, the sign of $\lambda_{r,T}$ is
the sign of $\alpha-r$ for all sufficiently large $T$. For each $(r,a)$-piece, we have a $z$-integral of the form
\begin{equation}
\frac{1}{2\pi i}
\int_{2\sigma_0-2iT'}^{2\sigma_0+2iT'}
(\log T)^{k-a}
e^{\lambda_{r,T}z}
z^{-a-k-1}
\tilde J_r^{(a)}(z;T)\,\dd z .
\label{eq:finite-z-piece}
\end{equation}
As stated in \eqref{eq:Jtilde-expansion-revised}, $\tilde J_r^{(a)}(z;T)$ is a polynomial in
$z\log T$ that is truncated due to Lemma~\ref{lem:two-centre-degree-bound}. Thus the above integrand is meromorphic
in $z$, with its only pole at $z=0$.

We first consider the case $r>\alpha$, which corresponds to $\lambda_{r,T}<0$ for all
sufficiently large $T$. We shift the finite segment
\begin{equation*}
\Re z=2\sigma_0,
\qquad
|\Im z|\le 2T',
\end{equation*}
to the segment
\begin{equation*}
\Re z=A,
\qquad
|\Im z|\le 2T',
\end{equation*}
where $A>1$ is fixed. No pole is crossed and on the new vertical side,
\begin{equation}
|e^{\lambda_{r,T}z}|
=
e^{A\lambda_{r,T}}
\ll
T^{-A(r-\alpha)+o(1)}.
\label{eq:right-shift-exp-decay}
\end{equation}
Using \eqref{eq:Jtilde-final-bound-revised}, we have
\begin{equation}
(\log T)^{k-a}|z|^{-a-k-1}
\left|\tilde J_r^{(a)}(z;T)\right|
\ll_k
\begin{cases}
(\log T)^k |z|^{-k-1}, & r=0,\\[3pt]
(\log T)^{2k-2r}|z|^{-2r-1}, & 1\le r\le k.
\end{cases}
\label{eq:shift-contour-decay}
\end{equation}
On the line $\Re z=A$, the $z$-dependent factor in
\eqref{eq:shift-contour-decay} is absolutely integrable in $\Im z$,
uniformly for $|\Im z|\le 2T'$, and therefore the vertical side is bounded by
the corresponding power of $\log T$ times the decay factor
$T^{-A(r-\alpha)+o(1)}$, and is consequently $o((\log T)^{2k})$. The two
horizontal sides lie at height $|\Im z|=2T'$ and have bounded length
depending only on $A$. By \eqref{eq:shift-contour-decay}, their combined
contribution is
\begin{equation}
\ll_A
\begin{cases}
(\log T)^k(T')^{-k-1}, & r=0,\\[3pt]
(\log T)^{2k-2r}(T')^{-2r-1}, & 1\le r\le k,
\end{cases}
=
o((\log T)^{2k}).
\end{equation}
Therefore the whole $(r,a)$-piece is $o((\log T)^{2k})$ when
$r>\alpha$.

We now consider the case $r<\alpha$, for which we eventually have $\lambda_{r,T}>0$ and we shift the finite segment
\begin{equation*}
\Re z=2\sigma_0,
\qquad
|\Im z|\le 2T',
\end{equation*}
to the segment
\begin{equation*}
\Re z=-A,
\qquad
|\Im z|\le 2T',
\end{equation*}
with the same fixed $A>1$. In this deformation, the only pole crossed is
the pole at $z=0$. On the new vertical side,
\begin{equation}
|e^{\lambda_{r,T}z}|
=
e^{-A\lambda_{r,T}}
\ll
T^{-A(\alpha-r)+o(1)}.
\label{eq:left-shift-exp-decay}
\end{equation}
Together with \eqref{eq:shift-contour-decay}, this shows that the new
vertical side is $o((\log T)^{2k})$. The two horizontal sides are again
$o((\log T)^{2k})$ by the same estimate as above. 
Combining the two cases, we arrive at
\begin{align}
\label{eq:z-residue-selection-finite}
&\frac{1}{2\pi i}
\int_{2\sigma_0-2iT'}^{2\sigma_0+2iT'}
(\log T)^{k-a}
\left(\frac{N}{T^r}\right)^z
z^{-a-k-1}
\tilde J_r^{(a)}(z;T)\,\dd z
\\
&\qquad\qquad\qquad\qquad\qquad =
\mathds 1_{r<\alpha}
(\log T)^{k-a}
\operatorname*{Res}_{z=0}
\left[
\left(\frac{N}{T^r}\right)^z
z^{-a-k-1}
\tilde J_r^{(a)}(z;T)
\right]
+
o((\log T)^{2k}). \nonumber
\end{align}
In order to evaluate the $z$-residue, we combine the exponential $z$-dependence from $\tilde{J}^{(a)}_{r}(z;T)$ with the factor $(N/T^r)^{z}$ to get
\begin{align}
\operatorname*{Res}_{z=0}
\biggl[
\exp\biggl(
 z(\log N-r\log T)
-z\log T
\biggl\{
\sum_{j=1}^{r}(x_j-1)+\sum_{j=r+1}^{k}x_j
\biggr\}
\Biggr)
 z^{-a-k-1}
\biggr].
\label{eq:residue-before-simplification}
\end{align}
Due to \eqref{eq:pretty-obvious}, the above expression simplifies and we can evaluate the residue explicitly
\begin{align*}
&\operatorname*{Res}_{z=0}
\biggl[
\exp\biggl(
 z\biggl(\log N-(\log T)\sum_{j=1}^{k}x_j\biggr)
\biggr)
 z^{-a-k-1}
\biggr]=
\frac{1}{(a+k)!}
\biggl(\log N-(\log T)\sum_{j=1}^{k}x_j\biggr)^{a+k}
\label{eq:elementary-z-residue}
\end{align*}
which gives
\begin{align}
& (\log T)^{k-a}
\operatorname*{Res}_{z=0}
\Biggl[
\biggl(\frac{N}{T^r}\biggr)^z
z^{-a-k-1}
\tilde J_r^{(a)}(z;T)
\Biggr]
\\
&\quad =
\frac{(\log T)^{k-a}}{(a+k)!(2\pi i)^{2k}}
\int_{c_j(r)}\int_{\hat c_\ell(r)}
\biggl(\log N-(\log T)\sum_{j=1}^{k}x_j\biggr)^{a+k}
\Delta(x)^2\Delta(y)^2
\prod_{j,\ell=1}^{k}(x_j-y_\ell-1)
\nonumber\\
&\qquad\qquad\qquad\times
\biggl(
\sum_{j=1}^{k}\frac1{x_j}
+
\sum_{\ell=1}^{k}\frac1{y_\ell+1}
\biggr)^a
\biggl(
\sum_{j=1}^{k}\frac1{x_j-1}
+
\sum_{\ell=1}^{k}\frac1{y_\ell}
\biggr)^k
\prod_{j=1}^{k}x_j^{-k}\,dx_j
\prod_{\ell=1}^{k}y_\ell^{-k}\,dy_\ell . \nonumber
\end{align}
As $\frac{\log N}{\log T}=\alpha+o(1)$ and the $x$-contours are fixed compact contours, we have 
\begin{align}
(\log T)^{k-a}
\Bigl(\log N-(\log T)\sum_{j=1}^{k}x_j\Bigr)^{a+k}
= (\log T)^{2k}
\Bigl(\alpha-\sum_{j=1}^{k}x_j\Bigr)^{a+k}
+ o\bigl((\log T)^{2k}\bigr)
\label{eq:logN-to-alpha-uniform}
\end{align}
uniformly for $x$ on these contours. Combining
\eqref{eq:z-residue-selection-finite}--\eqref{eq:logN-to-alpha-uniform}
with the preceding estimates gives
\begin{equation}
\mathcal I_k(T,N)
=
(\log T)^{2k}\mathcal M_k(\alpha)
+o\bigl((\log T)^{2k}\bigr),
\label{eq:Ik-final-asymptotic-before-limit}
\end{equation}
where
\begin{align*}
\mathcal M_k(\alpha)
&=
\frac{1}{k!^2}
\sum_{0 \leq r<\alpha}
\binom{k}{r}^{2}
\sum_{a=r}^{k}
\binom{k}{a}
\frac{1}{(a+k)!}
\nonumber\\
&\quad\times \frac{1}{(2\pi i)^{2k}}
\int_{c_j(r)}\int_{\hat c_\ell(r)}
\left(\alpha-\sum_{j=1}^{k}x_j\right)^{a+k}
\Delta(x)^2\Delta(y)^2
\prod_{j,\ell=1}^{k}(x_j-y_\ell-1)
\nonumber\\
&\quad\times
\left(
\sum_{j=1}^{k}\frac1{x_j}
+
\sum_{\ell=1}^{k}\frac1{y_\ell+1}
\right)^a
\left(
\sum_{j=1}^{k}\frac1{x_j-1}
+
\sum_{\ell=1}^{k}\frac1{y_\ell}
\right)^k
\prod_{j=1}^{k}x_j^{-k}\,dx_j
\prod_{\ell=1}^{k}y_\ell^{-k}\,dy_\ell .
\label{eq:Mkalpha-final-form}
\end{align*}
Here the condition $r<\alpha$ is precisely the selection rule obtained
from the finite contour shift in \eqref{eq:z-residue-selection-finite}; no
term with $r>\alpha$ contributes at the leading scale. Since
$\alpha\notin\mathbb Z$, there is no boundary case $r=\alpha$.
As $N=\lfloor T^\alpha\rfloor+\tfrac12$ differs from $T^\alpha$ by
$O(1)$, this is the leading asymptotic for the sharp cutoff at
height $T^\alpha$. Extracting a $(-1)^{k+a}$ factor from $(\alpha - \sum_j x_j)^{a+k}$ gives the form stated in \eqref{eq:prk}.
\paragraph{Integer values for $\alpha.$} It remains to discuss the case in which $\alpha$ is an integer. Let
$\alpha=m\in\mathbb Z_{>0}$. For $r\ne m$, the preceding contour
shift is unchanged, since
\begin{equation*}
\log N-r\log T=(m-r)\log T+o(\log T).
\end{equation*}
Thus the terms with $r<m$ contribute their residues at $z=0$, while the terms with $r>m$ are negligible. We call the terms with $r=m$ and $m\leq a\leq k$ the transition
block. If $m>k$, this block is empty and the formula follows
immediately; if $1\leq m\leq k$, we still need to treat this block. For this block, using $N=\lfloor T^m\rfloor+\tfrac12=T^m+O(1)$, we have 
\begin{equation*}
\log{N} = \log{(T^{m} + O(1)) = m\log{T} + \log{(1 + O(T^{-m}))}},
\end{equation*}
which gives
\begin{equation}
\label{eq:not-log}
\log N-m\log T=O(T^{-m}),
\end{equation}
using $\log{(1 + \varepsilon)} = O(\varepsilon)$. Hence the exponential
factor $(N/T^m)^z$ no longer gives decay in either horizontal direction as \eqref{eq:not-log} is not of order $\log{T}$.
Nevertheless this block is negligible at the scale $(\log T)^{2k}$.
Explicitly, after replacing $K_{T'}(z)$ by $1/z$, the $(m,a)$-piece is
\begin{equation*}
\frac{1}{2\pi i}
\int_{2\sigma_0-2iT'}^{2\sigma_0+2iT'}
(\log T)^{k-a}
\left(\frac{N}{T^m}\right)^z
z^{-a-k-1}
\tilde J_m^{(a)}(z;T)\,dz .
\end{equation*}
By \eqref{eq:shift-contour-decay},
\begin{equation*}
(\log T)^{k-a}|z|^{-a-k-1}
|\tilde J_m^{(a)}(z;T)|
\ll_k
(\log T)^{2k-2m}|z|^{-2m-1}.
\end{equation*}
Choose $A_T=T^\theta$, where $0<\theta<m$, and shift the finite segment
to the line $\Re z=A_T$. No pole is crossed, and since
\begin{equation*}
\left|
\left(\frac{N}{T^m}\right)^z
\right|
=
\exp\{O(T^{-m}A_T)\}
=
1+o(1),
\end{equation*}
the new vertical side is
\begin{equation*}
\ll
(\log T)^{2k-2m}
\int_{-2T'}^{2T'}(A_T+|y|)^{-2m-1}\,dy
\ll
(\log T)^{2k-2m}A_T^{-2m}
=
o((\log T)^{2k}).
\end{equation*}
The two horizontal sides are also negligible, since their total contribution is
\begin{equation*}
\ll
(\log T)^{2k-2m}A_T(T')^{-2m-1}
=
o((\log T)^{2k}).
\end{equation*}
Thus the transition block $r=m$ contributes $o((\log T)^{2k})$. Consequently the same formula for $\mathcal M_k(\alpha)$ remains valid at $\alpha=m$, with the convention that only terms with $r<m$ are included. Finally, by Proposition~\ref{cor:zero-order}, the zero-order of $P_{r,k}(\alpha)$ at $\alpha=r$ is at least $2r^2$. Thus, $\mathcal{M}_{k}(\alpha)$ is continuously differentiable $2r^2 -1$ times at $\alpha = r$ for $r=1, \ldots, k$.
\subsection{Proof of Proposition~\ref{eq:final-corollary}}
Throughout, we assume the hypotheses of Conjecture~\ref{eq:main-conj}; by its application, the normalised moments in \eqref{eq:first-corollary-statement} are equal to the conjectural main term
carrying an error term
\begin{equation*}
O_{k,c,\xi}\!\left(
(\log T)^{2k-1+\xi}
\bigl(1+\min(|z|,1)\log T\bigr)^{-\theta_{k,c}}
\right),
\end{equation*}
with $\theta_{k,c}$ as defined in \eqref{eq:theta_kc}. Following the two-centre decomposition, together with \eqref{eq:J-factor-revised}, the main term derived in the proof of Theorem~\ref{mkalpha} is a finite linear combination of terms of the form
\begin{equation*}
(\log T)^{k-a}z^{-a-k}J_r^{(a)}(-z;T)
=
(\log T)^{k-a}z^{-a-k}T^{-rz}\tilde J_r^{(a)}(z;T),
\end{equation*}
where $0\le a, r\le k$ and $0\le a\le k$. Here $J^{(a)}_r$ and $\tilde{J}^{(a)}_r$ are both as defined in \eqref{def:J} and \eqref{eq:Jtilde-def-revised} respectively. Since $\Re z>0$, we have $|T^{-rz}|\le1$. By Lemma~\ref{lem:two-centre-degree-bound} and through the steps outlined in \eqref{eq:Jtilde-zero-revised}-\eqref{eq:Jtilde-final-bound-revised}, we showed
\begin{equation*}
\tilde J_r^{(a)}(z;T)=0
\qquad (a<r),
\end{equation*}
and, for $a\ge r$,
\begin{equation}
\tilde J_r^{(a)}(z;T)
=
\begin{cases}
O_k\bigl((1+|z|\log T)^a\bigr), & r=0,\\[3pt]
O_k\bigl((1+|z|\log T)^{k+a-2r}\bigr), & 1\le r\le k.
\end{cases}
\label{eq:Jtilde-repeat}
\end{equation}
It suffices to consider $a\ge r$. Note that we make repeated use of the following two facts in the subsequent estimates:
\begin{equation}
\label{eq:lower-range-bound-minz1}
1+\min(|z|,1)\log T\asymp1
\qquad |z|\log T\leq1,
\end{equation}
and 
\begin{equation}
\label{eq:second-range-bound}
1+\frac{1}{|z|\log T}\asymp1
\qquad |z|\log T>1.
\end{equation}
Suppose $1\le r\le k$ for which, by the above, we have
\begin{equation}
(\log T)^{k-a}|z|^{-a-k}
\left|\tilde J_r^{(a)}(z;T)\right|
\ll_k
(\log T)^{k-a}|z|^{-a-k}
\bigl(1+|z|\log T\bigr)^{k+a-2r}.
\end{equation}
If $|z|\log T\le1$, then this is
\begin{equation*}
\ll_k
(\log T)^{2k}
\left(1+\frac{1}{|z|\log T}\right)^{2k},
\end{equation*}
which follows from rewriting $|z|^{-a-k} = (\log{T})^{a+k}(|z|\log{T})^{-a-k}$ and using $1+|z|\log T\asymp1$. This is bounded by the right-hand side of
\eqref{eq:first-corollary-statement}. In the case $|z|\log T>1$, then it is
\begin{equation*}
\ll_k
(\log T)^{2k}(|z|\log T)^{-2r},
\end{equation*}
using $1+|z|\log T\asymp |z|\log T$.
Since $r\ge1$, this is
\begin{equation*}
\ll_k
(\log T)^{2k}
\bigl(1+\min(|z|,1)\log T\bigr)^{-2},
\end{equation*}
and hence is again bounded by the right-hand side of
\eqref{eq:first-corollary-statement}.

It remains to handle $r=0$. By \eqref{eq:Jtilde-repeat},
\begin{equation*}
(\log T)^{k-a}|z|^{-a-k}
\left|\tilde J_0^{(a)}(z;T)\right|
\ll_k
(\log T)^{k-a}|z|^{-a-k}
\bigl(1+|z|\log T\bigr)^a.
\end{equation*}
If $|z|\log T\le1$, we obtain the required contribution in 
\eqref{eq:first-corollary-statement} immediately. If $|z|\log T>1$, then it is
\begin{equation*}
\ll_k
(\log T)^{2k}(|z|\log T)^{-k}.
\end{equation*}
For $k\ge2$, this is
\begin{equation*}
\ll_k
(\log T)^{2k}
\bigl(1+\min(|z|,1)\log T\bigr)^{-2}.
\end{equation*}
The only remaining case is when $k=1$ and $r=0$. In this case, direct evaluation from \eqref{eq:Jtilde-def-revised} gives
\begin{equation*}
\tilde J_0^{(0)}(z;T)=0,
\qquad
\tilde J_0^{(1)}(z;T)=1.
\end{equation*}
Thus the $r=0$ contribution is $O(|z|^{-2})$ for $k=1$. If $|z|\log T\le1$, then
\begin{equation*}
|z|^{-2}
\ll
(\log T)^2
\left(1+\frac{1}{|z|\log T}\right)^2,
\end{equation*}
while if $|z|\log T>1$, then
\begin{equation*}
|z|^{-2}
\ll
(\log T)^2
\bigl(1+\min(|z|,1)\log T\bigr)^{-2}.
\end{equation*}
Hence the $k=1,r=0$ contribution also satisfies
\eqref{eq:first-corollary-statement}. Therefore, the conjectural main term satisfies \eqref{eq:first-corollary-statement} for every $k\ge1$; we now bound the error term. If $|z|\log T\le1$, then \eqref{eq:lower-range-bound-minz1} holds and the error term is
\begin{equation*}
\ll_{k,c,\xi}(\log T)^{2k-1+\xi}.
\end{equation*}
As $\xi<1$, this is $o((\log T)^{2k})$, and is bounded by
\eqref{eq:first-corollary-statement}. Suppose now that
$|z|\log T>1$. Then the error term is
\begin{equation*}
(\log T)^{2k}
\bigl(1+\min(|z|,1)\log T\bigr)^{-2}
\cdot
(\log T)^{-1+\xi}
\bigl(1+\min(|z|,1)\log T\bigr)^{2-\theta_{k,c}}.
\end{equation*}
If $k\ge3$, then $\theta_{k,c}=k-1\ge2$, so the final factor is $O(1)$ because $0<\xi<1$. If $k=1,2$, then $\theta_{k,c}=1+c$, and the final
factor is
\begin{equation*}
\ll
(\log T)^{-1+\xi}(\log T)^{1-c}
=
(\log T)^{\xi-c}
\ll1,
\end{equation*}
because $0<\xi<c$. Combining the bounds for the main term and the error term proves
\eqref{eq:first-corollary-statement}. Finally, on any sub-range satisfying
$|z|\log T\ge c_0>0$, we have
\begin{equation*}
\left(1+\frac{1}{|z|\log T}\right)^{2k}
\ll_{k,c_0}1,
\end{equation*}
and \eqref{eq:second-corollary-statement} follows. 

\subsection{Proof of Proposition~\ref{cor:zero-order}}
To determine the zero-order of $P_{r,k}(\alpha)$ at $r=\alpha$, we count local degree with respect to the shifted variables, $x_j-1$ for $1\le j\le r$ and $y_\ell+1$ for $1\le \ell\le r$, as well as the unshifted variables, $x_j$ for $r<j\le k$ and $y_\ell$ for $r<\ell\le k$.
Since
\begin{equation*}
\sum_{j=1}^k x_j-\alpha
=
\sum_{j=1}^r(x_j-1)
+
\sum_{j=r+1}^k x_j
-
(\alpha-r),
\end{equation*}
we have
\begin{align*}
\Bigl( \, \sum_{j=1}^k x_j-\alpha \, \Bigr)^{k+m}
=
\sum_{d=0}^{k+m}
\binom{k+m}{d}
\Bigl(
\sum_{j=1}^r(x_j-1)
+
\sum_{j=r+1}^k x_j
\Bigr)^d
(-1)^{k+m-d}(\alpha-r)^{k+m-d}.
\end{align*}
We claim that the residue corresponding to the degree-$d$ term vanishes whenever $d>k+m-2r^2$. To prove this, first consider the factor
\begin{equation*}
\Delta(x)^2\Delta(y)^2
\prod_{i,j=1}^k(x_i-y_j-1)
\prod_{j=1}^kx_j^{-k}
\prod_{\ell=1}^k y_\ell^{-k},
\end{equation*}
from the integrand in \eqref{eq:prk}.
In all $2k$ variables, its lowest possible local total degree is
\begin{align*}
&2r(r-1)+2(k-r)(k-r-1)+2r(k-r)-2k(k-r)
 = 2r^2-2k,
\end{align*}
where the contribution $2r(k-r)$ arises from the mixed factors $x_i-y_j-1$ with $1\le i\le r, r<j\le k$ or $r<i\le k, 1\le j\le r$. Next,
\begin{equation*}
\biggl(
\sum_{j=1}^{k}\frac{1}{x_j}
+
\sum_{\ell=1}^{k}\frac{1}{y_\ell+1}
\biggr)^m \, \biggl(
\sum_{j=1}^{k}\frac{1}{x_j-1}
+
\sum_{\ell=1}^{k}\frac{1}{y_\ell}
\biggr)^k
\end{equation*}
has total pole order at most $m+k$ in these local variables. Therefore, after multiplying by the
degree-$d$ term from
\begin{equation*}
\biggl(\sum_{j=1}^k x_j-\alpha\biggr)^{k+m},
\end{equation*}
the local total degree of the full integrand is at least
\begin{align*}
d+(2r^2-2k)-m-k
=
d+2r^2-3k-m.
\end{align*}
A non-zero residue results only from monomials of the form
\begin{equation*}
\prod_{j=1}^{r}(x_j-1)^{-1}
\prod_{j=r+1}^{k}x_j^{-1}
\prod_{\ell=1}^{r}(y_\ell+1)^{-1}
\prod_{\ell=r+1}^{k}y_\ell^{-1},
\end{equation*}
which has total degree $-2k$. Hence the residue vanishes whenever
\begin{equation*}
d+2r^2-3k-m>-2k,
\end{equation*}
or equivalently whenever $d>k+m-2r^2$. The coefficient of $(\alpha-r)^q$ in the $m$th summand in \eqref{eq:prk} corresponds to
$d=k+m-q$. If $q<2r^2$, then
\begin{equation*}
d=k+m-q>k+m-2r^2,
\end{equation*}
so that coefficient vanishes. Thus every coefficient of $(\alpha-r)^q$
with $q<2r^2$ is zero. Hence each summand in the definition of
$P_{r,k}(\alpha)$, and therefore $P_{r,k}(\alpha)$ itself, is divisible by
$(\alpha-r)^{2r^2}$. Equivalently,
\begin{equation*}
P_{r,k}(\alpha)=O\left((\alpha-r)^{2r^2}\right)
\qquad(\alpha\to r),
\end{equation*}
as required.
\subsection{Proof of Theorem ~\ref{conj:derivative}}
 Let integer $k\geq 1$ and $H=X^{1-1/\alpha}$ for a fixed $\alpha>0$. For non-negative $g\in C_c^\infty((0,\infty))$ satisfying $\operatorname{supp} g\subset[1/2,4]$, define the smoothed variance
\begin{equation}\label{eq:smoothed-variance}
\tilde V_{g,k}(X,H)
:=
\frac{1}{HX}\int_{\mathbb{R}} g\!\left(\frac{x}{X}\right)\, \bigl|\Delta_k(x;H)\bigr|^2\,dx.
\end{equation}
The support condition localises $x$ to a fixed range that ensures $x\asymp X$ and avoids edge effects. The passage from the sharp average over $x\in[X,2X]$ to this smoothed average is justified by approximating $\mathds{1}_{[1,2]}(u)$ from above and below by non-negative smooth functions supported in $[\tfrac12,4]$. Thus it suffices to prove
\begin{align}
\label{eq:final-statement-to-prove}
    \lim_{X \rightarrow \infty}
    \frac{\tilde V_{g,k}(X,X^{1-1/\alpha})}
    {\bigl(\log(X^{1/\alpha})\bigr)^{2k-1}}
    =
    \int_{\mathbb{R}} g(u)\,du\,
    \frac{1}{k!^2}
    \sum_{0 \leq r < \alpha}
    {k \choose r}^{2} P'_{r,k}(\alpha),
\end{align}
for every fixed non-negative $g\in C_c^\infty((0,\infty))$ with $\operatorname{supp}g\subset[\tfrac12,4]$. 
We first use Lemma~\ref{lemma:titchmarsh} (Lemma 3.12 from \cite{kn:titchmarsh}) to show that a truncated Perron integral on the line $\Re s=c$ where $c= 1 + \frac{1}{\log{X}}$
approximates $\Psi_k(x;H)$ defined in \eqref{eq:psi} with a negligible mean-square error. Note that this choice of $c$ places the Perron contour just inside the region of absolute convergence for \eqref{eq:dirichlet-series} while keeping $X^{c -1}=O(1)$ and $-\zeta'/\zeta(c )\ll \log X$ by \eqref{eq:unfortunatelysuperimportant}. This ensures that the truncation errors remain only logarithmic in size. Define the truncated Perron kernel
\begin{equation*}
        K_T(u):=\frac1{2\pi i}\int_{c-iT}^{c+iT}u^s\frac{ds}{s}.
\end{equation*}
Then
\begin{equation}
\label{eq:full-variance-trunc-approx}
\begin{aligned}
&\frac1{2\pi i}\int_{c-iT}^{c+iT}
\left(-\Ldz(s)\right)^k
\frac{(x+H)^s-x^s}{s}\,ds                         =
\sum_{n=1}^{\infty}\Lambda^{*k}(n)
\left[
K_T\!\left(\frac{x+H}{n}\right)
-
K_T\!\left(\frac{x}{n}\right)
\right].
\end{aligned}
\end{equation}
By Lemma~\ref{lemma:titchmarsh}, equivalently by the usual truncated Perron
kernel estimate,
\begin{equation}
\label{eq:K-estimate}
        K_T(u)
        =
       \mathds{1}_{u>1}
        +
        O\left(
        u^c \min\left\{1,\frac1{T|\log u|}\right\}
        \right).
\end{equation}
This means
\begin{equation*}
        K_T\!\left(\frac{x+H}{n}\right)
        -
        K_T\!\left(\frac{x}{n}\right)
\end{equation*}
approximates the indicator of $x<n\le x+H$. Note that this is a direct
short-interval use of Perron's formula; we do not apply Perron's formula
separately at $x$ and at $x+H$. As the estimate $\eqref{eq:K-estimate}$ measures the error made by replacing the sharp cutoff $\mathds{1}_{u>1}$ with the truncated Perron kernel $K_T(u)$, we choose height $T$ large enough so that the averaged pointwise error contribution to $\tilde V_{g,k}(X,H)$ is negligible.  We define the height
\begin{equation*}
        T=\frac{X^{1+\epsilon}}{\sqrt H},
\end{equation*}
and choose $\epsilon>0$ so that, after division by $HX$, the
mean-square Perron truncation error is
$o\!\left((\log(X^{1/\alpha}))^{2k-1}\right)$, as required in
\eqref{eq:final-statement-to-prove}. More precisely, since $H=X^{1-1/\alpha}$, the condition $\frac{X}{TH}=o(1)$ is equivalent to
\begin{equation*}
        X^{\frac{1}{2\alpha}-\frac12-\epsilon}=o(1).
\end{equation*}
Thus it is enough to assume $\epsilon>
\max\left\{0,\frac{1}{2\alpha}-\frac12\right\}$.  
The only range in which the estimate \eqref{eq:K-estimate} requires care is the range where one of the arguments of $K_T$ may be close to the transition point $u=1$. On the support of $g(x/X)$, we have $x \in [X/2, 4X]$ which means  $x\asymp X$. Since $H= o(X)$, the interval $(x, x+H]$ lies in a constant multiple of $X$. Hence every integer $n$ that satisfies
\begin{equation*}
        x<n\le x+H,
\end{equation*}
must also satisfy
\begin{equation*}
    AX\le n\le BX
\end{equation*}
for suitable fixed constants
$A,B>0$, depending only on $\operatorname{supp} g$. We call $n\in[AX,BX]$ the near range and its complement the outer range.
Subtracting \eqref{eq:psi} from the truncated Perron expression in
\eqref{eq:full-variance-trunc-approx} and splitting the resulting sum
into the near and outer ranges gives
\begin{equation*}
\frac1{2\pi i}\int_{c-iT}^{c+iT}
\left(-\Ldz(s)\right)^k
\frac{(x+H)^s-x^s}{s}\,ds
-\Psi_k(x;H)
=
\mathcal E_{\rm near}(x)+\mathcal E_{\rm out}(x),
\end{equation*}
where
\begin{align}
\label{eq:near-range-E}
\mathcal E_{\rm near}(x)
&:=
\sum_{n\in[AX,BX]}\Lambda^{*k}(n)
\left[
K_T\!\left(\frac{x+H}{n}\right)
-
K_T\!\left(\frac{x}{n}\right)
-
\mathds{1}_{\{x<n\leq x+H\}}
\right],
\end{align}
and
\begin{align}
\label{eq:outer-range-E}
\mathcal E_{\rm out}(x)
:=
\sum_{n\notin[AX,BX]}\Lambda^{*k}(n)
\left[
K_T\!\left(\frac{x+H}{n}\right)
-
K_T\!\left(\frac{x}{n}\right)
\right].
\end{align}
In the outer range, the indicator $\mathds{1}_{x<n\le x+H}$ is zero, and the
ratios $x/n$ and $(x+H)/n$ are bounded away from $1$. Hence the Perron
kernel estimate gives
\begin{equation*}
\left|
K_T\!\left(\frac{x+H}{n}\right)
-
K_T\!\left(\frac{x}{n}\right)
\right|
\ll
\frac1T
\left[
\left(\frac{x+H}{n}\right)^c
+
\left(\frac{x}{n}\right)^c
\right]
\ll
\frac{X^c}{Tn^c},
\end{equation*}
using \eqref{eq:K-estimate} for the outer range. By the triangle inequality and the preceding estimate,
\begin{align*}
|\mathcal E_{\rm out}(x)|
&\leq
\sum_{n\notin[AX,BX]}
\Lambda^{*k}(n)
\left|
K_T\!\left(\frac{x+H}{n}\right)
-
K_T\!\left(\frac{x}{n}\right)
\right| \\
&\ll
\frac{X^c}{T}
\sum_{n\notin[AX,BX]}
\frac{\Lambda^{*k}(n)}{n^c} \\
&\leq
\frac{X^c}{T}
\sum_{n=1}^{\infty}
\frac{\Lambda^{*k}(n)}{n^c}.
\end{align*}
Since $\Lambda^{*k}(n)\ge0$ and $c=1+1/\log X$,
\begin{equation*}
        \sum_{n=1}^{\infty}\frac{\Lambda^{*k}(n)}{n^c}
        =
        \left(\sum_{n=1}^{\infty}\frac{\Lambda(n)}{n^c}\right)^k
        =
        \left(-\Ldz(c)\right)^k
        \ll_k
        (\log X)^k,
\end{equation*}
as in \eqref{eq:bounding-series}. The approximation in \eqref{eq:full-variance-trunc-approx} will be
inserted into the weighted mean square
\eqref{eq:smoothed-variance}. For $\delta := X/T$, the outer-range error and its contribution to \eqref{eq:smoothed-variance} are bounded by
\begin{equation*}
        \mathcal E_{\rm out}(x)\ll_k \delta(\log X)^k,
        \qquad
        \int_{\mathbb R}g\!\left(\frac{x}{X}\right)
        |\mathcal E_{\rm out}(x)|^2\,dx
        \ll_k X\delta^2(\log X)^{2k}.
\end{equation*}
It remains to handle the near range $n\asymp X$. There, the Perron kernel
estimate gives
\begin{equation}
\label{eq:full-perron-kernel-estimate}
\left|
K_T\!\left(\frac{x+H}{n}\right)
-
K_T\!\left(\frac{x}{n}\right)
-
{\mathds 1}_{x<n\le x+H}
\right| \ll
\min\left\{1,\frac{\delta}{|x-n|}\right\}
+
\min\left\{1,\frac{\delta}{|x+H-n|}\right\}.
\end{equation}
The near-range error arises from the truncated Perron kernel being least accurate near the jump points of  the approximations for the two indicators
${\mathds 1}_{n\le x}$ and ${\mathds 1}_{n\le x+H}$. Thus the error is controlled
by two boundary kernels, one measuring the distance from $x=n$ and the other
from $x=n-H$. Set
\begin{equation*}
\mathcal K_\delta(y):=
\begin{cases}
\displaystyle \min\left\{1,\frac{\delta}{|y|}\right\},
    & y\neq 0,\\[1mm]
1, & y=0.
\end{cases}
\end{equation*}
Using \eqref{eq:full-perron-kernel-estimate}, $\mathcal{E}_{\rm near}(x)$ as defined in \eqref{eq:near-range-E} and the triangle inequality, we have
\begin{align*}
|\mathcal E_{\rm near}(x)|
&\ll
\sum_{n\asymp X}\Lambda^{*k}(n)\mathcal K_\delta(x-n) +
\sum_{n\asymp X}\Lambda^{*k}(n)\mathcal K_\delta(x+H-n).
\end{align*}
These are the two boundary contributions, centred respectively at
$x=n$ and $x=n-H$. It is enough to estimate the first, since the
second is handled identically after translating $x$ by $H$. Uniformly for $n\asymp X$ and $x\asymp X$, respectively, we have
\begin{equation*}
\int_{x\asymp X}\mathcal K_\delta(x-n)\,dx
\ll\delta\log X,
\qquad
\sum_{n\asymp X}\mathcal K_\delta(x-n)
\ll(1+\delta)\log X.
\end{equation*}
By \eqref{eq:pointwise-dk}, together with the standard estimate $\sum_{n\le Y} d_k(n)^2\ll_k Y(\log Y)^{k^2-1}$ as stated in \cite[\S1.6, (1.80)]{iwaniec-kowalski}, we obtain
\begin{equation}
\label{eq:secondmoment-bound}
        \sum_{n\le Y}|\Lambda^{*k}(n)|^2
        \ll_k
        Y(\log Y)^{B_k}
\end{equation}
for some constant $B_k>0$. By the Schur test
\cite[Theorem~5.2, p.~22]{kn:halmos-sunder}, applied with counting
measure in $n$ and Lebesgue measure in $x$, we obtain
\begin{align*}
&\int_{\mathbb R}g\!\left(\frac{x}{X}\right)
\left|
\sum_{n\asymp X}\Lambda^{*k}(n)\mathcal K_\delta(x-n)
\right|^2dx \\
&\qquad\qquad\qquad\qquad\ll_g
\left(
\sup_{x\asymp X}
\sum_{n\asymp X}\mathcal K_\delta(x-n)
\right)
\left(
\sup_{n\asymp X}
\int_{x\asymp X}\mathcal K_\delta(x-n)\,dx
\right) 
\sum_{n\asymp X}|\Lambda^{*k}(n)|^2                \\
&\qquad\qquad\qquad\qquad\ll
(1+\delta)\log X \cdot \delta\log X
\sum_{n\asymp X}|\Lambda^{*k}(n)|^2                                      \\
&\qquad\qquad\qquad\qquad\ll_k
X(\delta+\delta^2)(\log X)^{B_k},
\end{align*}
where in the last line we applied \eqref{eq:secondmoment-bound} and also enlarged
$B_k$. Similarly,
\begin{align*}
\int_{\mathbb R}g\!\left(\frac{x}{X}\right)
\left|
\sum_{n\asymp X}\Lambda^{*k}(n)\mathcal K_\delta(x+H-n)
\right|^2dx
&\ll_k
X(\delta+\delta^2)(\log X)^{B_k}.
\end{align*}
Collectively, this leads to
\begin{equation*}
\int_{\mathbb R}g\!\left(\frac{x}{X}\right)
\left|\Psi_k(x;H)-\frac1{2\pi i}\int_{c-iT}^{c+iT}
\left(-\Ldz(s)\right)^k
\frac{(x+H)^s-x^s}{s}\,ds    \right|^2dx
\ll_k
X(\delta+\delta^2)(\log X)^{B_k}.
\end{equation*}
Dividing through by $HX$, and using $\delta=\frac XT=\sqrt H\,X^{-\epsilon}$, 
we get
\begin{equation*}
\frac1{HX}\int_{\mathbb R}g\!\left(\frac{x}{X}\right)
\left|\Psi_k(x;H)-\frac1{2\pi i}\int_{c-iT}^{c+iT}
\left(-\Ldz(s)\right)^k
\frac{(x+H)^s-x^s}{s}\,ds  \right|^2dx
=
o\!\left((\log X)^{2k-1}\right).
\end{equation*}
We next shift from the line $\Re s=c$ to $        \Re s=\sigma_*=\frac12+\frac1{\log T}$.
Assuming RH, the poles of $-\zeta'/\zeta(s)$ coming from non-trivial zeros of
$\zeta(s)$ lie on $\Re s=1/2$, and are therefore not crossed. The only pole
crossed is at $s=1$, giving
\begin{align}
\label{eq:res-p}
\frac1{2\pi i}\int_{c-iT}^{c+iT}
\left(-\Ldz(s)\right)^k
\frac{(x+H)^s-x^s}{s}\,ds
        =
        &\operatorname*{Res}_{s=1}
        \left[
        \left(-\Ldz(s)\right)^k
        \frac{(x+H)^s-x^s}{s}
        \right] \\
        &\qquad \qquad \qquad \qquad \qquad +
        \mathcal P_k(x;H)
        +
        R_{\rm hor}(x), \nonumber
\end{align}
where
\begin{equation}
\label{eq:P-should-have-a-ref}
        \mathcal P_k(x;H)
        :=
        \frac1{2\pi}
        \int_{-T}^{T}
        \left(-\Ldz(\sigma_*+it)\right)^k
        \frac{(x+H)^{\sigma_*+it}-x^{\sigma_*+it}}
             {\sigma_*+it}\,dt,
\end{equation}
and
\begin{align}
\label{eq:Rhor-definition}
R_{\rm hor}(x)
:=
\frac1{2\pi i}
\left( \,
\int_{\sigma_*+iT}^{c+iT}
+
\int_{c-iT}^{\sigma_*-iT}
\, \right)
\left(-\Ldz(s)\right)^k
\frac{(x+H)^s-x^s}{s}\,ds .
\end{align}
We now bound the horizontal contribution $R_{\rm hor}$. On the
horizontal sides $s=\sigma\pm iT$, with $\sigma_*\le\sigma\le c$, the
standard RH bound in \eqref{eq:titchmarsh} gives $\Ldz(s)\ll(\log T)^2$. For $x \asymp X$, we have
\begin{equation}
\label{perron-diff-kernel}
\left|
\frac{(x+H)^s-x^s}{s}
\right|
\ll
\min\left\{
        H X^{\sigma-1},
        \frac{X^\sigma}{T}
\right\}.
\end{equation}
The first bound follows from $\frac{(x+H)^s-x^s}{s}=\int_x^{x+H}u^{s-1}\,du$ while the second follows from the triangle inequality and $|s|\asymp T$. Since $X/(TH)=o(1)$, the second term is the smaller one for all sufficiently
large $X$. This means that each horizontal side is
\begin{equation*}
        \ll
        \int_{\sigma_*}^{c}
        (\log T)^{2k}\frac{X^\sigma}{T}\,d\sigma
        \ll
        \frac{X}{T}(\log T)^{2k+1}.
\end{equation*}
Therefore
\begin{equation*}
\frac1{HX}\int_{\mathbb R}g\!\left(\frac{x}{X}\right)
|R_{\rm hor}(x)|^2dx
\ll
\frac1H\left(\frac XT\right)^2(\log T)^{4k+2}
=
X^{-2\epsilon}(\log T)^{4k+2}
=
o\!\left((\log X)^{2k-1}\right).
\end{equation*}
The remainder $R(x)$ combines the horizontal-contour error $R_{\text{hor}}(x)$ with the difference between $\Psi_k$ and the truncated Perron integral. From the
definition of $\Delta_{k}(x;H)$ and using \eqref{eq:res-p}, we obtain
\begin{equation}
\label{eq:back-to-reality}
        \Delta_{k}(x;H)
        =
        \mathcal P_k(x;H)+R(x),
\end{equation}
where
\begin{equation}
\label{eq:remainder-small}
        \frac1{HX}\int_{\mathbb R}g\!\left(\frac{x}{X}\right)
        |R(x)|^2dx
        =
        o\!\left((\log X)^{2k-1}\right).
\end{equation}
To pass from $\Delta_{k}$ to $\mathcal P_k$ inside the variance, we use
\begin{equation*}
        |\Delta_{k}|^2
        =
        |\mathcal P_k|^2
        +
        2\Re(\mathcal P_k\overline R)
        +
        |R|^2.
\end{equation*}
Due to \eqref{eq:remainder-small}, the term involving $|R|^2$ is negligible.  For
the cross-term, we use the Cauchy-Schwarz inequality
\begin{align}
\label{eq:Cauchy-schwarz}
&\biggl| \frac1{HX}\int_{\mathbb R}g\!\left(\frac{x}{X}\right)
\mathcal P_k(x;H)\overline{R(x)}\,dx \biggr| \\ &\qquad \qquad \qquad \quad \leq \biggl(         \frac1{HX}\int_{\mathbb R}g\!\left(\frac{x}{X}\right)
        |\mathcal P_k(x;H)|^2dx \biggr)^{1/2} \biggl(  \frac1{HX}\int_{\mathbb R}g\!\left(\frac{x}{X}\right)
        |R(x)|^2dx\biggr)^{1/2}\nonumber .
\end{align}

We now turn to the mean square of $\mathcal P_k$, the analysis of which serves two purposes. First, the subsequent tail estimates and central-range evaluation isolate the contributions to the main asymptotic for $\tilde{V}_{g,k}(X,H)$; second, it supplies the bound
\begin{equation*}
        \frac1{HX}\int_{\mathbb R}g\!\left(\frac{x}{X}\right)
        |\mathcal P_k(x;H)|^2dx
        \ll
        (\log X)^{2k-1},
\end{equation*}
needed to dispose of the
cross-term in \eqref{eq:Cauchy-schwarz}. Once this has been proven, the replacement of
$\Delta_k$ by $\mathcal P_k$ in the variance, at leading order, follows from
\eqref{eq:remainder-small}. Thus, the reduction we aim to prove is
\begin{align}
\nonumber \tilde V_{g,k}(X,H)
= &\frac{1}{HX}\int_{\mathbb{R}} g \left(\frac{x}{X}\right)\, \biggl|\frac{1}{2\pi }\int_{-T}^{T} \biggl(-\Ldz(\sigma_{*} + it)\biggr)^{k}\,{\frac{(x+H)^{\sigma_{*} + it}-x^{\sigma_{*} + it}}{\sigma_{*} + it}}\,dt\biggr|^2\,dx \\
&\qquad \;+\; o((\log{X})^{2k-1}). 
\label{eq:V-contour-reduction}
\end{align} 
To isolate the main contribution to the mean square in
\eqref{eq:V-contour-reduction}, set
\begin{equation*}
Z:=(\log\log X)^A,\qquad A>2.
\end{equation*}
We shall show that, after normalisation by $HX$, only the portion of the
$t$-integral within $\frac{X}{HZ}\leq |t|\leq\frac{XZ}{H}$ contributes to the leading term. Divide the excluded range into
\begin{equation*}
        E_1:=\left[0,\frac{X}{HZ}\right],
        \qquad
        E_2:=\left[\frac{XZ}{H},\frac{XZ}{H}\log T\right],
        \qquad
        E_3:=\left[\frac{XZ}{H}\log T,T\right].
\end{equation*}
The choice of $\epsilon$ gives
\begin{equation*}
\frac{(Z/h)\log T}{T}
=
\frac{Z\log T}{X^{\epsilon}\sqrt H}
=o(1),
\end{equation*}
so these ranges are contained in $[0,T]$ for all sufficiently large
$X$. We will show that these regions offer a negligible contribution, i.e.
\begin{align}
\label{eq:tail-negligible}
&\frac{1}{HX}\int_{\mathbb R} g\!\left(\frac{x}{X}\right)
 \biggl|
 \frac{1}{2\pi}
 \int_{|t|\in E_1 \cup E_2 \cup E_3}
 \biggl(-\Ldz(\sigma_*+it)\biggr)^k
 \frac{(x+H)^{\sigma_*+it}-x^{\sigma_*+it}}
      {\sigma_*+it}\,dt
 \biggr|^2dx  \nonumber                                   
 \\
&\qquad\qquad  =
o\!\left((\log X)^{2k-1}\right). 
\end{align}
Using a Saffari--Vaughan-type reduction, i.e. proceeding as in
\cite[Page 25]{kn:safvau} (see also \cite[Page 5]{kn:matrad-liouville}, \cite[Page 12]{bettinconrey}), it is enough for
\eqref{eq:tail-negligible} to prove
\begin{align}
\label{eq:saffari-vaughan}
&\frac{X^{2/\log T}}{H}\int_{\mathbb R} g\!\left(\frac{x}{X}\right)
 \biggl|
 \frac{1}{2\pi}
 \int_{|t|\in E_1 \cup E_2 \cup E_3}
 \biggl(-\Ldz(\sigma_*+it)\biggr)^k
 x^{it}
 \frac{(1+u)^{\sigma_*+it}-1}
      {\sigma_* +it}\,dt
 \biggr|^2dx                                 \nonumber     \\
&\qquad =
o\!\left((\log X)^{2k-1}\right)
\end{align}
uniformly for all $u \ll H/X$. Note that this reduction is beneficial for giving upper bounds on the $|t| \in [0, \frac{X}{HZ}] \cup [\frac{XZ}{H}, T]$ regions; for the exact asymptotic in the central range $|t| \in [\frac{X}{HZ}, \frac{XZ}{H}]$ obtained in Section ~\ref{subsec:central-range}, we will continue with the exact kernel. Expanding the squared modulus 
in \eqref{eq:saffari-vaughan} gives a double integral in variables
$t_1,t_2$. After setting $s=t_1-t_2$ and $z=t_1$, we have
\begin{align}
\label{eq:LambdaMellin-tail}
&\frac{X^{1+2/\log T}}{H(2\pi)^2}
\iint_{\substack{|z|,\,|z-s|\in E_1\cup E_2\cup E_3}}
X^{is}\hat g(1+is) \\
\quad \times
\Bigl(&-\Ldz(\sigma_*+iz)\Bigr)^k
\Bigl(-\Ldz(\sigma_*+i(s-z))\Bigr)^k
\frac{\bigl((1+u)^{\sigma_*+iz}-1\bigr)
      \bigl((1+u)^{\sigma_*+i(s-z)}-1\bigr)}
     {(\sigma_*+iz)(\sigma_*+i(s-z))}\,dz\,ds . \nonumber
\end{align}
The $x$-integral is evaluated through the Mellin transform of $g$:
\begin{equation}
\label{eq:mellin-decay}
        \hat g(w):=\int_0^\infty g(v)v^{w-1}\,dv.
\end{equation}
After the change of variables $x=Xv$, we have
\begin{equation*}
        \int_{\mathbb R}g\!\left(\frac{x}{X}\right)x^{is}\,dx
        =
        X^{1+is}
        \hat g\!\left(1+is\right),
\end{equation*}
the difference variable $s=t_1-t_2$ appears with the rapidly decaying weight $\hat g(1+is)$. We now obtain simple bounds for the short-interval numerator in \eqref{eq:LambdaMellin-tail}. Define $h:=\frac HX$ and
\begin{equation}
\label{eq:F-kernel}
        F(z,s,u):=
        \bigl((1+u)^{\sigma_*+iz}-1\bigr)
        \bigl((1+u)^{\sigma_*-iz+is}-1\bigr).
\end{equation}
The identity $ (1+u)^w-1
        =
        w\int_0^u(1+r)^{w-1}\,dr$
implies that
\begin{equation*}
\left|\frac{(1+u)^w-1}{w}\right|
\leq
\int_0^u(1+r)^{\sigma_*-1}\,dr
\ll u.
\end{equation*}
On the other hand, we have
\begin{equation*}
\left|\frac{(1+u)^w-1}{w}\right|
\leq
\frac{(1+u)^{\sigma_*}+1}{|w|}
\ll |w|^{-1},
\end{equation*}
because $u \ll h =o(1)$ and $\sigma_*$ remains bounded. Combining these two bounds yields
\begin{equation*}
        \left|\frac{(1+u)^w-1}{w}\right|
        \ll
        \min(u,|w|^{-1}),   \qquad 0\le u\ll h,
\end{equation*}
uniformly over the $w$-interval. Consequently, if $a=\sigma_*+iz$ and $b=\sigma_*-iz+is$, then
\begin{align}
\label{eq:AB-split-tail}
&\left|\Ldz(a)\right|^k
 \left|\Ldz(b)\right|^k
 \left|\frac{F(z,s,u)}{ab}\right|                                      \notag\\
&\qquad\ll
\left|\Ldz(a)\right|^k\left|\Ldz(b)\right|^k
\min(u,|a|^{-1})\min(u,|b|^{-1})                                      \notag\\
&\qquad\ll
\left|\Ldz(a)\right|^{2k}
\min\!\left(u^2,|a|^{-2}\right)
+
\left|\Ldz(b)\right|^{2k}
\min\!\left(u^2,|b|^{-2}\right),
\end{align}
where we used $AB\le(A^2+B^2)/2$ in the last line. Before estimating separately the contributions from $E_1$, $E_2$, and $E_3$ to the double integral in \eqref{eq:LambdaMellin-tail}, we first restrict this integral to cases in which $t_1$ and $t_2$ are close; this is because the rapid decay in \eqref{eq:mellin-decay} makes pairs with large $|t_1-t_2|$ negligible. Fix $\eta>0$ and set \begin{equation*} S:=(\log X)^\eta. \end{equation*} We now estimate the portion with $|s|>S$; once this has been shown to be negligible, the remaining analysis may be restricted to $|s|\leq S$.

\subsection{Bounds on tail regions} 
Taking absolute values in \eqref{eq:LambdaMellin-tail} and applying
\eqref{eq:AB-split-tail}, the contribution from $|s|>S$ is bounded
by the two terms on the right-hand side of
\eqref{eq:AB-split-tail}. For the second term, we make the substitution
$y=s-z$ for each fixed $s$. This interchanges the conditions
\begin{equation*}
|z|\in E_1\cup E_2\cup E_3,
\qquad
|z-s|\in E_1\cup E_2\cup E_3,
\end{equation*}
from \eqref{eq:LambdaMellin-tail} and transforms the second term into
the first. For each fixed $s$, the substitution preserves the measure $dz$ and
interchanges the two terms. Hence their $z$-integrals are equal, and
their combined contribution to the $|s|>S$ part of
\eqref{eq:LambdaMellin-tail} is twice either contribution. Since
the resulting integrand is non-negative, we may omit the condition
$|z-s|\in E_1\cup E_2\cup E_3$ and enlarge the remaining $z$-range
to $|z|\leq T$. After absorbing the resulting factor of $2$ and using
\begin{equation*}
\frac{X^{1+2/\log T}}{H}\asymp\frac1h,
\end{equation*}
it is sufficient to prove
\begin{equation}
\label{eq:s-tail-one-variable-reduction}
\frac1h
\int_{|s|>S}|\hat g(1+is)|
\int_{|z|\leq T}
\left|\Ldz(\sigma_*+iz)\right|^{2k}
\min\!\left(u^2,|\sigma_*+iz|^{-2}\right)
\,dz\,ds = o\!\left((\log X)^{2k-1}\right).
\end{equation}
As $g\in C_c^\infty((0,\infty))$, repeated integration by parts
in the definition \eqref{eq:mellin-decay} gives
\begin{equation}
\label{eq:mellin-rapid-decay}
|\hat g(1+is)|\ll_{B,g}(1+|s|)^{-B}
\end{equation}
for every $B>0$. In particular, for $B>1$,
\begin{equation}
\label{eq:mellin-L1-bound}
\int_{\mathbb R}|\hat g(1+is)|\,ds
\ll_g 1,
\end{equation}
and
\begin{equation}
\label{eq:mellin-tail-bound}
\int_{|s|>S}|\hat g(1+is)|\,ds
\ll_{B,g}S^{1-B}.
\end{equation}
We first dispose of the bounded range $|z|<1$. Since $\sigma_*=\frac12+\frac1{\log T}$ for all sufficiently large $T$, we have
$\frac12<\sigma_*\leq\frac34$. Therefore, for real $z$ with
$|z|<1$, the point $\sigma_*+iz$ lies in the fixed rectangle
\begin{equation*}
\left\{w\in\mathbb C:
\frac12\leq\Re w\leq\frac34,\quad |\Im w|\leq1
\right\}.
\end{equation*}
As $\zeta(w)$ has no zeros or poles in this rectangle,
$\zeta'/\zeta$ is continuous and bounded there. Hence
\begin{equation}
\label{eq:bounded-height-log-derivative}
\Ldz(\sigma_*+iz)\ll1,
\qquad |z|<1,
\end{equation}
with an implied constant independent of $T$. Since
$\min\!\left(u^2,|\sigma_*+iz|^{-2}\right)\leq u^2$,
\eqref{eq:mellin-tail-bound} gives
\begin{align}
\label{eq:s-tail-bounded-z}
&\frac1h
\int_{|s|>S}|\hat g(1+is)|
\int_{|z|<1}
\left|\Ldz(\sigma_*+iz)\right|^{2k}
\min\!\left(u^2,|\sigma_*+iz|^{-2}\right)
\,dz\,ds \notag\\
&\qquad\ll_{B,g}
\frac{u^2}{h}S^{1-B}
\ll_{B,g}hS^{1-B},
\end{align}
where $u\ll h$ gives $u^2/h\ll h$. Since
$h=X^{-1/\alpha}$ and $S=(\log X)^\eta$, the last expression is
$o\!\left((\log X)^{2k-1}\right)$.
For the range $1\leq|z|\leq T$, the inner integrand in
\eqref{eq:s-tail-one-variable-reduction} is even in the real
variable $z$, since
\begin{equation*}
\left|\Ldz(\sigma_*-iz)\right|
=
\left|\Ldz(\sigma_*+iz)\right|,
\qquad
|\sigma_*-iz|=|\sigma_*+iz|.
\end{equation*}
Thus the positive and negative $z$-ranges contribute equally.
Moreover, for $z\geq1$, $|\sigma_*+iz|^{-2}\leq z^{-2}$ and, after absorbing the factor of $2$, it therefore suffices to bound
\begin{equation}
\label{eq:s-tail-one-variable-moment}
\int_1^T
\left|\Ldz(\sigma_*+iz)\right|^{2k}
\min\!\left(u^2,z^{-2}\right)\,dz.
\end{equation}
Define the cumulative moment
\begin{equation}
\label{eq:definition-of-A}
A(Y):=
\int_1^Y
\left|\Ldz(\sigma_*+iz)\right|^{2k}\,dz.
\end{equation}
We claim that, uniformly for $1\leq Y\leq T$,
\begin{equation}
\label{eq:A-bound}
A(Y)\ll Y(\log X)^{2k}.
\end{equation}
The bound \eqref{eq:first-corollary-statement} from
Proposition~\ref{eq:final-corollary} holds once its height parameter is
sufficiently large. We therefore fix a constant $Y_0\geq2$, independent of $X$,
$T$, and $Y$, sufficiently large that
\eqref{eq:first-corollary-statement} is valid whenever its averaging
height is at least $Y_0$. For $Y_0\leq Y\leq T$, we apply \eqref{eq:first-corollary-statement} with averaging height $Y$
and shifts
\begin{equation*}
z_1=\sigma_*-\frac12=\frac1{\log T},
\qquad
z_2=\frac12-\sigma_*=-\frac1{\log T}.
\end{equation*}
Here, applying the proposition with averaging height $Y$ means that the
height parameter appearing in the proposition is replaced by $Y$; the
parameter $T$ in $\sigma_*$ and in the shifts above continues to
denote the original Perron height. The shifts are real, and the
inequality $Y\leq T$ gives
\begin{equation*}
0<-\Re z_2=\frac1{\log T}\leq\frac3{\log Y},
\qquad
0<\Re(z_1-z_2)=\frac2{\log T}\leq\frac3{\log Y}.
\end{equation*}
The non-negativity of the integrand in
\eqref{eq:definition-of-A} now gives
\begin{align*}
A(Y)
&\leq
\int_0^Y
\left|\Ldz(\sigma_*+iz)\right|^{2k}\,dz  \ll
Y\left(\log Y+\frac12\log T\right)^{2k} \ll
Y(\log X)^{2k},
\end{align*}
where we use $Y\leq T$ and
$\log T\asymp\log X$. We next cover the bounded range $1\leq Y\leq Y_0$, where the
asymptotic proposition is not required. The standard local
partial-fraction expansion for $\zeta'/\zeta$, together with RH,
gives, uniformly for $1\leq z\leq Y_0$,
\begin{align*}
\left|\Ldz(\sigma_*+iz)\right|
\ll_{Y_0}
1+
\sum_{\substack{\rho=\frac12+i\gamma\\|\gamma-z|\leq1}}
\frac{1}{|\sigma_*+iz-\rho|}\ll_{Y_0}\log T.
\end{align*}
Only finitely many zeros occur in this fixed bounded range, while RH
gives $|\sigma_*+iz-\rho|
\geq
\sigma_*-\frac12
=
\frac1{\log T}.$ It follows that, for $1\leq Y\leq Y_0$,
\begin{equation*}
A(Y)
\ll_{k,Y_0}
Y(\log T)^{2k}
\ll_{k,Y_0}
Y(\log X)^{2k}.
\end{equation*}
The constant $Y_0$ is fixed, so its dependence may be absorbed into
the implied constant. This proves \eqref{eq:A-bound} throughout
$1\leq Y\leq T$.

 By Stieltjes integration, we write \eqref{eq:s-tail-one-variable-moment} as
\begin{equation}
\label{eq:stieljes-estimate-unshifted}
\int_1^T
\left|\Ldz(\sigma_*+iz)\right|^{2k}
\min(u^2,z^{-2})\,dz
=
\int_1^T \min(u^2,z^{-2})\,dA(z).
\end{equation}
The case $u=0$ is immediate, as the integral vanishes. So we
assume $u>0$; the bound $u\ll h=X^{-1/\alpha}$ gives $u<1$ for all
sufficiently large $X$, so that $1/u>1$. We split at the point $z=1/u$, at which the two parts of the minimum are equal. If $1/u\le T$, then we have
\begin{equation*}
\int_1^T \min(u^2,z^{-2})\,dA(z)
=
u^2A(1/u)+\int_{1/u}^T z^{-2}\,dA(z).
\end{equation*}
By the bound in \eqref{eq:A-bound}, $u^2A(1/u)\ll u(\log X)^{2k}$. As for the second term, partial summation and the bound for $A(z)$ give
\begin{equation*}
\int_{1/u}^T z^{-2}\,dA(z)
\ll
T^{-2}A(T)+u^2A(1/u)+\int_{1/u}^T A(z)z^{-3}\,dz
\ll
u(\log X)^{2k},
\end{equation*}
using $T^{-1}\le u$ and
$\int_{1/u}^T z^{-2}\,dz \ll u$. If instead $1/u>T$, then $\min(u^2,z^{-2})=u^2$ throughout
$1\le z\le T$, and therefore
\begin{equation*}
\int_1^T \min(u^2,z^{-2})\,dA(z)
=
u^2A(T)
\ll
u^2T(\log X)^{2k}
\ll
u(\log X)^{2k}.
\end{equation*}
As $u\ll h$, both cases give
\begin{equation}
\label{eq:s-tail-moment-bound}
\int_1^T
\left|\Ldz(\sigma_*+iz)\right|^{2k}
\min\!\left(u^2,z^{-2}\right)\,dz
\ll
h(\log X)^{2k}.
\end{equation}
Applying \eqref{eq:mellin-tail-bound} and
\eqref{eq:s-tail-moment-bound} to the part of
\eqref{eq:s-tail-one-variable-reduction} with
$1\leq |z|\leq T$ gives
\begin{align*}
&\frac1h
\int_{|s|>S}|\hat g(1+is)|
\int_{1\leq|z|\leq T}
\left|\Ldz(\sigma_*+iz)\right|^{2k}
\min\!\left(u^2,|\sigma_*+iz|^{-2}\right)
\,dz\,ds \\
&\qquad\ll_{B,g}
(\log X)^{2k}S^{1-B}.
\end{align*}
Together with the bounded-$z$ estimate
\eqref{eq:s-tail-bounded-z}, this shows that the full expression in
\eqref{eq:s-tail-one-variable-reduction} is
\begin{equation*}
\ll_{B,g}
hS^{1-B}
+
(\log X)^{2k}S^{1-B}
\ll_{B,g}
(\log X)^{2k}S^{1-B}.
\end{equation*}
Selecting $B$ so that $\eta(B-1)>1$, and recalling that
$S=(\log X)^\eta$, we obtain
\begin{equation*}
(\log X)^{2k}S^{1-B}
=
(\log X)^{2k-\eta(B-1)}
=
o\!\left((\log X)^{2k-1}\right).
\end{equation*}
This proves \eqref{eq:s-tail-one-variable-reduction}. Hence, when
estimating the contributions from $E_1$, $E_2$, and $E_3$, we may
restrict the difference variable to $|s|\leq S=(\log X)^\eta$.
For each fixed $s$ and each $z$-range $E_j$, the additional condition
\begin{equation*}
|z-s|\in E_1\cup E_2\cup E_3
\end{equation*}
restricts the $z$-domain to a finite union of sub-intervals. We suppress this condition from the notation moving forward. When the
shifted product is retained, this suppression is purely notational:
partial summation is applied separately on each actual sub-interval,
and the resulting estimates are uniform in its endpoints. We enlarge
the domain only after \eqref{eq:AB-split-tail} has produced a
non-negative majorant. 

\paragraph{Bounding $E_1$.} 
We now estimate the contribution to the double integral in
\eqref{eq:LambdaMellin-tail} from the first range; we must show that the contribution is $o\!\left((\log X)^{2k-1}\right)$ uniformly for $u\ll h$ in order to prove
\eqref{eq:tail-negligible}. Recall that
\begin{equation*}
E_1=\left[0,\frac{X}{HZ}\right]
=\left[0,\frac1{hZ}\right],
\qquad
Y_1:=\frac1{hZ}.
\end{equation*}
After taking absolute values and applying
\eqref{eq:AB-split-tail}, the $E_1$-contribution consists of the
first term restricted by $|z|\leq Y_1$ and the second term restricted
by $|z-s|\leq Y_1$. These two contributions are interchanged by the
substitution $z\mapsto s-z$, so it is enough to treat the first. We
begin with its sub-range $|z|\leq h^{-1/2}$.
Omitting the remaining domain condition from
\eqref{eq:LambdaMellin-tail}, which can only enlarge this
non-negative upper bound, it is enough to estimate
\begin{align*}
\frac1h
\int_{|s|\leq S}|\hat g(1+is)|
\int_{|z|\leq h^{-1/2}}
\left|\Ldz(\sigma_*+iz)\right|^{2k}
\min\!\left(u^2,|\sigma_*+iz|^{-2}\right)
\,dz\,ds .
\end{align*}
On $|z|<1$, the logarithmic derivative is bounded by
\eqref{eq:bounded-height-log-derivative}. On
$1\leq |z|\leq h^{-1/2}$, the definition
\eqref{eq:definition-of-A}, the bound \eqref{eq:A-bound}, and the
evenness of the integrand give
\begin{equation*}
\int_{1\leq |z|\leq h^{-1/2}}
\left|\Ldz(\sigma_*+iz)\right|^{2k}\,dz
\ll A(h^{-1/2})
\ll h^{-1/2}(\log X)^{2k}.
\end{equation*}
Using $\min\!\left(u^2,|\sigma_*+iz|^{-2}\right)\leq u^2$ and the bound \eqref{eq:mellin-L1-bound}, the preceding double integral is
\begin{align*}
&\ll
\frac{u^2}{h}
\left(1+A(h^{-1/2})\right)
\int_{\mathbb R}|\hat g(1+is)|\,ds \\
&\ll
h^{1/2}(\log X)^{2k}
=
o\!\left((\log X)^{2k-1}\right),
\end{align*}
where $u\ll h$ and $h^{1/2}=X^{-1/(2\alpha)}$ were used. Next, we focus on the range $h^{-1/2}\leq |z|\leq Y_1$. Applying
\eqref{eq:AB-split-tail} would give only
\begin{equation*}
O\!\left(\frac{(\log X)^{2k}}{Z}\right),
\end{equation*}
which is not sufficient to obtain
$o((\log X)^{2k-1})$. For this range, we work directly with the shifted product in
\eqref{eq:LambdaMellin-tail}, allowing
\eqref{eq:first-corollary-statement} to provide the required decay
in $s$. For fixed $s$, define
\begin{equation}
\label{eq:shifted-moment-bound-A}
A_s(Y):=
\int_0^Y
\Bigl(-\Ldz(\sigma_*+iz)\Bigr)^k
\Bigl(-\Ldz(\sigma_*-iz+is)\Bigr)^k\,dz.
\end{equation}
The lower limit agrees with the averaging interval in
\eqref{eq:first-corollary-statement}. We apply that estimate at
averaging height $Y$ with $z_1=\sigma_*-\frac12=\frac1{\log T}$ and $z_2=\frac12-\sigma_*-is=-\frac1{\log T}-is$. For $h^{-1/2}\leq Y\leq Y_1$, we have
$\log Y\asymp\log X$, and hence
\begin{equation*}
|s|\leq(\log X)^\eta\ll(\log Y)^\eta.
\end{equation*}
Thus the imaginary-shift conditions in
\eqref{eq:first-corollary-statement} are satisfied; the real-part
conditions follow from $Y\leq Y_1\leq T$. Additionally, we have
\begin{equation}
\label{eq:two-difference-inequalities}
|z_1-z_2|
=
\left|\frac2{\log T}+is\right|
\geq |s|,
\qquad
|z_1-z_2|\geq\frac2{\log T}.
\end{equation}
We will now use the above inequalities to bound the right-hand side of \eqref{eq:first-corollary-statement}. The second inequality in \eqref{eq:two-difference-inequalities} gives
\begin{equation}
\label{bounds-for-rhs}
(\log Y)^{2k}
\left(
1+\frac1{|z_1-z_2|\log Y}
\right)^{2k}
=
\left(
\log Y+\frac1{|z_1-z_2|}
\right)^{2k}
\ll
(\log X)^{2k},
\end{equation}
while the first inequality implies $\min\!\left(|z_1-z_2|,1\right)
\geq
\min\!\left(|s|,1\right)$. Hence the decay factor in \eqref{eq:first-corollary-statement}
satisfies
\begin{equation*}
\left(
1+\min\!\left(|z_1-z_2|,1\right)\log Y
\right)^{-2}
\ll
\left(
1+\min\!\left(|s|,1\right)\log(2Y)
\right)^{-2},
\end{equation*}
where we used $\log(2Y)\asymp\log Y$. Together with \eqref{bounds-for-rhs},
this gives
\begin{equation}
\label{eq:corollary-bounds}
|A_s(Y)|
\ll
Y(\log X)^{2k}D_Y(s)^{-2},
\qquad
D_Y(s):=1+\min\!\left(1,|s|\right)\log(2Y),
\end{equation}
uniformly for $|s|\leq(\log X)^\eta$ and
$h^{-1/2}\leq Y\leq Y_1$. To apply partial summation with respect to $A_s$, we must consider the factor
\begin{equation*}
\frac{F(z,s,u)}
{(\sigma_*+iz)(\sigma_*-iz+is)},
\end{equation*}
where $F$ is defined in \eqref{eq:F-kernel}. We now bound this factor
and its derivative with respect to $z$. Using $\frac{(1+u)^w-1}{w}=\int_0^u(1+r)^{w-1}\,dr$,
we obtain
\begin{equation}
\label{eq:E1-kernel-bound}
\left|
\frac{F(z,s,u)}
{(\sigma_*+iz)(\sigma_*-iz+is)}
\right|
\ll u^2.
\end{equation}
Differentiating under the integral sign gives one extra factor $\log(1+r)\ll r$, hence
\begin{equation}
\label{eq:E1-kernel-derivative}
\left|
\partial_z
\left(
\frac{F(z,s,u)}
{(\sigma_*+iz)(\sigma_*-iz+is)}
\right)
\right|
\ll u^3.
\end{equation}
Thus the short-interval kernel is slowly varying on the $E_1$-range. Applying partial summation with respect to $A_s$,
in the same form as \eqref{eq:stieljes-estimate-unshifted}, gives
\begin{align}
\label{eq:partial-summation-stieljes-IBP}
&\int_{h^{-1/2}}^{Y_1}
\Bigl(-\Ldz(\sigma_*+iz)\Bigr)^k
\Bigl(-\Ldz(\sigma_*-iz+is)\Bigr)^k
\frac{F(z,s,u)}
{(\sigma_*+iz)(\sigma_*-iz+is)}
\,dz \\
&=
\left[
\frac{F(z,s,u)}
{(\sigma_*+iz)(\sigma_*-iz+is)}
A_s(z)
\right]_{h^{-1/2}}^{Y_1} -
\int_{h^{-1/2}}^{Y_1}
A_s(z)
\partial_z
\left(
\frac{F(z,s,u)}
{(\sigma_*+iz)(\sigma_*-iz+is)}
\right)
\,dz . \notag
\end{align}
Using \eqref{eq:corollary-bounds},
\eqref{eq:E1-kernel-bound}, and
\eqref{eq:E1-kernel-derivative}, the absolute value of this
expression is
\begin{align*}
\ll\;&
u^2(\log X)^{2k}
\left(
Y_1D_{Y_1}(s)^{-2}
+
h^{-1/2}D_{h^{-1/2}}(s)^{-2}
\right) +
u^3(\log X)^{2k}
\int_{h^{-1/2}}^{Y_1}YD_Y(s)^{-2}\,dY.
\end{align*}
After multiplying by $|\hat g(1+is)|$, we use
\begin{equation}
\label{eq:mellin-corollary-decay}
\int_{\mathbb R}
|\hat g(1+is)|D_Y(s)^{-2}\,ds
\ll
\frac1{\log(2Y)}.
\end{equation}
This follows from \eqref{eq:mellin-rapid-decay} by splitting the
$s$-integral at $|s|=1/\log(2Y)$ and $|s|=1$. Applying \eqref{eq:mellin-corollary-decay} to the endpoint and
integral terms in
\eqref{eq:partial-summation-stieljes-IBP}, multiplying by
$h^{-1}|\hat g(1+is)|$, and integrating over $|s|\leq S$, we find
that the part of \eqref{eq:LambdaMellin-tail} with
$h^{-1/2}\leq z\leq Y_1$ is bounded by
\begin{align*}
\ll\;&
\frac{1}{h}u^2(\log X)^{2k}
\left(
\frac{Y_1}{\log(2Y_1)}
+
\frac{h^{-1/2}}{\log(2h^{-1/2})}
\right) +
\frac{1}{h}u^3(\log X)^{2k}
\int_{h^{-1/2}}^{Y_1}
\frac{Y}{\log(2Y)}\,dY.
\end{align*}
Including the previously estimated range $|z|\leq h^{-1/2}$, this bounds the entire positive range $0\le z\le Y_1$. Thus the
positive-$z$ part of the $E_1$ contribution to \eqref{eq:LambdaMellin-tail} is
\begin{align*}
\ll\;&
\frac{1}{h}u^2(\log X)^{2k}
\left(
\frac{Y_1}{\log(2Y_1)}
+
\frac{h^{-1/2}}{\log(2h^{-1/2})}
\right) \\
&\qquad \qquad\qquad\quad+
\frac{1}{h}u^3(\log X)^{2k}
\int_{h^{-1/2}}^{Y_1}
\frac{Y}{\log(2Y)}\,dY
+
o\!\left((\log X)^{2k-1}\right).
\end{align*}
Using
\begin{equation*}
\int_{h^{-1/2}}^{Y_1}
\frac{Y}{\log(2Y)}\,dY
\ll
\frac{Y_1^2}{\log(2Y_1)},
\end{equation*}
together with $u\ll h$, $Y_1=1/(hZ)$, and $\log(2Y_1)\asymp\log(2h^{-1/2})\asymp\log X$,
the above bound is
\begin{align*}
\ll\;&
\frac{(\log X)^{2k-1}}{Z}
+
h^{1/2}(\log X)^{2k-1}
+
\frac{(\log X)^{2k-1}}{Z^2}
+
o\!\left((\log X)^{2k-1}\right) =
o\!\left((\log X)^{2k-1}\right),
\end{align*}
because $Z\to\infty$ and $h\to0$.

Therefore, the positive-$z$ part of $E_1$ has the required bound. The negative range $-Y_1\leq z\leq-h^{-1/2}$ is treated identically, by conjugation. As all bounds depend only on
$|s|$ and $g$ is real-valued, the same estimate applies.
\paragraph{Bounding $E_2$.} As with $E_1$, we have already restricted
to $|s|\le(\log X)^\eta$; on the region $E_2:=\left[\frac{Z}{h},\frac{Z}{h}\log T\right]$, we have $ |\sigma_*+iz|\asymp |\sigma_*-iz+is|\asymp |z|$. We now estimate the contribution to the double integral in
\eqref{eq:LambdaMellin-tail} from $E_2$. We split the remaining $s$-range at
\begin{equation*}
        s_0:=\frac{(\log\log X)^2}{\log X}.
\end{equation*}

First suppose $|s|<s_0$. Applying \eqref{eq:AB-split-tail} and
repeating the arguments made in the derivation of
\eqref{eq:s-tail-one-variable-reduction}, with the $E_2$-restriction
imposed on $|z|$ in the first term and on $|z-s|$ in the second, we
can bound the contribution from this range by
\begin{align*}
\ll
\frac1h
\int_{|s|<s_0}
|\hat g(1+is)|
\int_{\frac{Z}{h}}^{\frac{Z}{h}\log T}
\left|\Ldz(\sigma_*+iz)\right|^{2k}z^{-2}\,dz\,ds.
\end{align*}
By the same Stieltjes argument as in the calculation following
\eqref{eq:stieljes-estimate-unshifted}, together with
\eqref{eq:A-bound}, we have
\begin{align}
\int_{\frac{Z}{h}}^{\frac{Z}{h}\log{T}}
\left|\Ldz(\sigma_*+iz)\right|^{2k}z^{-2}\,dz
&=
\int_{\frac{Z}{h}}^{\frac{Z}{h}\log{T}} z^{-2}\,dA(z) \nonumber      \\
&\ll
(\log X)^{2k}
\left(
\frac{h}{Z}
+
\int_{\frac{Z}{h}}^{\frac{Z}{h}\log{T}}z^{-2}\,dz
\right)                                     \nonumber  \\
&\ll
\frac{h}{Z}(\log X)^{2k}.
\label{eq:stieljes-partialsummation-approach}
\end{align}
Therefore
\begin{equation*}
\ll
\frac1h
\frac{(\log\log X)^2}{\log X}
\frac{h}{Z}(\log X)^{2k}
=
\frac{(\log X)^{2k-1}(\log\log X)^2}{Z}
=
o((\log X)^{2k-1}),
\end{equation*}
since $Z/(\log\log X)^2\to\infty$. For $s_0\le |s|\le(\log X)^\eta$,
the decoupling from \eqref{eq:AB-split-tail} would lose the required saving in $s$. Instead, we apply the form of Proposition~\ref{eq:final-corollary} stated in \eqref{eq:second-corollary-statement} for this range and then use partial summation, as in the treatment of $E_1$. On the range $E_2$, we have $F(z,s,u)=O(1)$, and the kernel factor satisfies
\begin{equation*}
        \left|
        \frac{F(z,s,u)}
             {(\sigma_*+iz)(\sigma_*-iz+is)}
        \right|
        \ll |z|^{-2}.
\end{equation*}
Furthermore,
\begin{equation*}
\partial_zF(z,s,u)
=
i\log(1+u)
\left((1+u)^{\sigma_*-iz+is}
      -(1+u)^{\sigma_*+iz}\right),
\end{equation*}
and because $u\ll h$, we have $\partial_zF(z,s,u)\ll h$. It follows that
\begin{equation*}
        \partial_z
        \left(
        \frac{F(z,s,u)}
             {(\sigma_*+iz)(\sigma_*-iz+is)}
        \right)
        \ll
        h|z|^{-2}+|z|^{-3}.
\end{equation*}
For fixed $s$, we use the cumulative shifted moment $A_s(Y)$ as defined in \eqref{eq:shifted-moment-bound-A}. Then we apply the proposition and, uniformly on $\frac{Z}{h}\le Y\le \frac{Z}{h}\log T$, we have
\begin{equation}
\label{eq:final-As-bound}
        A_s(Y)
        \ll
        Y(\log X)^{2k}
        \left(1+\min(1,|s|)\log X\right)^{-2},
\end{equation}
as $\log Y\asymp\log X$ throughout $E_2$. Unlike the full estimate \eqref{eq:first-corollary-statement}, the
simplified estimate \eqref{eq:second-corollary-statement} omits the
factor
\begin{equation*}
\left(
1+\frac{1}{|z_1-z_2|\log Y}
\right)^{2k},
\end{equation*}
and therefore requires $|z_1-z_2|\log Y$ to be bounded away from
zero. On the present range,
\begin{equation*}
|z_1-z_2|\log Y
\geq |s|\log Y
\geq s_0\log Y
\asymp(\log\log X)^2.
\end{equation*}
Thus this additional condition is satisfied, and
\eqref{eq:second-corollary-statement} gives \eqref{eq:final-As-bound}
uniformly over the stated ranges of $s$ and $Y$. Partial summation then gives
\begin{align*}
&\int_{\frac{Z}{h}}^{\frac{Z}{h}\log{T}}
\Bigl(-\Ldz(\sigma_*+iz)\Bigr)^k
\Bigl(-\Ldz(\sigma_*-iz+is)\Bigr)^k
\frac{F(z,s,u)}
     {(\sigma_*+iz)(\sigma_*-iz+is)}
\,dz                                      \\
& =
\left[
\frac{F(z,s,u)}
     {(\sigma_*+iz)(\sigma_*-iz+is)}
A_s(z)
\right]_{\frac{Z}{h}}^{\frac{Z}{h}\log{T}}                                      -
\int_{\frac{Z}{h}}^{\frac{Z}{h}\log{T}}
A_s(z)\,
\partial_z
\left(
\frac{F(z,s,u)}
     {(\sigma_*+iz)(\sigma_*-iz+is)}
\right)
\,dz .
\end{align*}
Using the bounds for $F(z,s,u)$ and its derivative, the endpoint terms provide
\begin{equation*}
\ll
(\log X)^{2k}
\left(1+\min(1,|s|)\log X\right)^{-2}
\frac{h}{Z},
\end{equation*}
while the integral term contributes
\begin{align*}
&\ll
(\log X)^{2k}
\left(1+\min(1,|s|)\log X\right)^{-2}
\int_{\frac{Z}{h}}^{\frac{Z}{h}\log{T}}
z\left(hz^{-2}+z^{-3}\right)\,dz       \\
&\ll
(\log X)^{2k}
\left(1+\min(1,|s|)\log X\right)^{-2}
\left(
        h\log\log X+\frac{H}{XZ}
\right).
\end{align*}
We then have
\begin{align*}
&\int_{\frac{Z}{h}}^{\frac{Z}{h}\log{T}}
\Bigl(-\Ldz(\sigma_*+iz)\Bigr)^k
\Bigl(-\Ldz(\sigma_*-iz+is)\Bigr)^k
\frac{F(z,s,u)}
     {(\sigma_*+iz)(\sigma_*-iz+is)}
\,dz                                      \\
&\qquad \qquad\qquad\qquad\ll
(\log X)^{2k}
\left(1+\min(1,|s|)\log X\right)^{-2}
\left(
        \frac{h}{Z}+h\log\log X
\right). \nonumber
\end{align*}
Multiplying by $|\hat g(1+is)|$, integrating over $\frac{(\log\log X)^2}{\log X} \le |s|\le(\log X)^\eta$ and including the outside factor $1/h$, we obtain
\begin{align}
&\ll
\frac1h
(\log X)^{2k}
\left(
        \frac{h}{Z}+h\log\log X
\right)   \nonumber     \\
&\qquad\qquad \times
\int_{|s|\ge(\log\log X)^2/\log X}
|\hat g(1+is)|
\left(1+\min(1,|s|)\log X\right)^{-2}\,ds .
\label{eq:bound-on-average-g}
\end{align}
To bound the above, we split the $s$-integral into
$s_0\le |s|\le1$ and $|s|>1$.  As $\hat g$ is bounded and rapidly decaying,
\begin{equation*}
\int_{|s|\ge s_0}|\hat g(1+is)|
(1+\min(1,|s|)\log X)^{-2}\,ds
\ll
\int_{s_0}^{1}(1+s\log X)^{-2}\,ds+(\log X)^{-2}.
\end{equation*}
Changing variables $v=s\log X$, we have
\begin{equation}
\int_{s_0}^{1}(1+s\log X)^{-2}\,ds
=
\frac1{\log X}
\int_{(\log\log X)^2}^{\log X}(1+v)^{-2}\,dv
\ll
\frac1{\log X(\log\log X)^2}.
\label{eq:E2-changeofvariables}
\end{equation}
From this, we can bound \eqref{eq:bound-on-average-g} and arrive at
\begin{align*}
&\ll
(\log X)^{2k}
\left(
        \frac1Z+\log\log X
\right)
\frac{1}{\log X(\log\log X)^2}        \\
&\ll
\frac{(\log X)^{2k-1}}{\log\log X}
+
\frac{(\log X)^{2k-1}}{Z(\log\log X)^2} =  o \Bigl((\log X)^{2k-1}\Bigr).
\end{align*}
The two ranges of $s$ show that the positive $z$-part of $E_2$ is
negligible. The same estimate applies to the negative $z$-interval
by conjugation, so the full $E_2$ contribution is
$o((\log X)^{2k-1})$.
\paragraph{Bounding $E_3$.} Now we consider the range $E_3 = \left[\frac{Z}{h}\log T,T\right]$ with $|s|\le(\log X)^\eta$. On this range $|z|\ge \frac{Z}{h}\log T$ and, as in $E_2$, we get $|\sigma_*+iz|\asymp |\sigma_*+i(s-z)|\asymp |z|$
due to the truncation in $s$. On $E_3$, we can use
\eqref{eq:AB-split-tail} and proceed as we did for the near-diagonal part of
$E_2$, i.e. reducing to a one-variable integral by the
change of variables $y=s-z$, followed by enlarging the domain. Thus the
$E_3$-contribution to \eqref{eq:saffari-vaughan} is bounded by
\begin{equation*}
\ll
\frac1h
\int_{|s|\le(\log X)^\eta}
|\hat g(1+is)|
\int_{\frac{Z}{h}\log T}^{T}
\left|\Ldz(\sigma_*+iz)\right|^{2k}z^{-2}\,dz\,ds .
\end{equation*}
By \eqref{eq:definition-of-A}, followed by Stieltjes integration by
parts using \eqref{eq:A-bound}, we have
\begin{equation*}
\int_{\frac{Z}{h}\log T}^{T}
\left|\Ldz(\sigma_*+iz)\right|^{2k}z^{-2}\,dz
=
\int_{\frac{Z}{h}\log T}^{T}z^{-2}\,dA(z)
\ll
\frac{h(\log X)^{2k}}{Z\log T}.
\end{equation*}
Moreover, \eqref{eq:mellin-L1-bound} gives
\begin{equation*}
\int_{|s|\leq(\log X)^\eta}
|\hat g(1+is)|\,ds
\ll_g 1.
\end{equation*}
This means that the region $E_3$ contributes
\begin{equation*}
\ll
\frac1h
\frac{h(\log X)^{2k}}{Z\log T}
\int_{|s|\le(\log X)^\eta}
|\hat g(1+is)|\,ds
\ll
\frac{(\log X)^{2k}}{Z\log T}.
\end{equation*}
As $\log T\asymp\log X$, we obtain
\begin{equation*}
        \ll \frac{(\log X)^{2k-1}}{Z}
        =
        o((\log X)^{2k-1}).
\end{equation*}
The negative $z$-range is treated identically, by conjugation. The choice $Z=(\log\log X)^A$, with $A>2$, is now sufficient for all three
tail ranges: the $E_1$ and $E_3$ bounds gain a factor $1/Z$, while the near part of $E_2$
gains $Z^{-1}(\log\log X)^2=o(1)$ and the remaining part of $E_2$ gains
$1/\log\log X$. Hence, we have proven \eqref{eq:saffari-vaughan} and the total
contribution of $z \in [0, \frac{1}{hZ}] \cup [\frac{Z}{h}, T]$ to the smoothed variance \eqref{eq:V-contour-reduction} is $o((\log X)^{2k-1})$.
\subsection{The central range} \label{subsec:central-range} We now focus on the central contribution given by
\begin{equation}
\label{eq:central-region-analysis}
 \frac1{HX}\int_{\mathbb R}g\!\left(\frac{x}{X}\right)
        \bigg|\frac1{2\pi}
        \int_{\frac1{hZ}\le |t|\le \frac Zh}
        \left(-\Ldz(\sigma_*+it)\right)^k
        \frac{(x+H)^{\sigma_*+it}-x^{\sigma_*+it}}
             {\sigma_*+it}\,dt\bigg|^2dx .
\end{equation}
Expanding the squared modulus and putting $s=t_1-t_2$, $z=t_1$, we obtain
\begin{equation}
\label{eq:expanded-central-region}
\frac{1}{H(2\pi)^2}
\iint_{\frac1{hZ}\le |z|, |z-s| \le \frac Zh}
\Bigl(-\Ldz(\sigma_*+iz)\Bigr)^k
\Bigl(-\Ldz(\sigma_*-iz+is)\Bigr)^k
\mathcal L(z,s)\,dz\,ds ,
\end{equation}
where
\begin{align}
\label{eq:L-kernel-def}
\mathcal L(z,s)
&:=
\frac{1}{X}\int_{\mathbb R} g\!\left(\frac{x}{X}\right)
\frac{(x+H)^{\sigma_*+iz}-x^{\sigma_*+iz}}{\sigma_*+iz}\,
\frac{(x+H)^{\sigma_*-iz+is}-x^{\sigma_*-iz+is}}{\sigma_*-iz+is}\,dx       \\
&=
\frac{1}{X}\int_{\mathbb R} g\!\left(\frac{x}{X}\right)
\frac{F(z,s,H/x)}
     {(\sigma_*+iz)(\sigma_*-iz+is)}
x^{2\sigma_*+is}\,dx .
\nonumber
\end{align}
Here $F$ is the function defined in \eqref{eq:F-kernel}. To bound $\mathcal L$, note first that on the support of $g(x/X)$ we have
$H/x\asymp h$. Writing $u=H/x$, we have $x\partial_x=-u\partial_u$ and each application of $u\partial_u$ to $(1+u)^a$ or $(1+u)^b$ offers at most $O(u|a|)$ or $O(u|b|)$, respectively. Hence, for every fixed integer $m\ge0$,
\begin{equation}
\label{eq:preceding-bound-1}
        \bigl|(x\partial_x)^mF(z,s,H/x)\bigr|
        \ll_m
        \bigl(1+h(1+|z|+|z-s|)\bigr)^m.
\end{equation}
Moreover, differentiating $F$ with respect to $z$ introduces an additional
factor $\log(1+u)\ll u\ll h$. Therefore
\begin{equation}
\label{eq:preceding-bound-2}
        \bigl|(x\partial_x)^m\partial_zF(z,s,H/x)\bigr|
        \ll_m
        h\bigl(1+h(1+|z|+|z-s|)\bigr)^m.
\end{equation}
We integrate by parts $\ell \geq 1$ times in \eqref{eq:L-kernel-def} using $x\partial_x \, x^{2\sigma_* + is}=(2\sigma_* +is)x^{2\sigma_* +is}$; each integration by parts gains a factor $(1+|s|)^{-1}$, while the adjoint of
$x\partial_x$ only produces a fixed linear combination of
$(x\partial_x)^j$, $j\le \ell$.  The boundary terms vanish because $g(x/X)$ is smooth and compactly
supported. The bound \eqref{eq:preceding-bound-1} and the smooth compact support of
$g(x/X)$ give
\begin{equation*}
        \mathcal L(z,s)
        \ll_\ell
        \frac{X\bigl(1+h(1+|z|+|z-s|)\bigr)^\ell}
        {(1+|s|)^\ell
        |(\sigma_*+iz)(\sigma_*-iz+is)|}.
\end{equation*}
Here the factor $X$ comes from the $x$-integral, since
$X^{-1}\int g(x/X)x^{2\sigma_*}\,dx\ll X^{2\sigma_*}\asymp X$. Similarly, differentiating the factors depending explicitly on $z$, we obtain
\begin{equation*}
        \partial_z\mathcal L(z,s)
        \ll_\ell
        \frac{X\bigl(1+h(1+|z|+|z-s|)\bigr)^\ell}
        {(1+|s|)^\ell
        |(\sigma_*+iz)(\sigma_*-iz+is)|}
        \left(
        h+\frac1{1+|z|}+\frac1{1+|z-s|}
        \right).
\end{equation*}
These estimates hold for every integer $\ell\geq0$, with the case
$\ell=0$ following directly from \eqref{eq:L-kernel-def}. In the central range, we may assume $1+h(1+|z|+|z-s|)\ll Z$ and $ h+\frac1{1+|z|}+\frac1{1+|z-s|}\ll hZ$. Therefore
\begin{equation}
\label{eq:L-kernel-bounds}
        \mathcal L(z,s)
        \ll_\ell
        \frac{Z^\ell X}
        {(1+|s|)^\ell
        |(\sigma_*+iz)(\sigma_*-iz+is)|},
\end{equation}
and
\begin{equation}
\label{eq:L-kernel-derivative-bounds}
        \partial_z\mathcal L(z,s)
        \ll_\ell
        \frac{Z^{\ell+1}H}
        {(1+|s|)^\ell
        |(\sigma_*+iz)(\sigma_*-iz+is)|}.
\end{equation}
We first estimate the portion of
\eqref{eq:expanded-central-region} for which
$|s|>(\log X)^\eta$. Using \eqref{eq:L-kernel-bounds} and the same decoupling used in \eqref{eq:AB-split-tail}, we have
\begin{equation*}
\ll_\ell
\frac{Z^\ell X}{H(1+|s|)^\ell} \biggl( \bigg|\Ldz(\sigma_*+iz)\bigg|^{2k}|\sigma_*+iz|^{-2} + \bigg|\Ldz(\sigma_*-iz+is)\bigg|^{2k}|\sigma_*-iz+is|^{-2}  \biggr)
\end{equation*}
The second term is treated by a variable change and enlarging the resulting domain as in the large-$|s|$ tail truncation. As before, it suffices to treat the positive $z$-range. Thus, for fixed $s$ in this range, the corresponding $z$-integral
in \eqref{eq:expanded-central-region} is
\begin{equation*}
\ll_\ell
\frac{Z^\ell X}{H(1+|s|)^\ell}
\int_{1/(hZ)}^{Z/h}
\left|\Ldz(\sigma_*+it)\right|^{2k}t^{-2}\,dt .
\end{equation*}
By a partial summation similar to \eqref{eq:stieljes-partialsummation-approach}, we get
\begin{equation}
\int_{1/(hZ)}^{Z/h}
\left|\Ldz(\sigma_*+it)\right|^{2k}t^{-2}\,dt
=
\int_{1/(hZ)}^{Z/h}t^{-2}\,dA(t)
\ll
hZ(\log X)^{2k},
\end{equation}
which follows from \eqref{eq:A-bound}. Therefore the contribution of \(|s|>(\log X)^\eta\) is
\begin{equation*}
\ll_\ell
Z^{\ell+1}(\log X)^{2k}
\int_{|s|>(\log X)^\eta}(1+|s|)^{-\ell}\,ds .
\end{equation*}
Choosing $\ell$ so large that $\eta(\ell-1)>1$, this is $o\bigl((\log X)^{2k-1}\bigr)$.

Next, we analyse the range
\begin{equation*}
(\log X)^{-1+\varepsilon}<|s|\leq(\log X)^\eta,
\end{equation*}
where $\varepsilon$ is fixed with $0<\varepsilon<1$. On this range,
$|s|\ll|z|$. As before, it is enough to treat the positive
$z$-range; the simultaneous substitution
$(z,s)\mapsto(-z,-s)$ handles the negative range by conjugation.
Taking $\ell=0$ in the bounds above gives
\begin{equation}
\label{eq:kernel-bounds-L}
\mathcal L(z,s)\ll \frac{X}{z^2},
\qquad
\partial_z\mathcal L(z,s)
\ll
X\left(hz^{-2}+z^{-3}\right).
\end{equation}
The argument establishing \eqref{eq:corollary-bounds} also applies
throughout the central range. For
\begin{equation*}
\frac1{hZ}\leq Y\leq\frac Zh,
\end{equation*}
we have $\log Y\asymp\log X$, $Y\leq Z/h=o(T)$. Thus the real-part conditions in \eqref{eq:first-corollary-statement} continue to hold. In addition, as $|s|\leq(\log X)^\eta\ll(\log Y)^\eta$, the imaginary-shift conditions also hold. Hence
\begin{equation}
\label{eq:corollary-bounds-central}
|A_s(Y)|
\ll
Y(\log X)^{2k}
\left(1+\min(1,|s|)\log X\right)^{-2},
\end{equation}
uniformly for $|s|\leq(\log X)^\eta$ and
$\frac1{hZ}\leq Y\leq\frac Zh$.

As in the treatment of $E_2$ for the range $(\log\log X)^2(\log X)^{-1}
\leq |s|\leq(\log X)^\eta$,
the decoupling in \eqref{eq:AB-split-tail} would lose the required
decay in $s$. We retain the shifted product and apply partial
summation with respect to $A_s$, in the form
\eqref{eq:partial-summation-stieljes-IBP}. Following the convention
introduced before bounding $E_1$, we suppress the additional
condition on $|z-s|$ from \eqref{eq:expanded-central-region}; partial
summation is understood to be applied on the actual sub-interval,
uniformly in its endpoints. Using \eqref{eq:kernel-bounds-L} and
\eqref{eq:corollary-bounds-central}, we obtain
\begin{equation*}
\int_{\frac1{hZ}}^{\frac Zh}
\Bigl(-\Ldz(\sigma_*+iz)\Bigr)^k
\Bigl(-\Ldz(\sigma_*-iz+is)\Bigr)^k
\mathcal L(z,s)\,dz
\ll
HZ(\log X)^{2k}
\left(1+\min(1,|s|)\log X\right)^{-2},
\end{equation*}
uniformly for $s$ in the stated range. With the prefactor $1/H$, the intermediate-$s$ contribution is
\begin{align*}
&\ll
Z(\log X)^{2k}
\int_{(\log X)^{-1+\varepsilon}<|s|\le(\log X)^\eta}
\left(1+\min(1,|s|)\log X\right)^{-2}\,ds \\
&\ll
Z(\log X)^{2k}
\Bigl( (\log X)^{-1-\varepsilon}+(\log X)^{-2+\eta} \Bigr).
\end{align*}
As we chose $Z=(\log\log X)^A$ and we can fix $\eta<1$, this contribution is $o\bigl((\log X)^{2k-1}\bigr)$. To evaluate the remaining part, which is $|s|\le(\log X)^{-1+\varepsilon}$, we use $\kappa=\sigma_*-\frac12=\frac1{\log T}$ and apply Conjecture~\ref{eq:main-conj} with the earlier choices $z_1=\kappa$, $z_2=-\kappa-is$, and $z_1-z_2=2\kappa+is$, 
for a fixed $0<\xi<1$ if $k\geq 3$ and $0 < \xi < c$ for $k=1$ and $k=2$. For each
\begin{equation*}
\frac1{hZ}\leq Y\leq\frac Zh,
\end{equation*}
we apply Conjecture~\ref{eq:main-conj} with its averaging height
$T$ replaced by $Y$. To distinguish the auxiliary shift parameter
denoted by $Y$ in the statement of the conjecture from this averaging
height, denote it here by $\mathcal Y$ and take $\mathcal Y=1$. This
choice is admissible because, for all sufficiently large $X$,
\begin{equation*}
|\Im z_2|
=
|\Im(z_1-z_2)|
=
|s|
\leq(\log X)^{-1+\varepsilon}
<1.
\end{equation*}
Furthermore, $Y\leq Z/h=o(T)$ gives
\begin{equation*}
0<-\Re z_2=\frac1{\log T}\leq\frac3{\log Y},
\qquad
0<\Re(z_1-z_2)=\frac2{\log T}\leq\frac3{\log Y},
\end{equation*}
so the real-part conditions in Conjecture~\ref{eq:main-conj} are
also satisfied. For the cumulative moment $A_s(Y)$ defined in
\eqref{eq:shifted-moment-bound-A}, the conjecture gives
\begin{align}
\label{eq:central-conj-main-expanded}
A_s(Y)
&=
\frac{Y}{k!^2(2\pi i)^{2k}}
\int_{\mathcal C_u^k}\int_{\mathcal C_v^k}
Y^{\frac12(\sum_{j=1}^k u_j-\sum_{\ell=1}^k v_\ell)}
\prod_{j,\ell=1}^k(v_\ell-u_j-2\kappa-is)                         \notag\\
&\quad\times
\prod_{1\le i<j\le k}(u_i-u_j)^2
\prod_{1\le i<j\le k}(v_i-v_j)^2
\biggl(
\sum_{j=1}^k\frac1{u_j}
+
\sum_{\ell=1}^k\frac1{v_\ell-2\kappa-is}
-\frac{\log Y}{2}
\biggr)^k                                      \notag\\
&\quad\times
\biggl(
\sum_{j=1}^k\frac1{u_j+2\kappa+is}
+
\sum_{\ell=1}^k\frac1{v_\ell}
+\frac{\log Y}{2}
\biggr)^k
(2\kappa+is)^{-k^2}
\prod_{j=1}^k u_j^{-k}\,du_j
\prod_{\ell=1}^k v_\ell^{-k}\,dv_\ell                         \notag\\
&\quad+
O\!\left(
Y(\log Y)^{2k-1+\xi}
\left(1+\min(1,|s|)\log Y\right)^{-\theta_{k,c}}
\right).
\end{align}
 We now compare $\mathcal L(z,s)$ with $\mathcal L(z,0)$; for this
purpose, we define
\begin{equation*}
Q(w,u):=\frac{(1+u)^w-1}{w}.
\end{equation*}
After setting $x=Xv$ in \eqref{eq:L-kernel-def}, we have
\begin{equation}
\label{eq:L-kernel-scaled}
X^{-is}\mathcal L(z,s)
=
X^{2\sigma_*}
\int_{\mathbb R}
g(v)v^{2\sigma_*+is}
Q\!\left(\sigma_*+iz,\frac{h}{v}\right)
Q\!\left(\sigma_*-iz+is,\frac{h}{v}\right)\,dv .
\end{equation}
This form isolates the factor $X^{is}$ arising from
$x^{is}=X^{is}v^{is}$. On the support of $g$, we have $v\asymp1$,
so $v^{is}$ varies by only $O_g(|s|)$ from its value at $s=0$.
The second factor involving $Q$ also varies slowly with $s$. To make this precise, we note that $h/v\asymp h$ on the support of
$g$. Furthermore, in the central range, uniformly for
$|s|\leq(\log X)^{-1+\varepsilon}$, $|\sigma_*+iz|
\asymp
|\sigma_*-iz+is|
\asymp |z|$. As $\log(1+u)\ll h$ for $u\asymp h$, direct differentiation of
$Q(w,u)$ gives
\begin{align}
\label{eq:derivatives-of-Q}
Q(w,u)&\ll |z|^{-1},\\
\partial_wQ(w,u)&\ll h|z|^{-1}+|z|^{-2}, \notag \\
\partial_w^2Q(w,u)&\ll
h^2|z|^{-1}+h|z|^{-2}+|z|^{-3}, \notag
\end{align}
whenever $\Re{w} = \sigma_*$, $|w|\asymp|z|$ and $u \asymp h$. Differentiating the integrand in
\eqref{eq:L-kernel-scaled} with respect to $s$ gives
\begin{align*}
&\partial_s\!\left[
v^{is}Q\!\left(\sigma_*-iz+is,\frac{h}{v}\right)
\right] =
iv^{is}
\left[
(\log v)Q\!\left(\sigma_*-iz+is,\frac{h}{v}\right)
+
\partial_wQ\!\left(\sigma_*-iz+is,\frac{h}{v}\right)
\right].
\end{align*}
Using $\log v=O_g(1)$ and
\begin{equation*}
X^{2\sigma_*}
\int_{\mathbb R}|g(v)|v^{2\sigma_*}\,dv
\ll_g X,
\end{equation*}
we obtain, uniformly in the small-$s$ central range,
\begin{equation}
\partial_s\!\left(X^{-is}\mathcal L(z,s)\right)
\ll \frac{X}{z^2}.
\end{equation}
Differentiating this expression once with respect to $z$ and using
the bounds for $\partial_wQ$ and $\partial_w^2Q$ in \eqref{eq:derivatives-of-Q} similarly gives
\begin{equation}
\partial_z\partial_s
\!\left(X^{-is}\mathcal L(z,s)\right)
\ll
X\left(hz^{-2}+z^{-3}\right).
\end{equation}
Integrating these estimates with respect to $s$ from $0$ to $s$
proves
\begin{equation}
\label{eq:L-small-s-replacement}
\mathcal L(z,s)
=
X^{is}\mathcal L(z,0)
+
O\!\left(|s|\frac{X}{z^2}\right),
\end{equation}
and
\begin{equation}
\label{eq:L-small-s-replacement-derivative}
\left|
\partial_z\!\left(
\mathcal L(z,s)-X^{is}\mathcal L(z,0)
\right)
\right|
\ll
|s|X\left(hz^{-2}+z^{-3}\right).
\end{equation}
To estimate the error introduced in
\eqref{eq:expanded-central-region} by replacing
$\mathcal L(z,s)$ with $X^{is}\mathcal L(z,0)$, we apply partial summation with
respect to $A_s$, defined in \eqref{eq:shifted-moment-bound-A}, using
the bound \eqref{eq:corollary-bounds-central}. Write
\begin{equation*}
\mathcal{R}(z,s):=\mathcal L(z,s)-X^{is}\mathcal L(z,0).
\end{equation*}
The shifted logarithmic derivatives need not be estimated
pointwise: partial summation transfers the $z$-derivative to
$\mathcal{R}(z,s)$, for which
\eqref{eq:L-small-s-replacement} and
\eqref{eq:L-small-s-replacement-derivative} are available.
Applying partial summation in the form of
\eqref{eq:partial-summation-stieljes-IBP}, together with
\eqref{eq:corollary-bounds-central}, gives
\begin{align*}
&\left|
\int_{1/(hZ)}^{Z/h}
\Bigl(-\Ldz(\sigma_*+iz)\Bigr)^k
\Bigl(-\Ldz(\sigma_*-iz+is)\Bigr)^k
\mathcal{R}(z,s)\,dz
\right|                                                     \\
&\qquad\ll
|s|X(\log X)^{2k}
\left(1+\min(1,|s|)\log X\right)^{-2}    
\left(
hZ+\frac{h}{Z}
+h\int_{1/(hZ)}^{Z/h}\frac{dz}{z}
+\int_{1/(hZ)}^{Z/h}\frac{dz}{z^2}
\right)                                                       \\
&\qquad\ll
|s|HZ(\log X)^{2k}
\left(1+\min(1,|s|)\log X\right)^{-2}.
\end{align*}
Here the first two terms in parentheses arise from the two endpoints
in partial summation, while the final two arise from integrating the
two terms in
\eqref{eq:L-small-s-replacement-derivative}. The negative
$z$-range satisfies the same estimate by conjugation.

After including the outside factor $1/H$ and integrating over the
small-$s$ range, the error made by replacing $\mathcal L(z,s)$ with
$X^{is}\mathcal L(z,0)$ is
\begin{align*}
&\ll
Z(\log X)^{2k}
\int_{|s|\leq(\log X)^{-1+\varepsilon}}
|s|\left(1+|s|\log X\right)^{-2}\,ds                          \\
&=
\frac{2Z(\log X)^{2k}}{(\log X)^2}
\int_0^{(\log X)^\varepsilon}
\frac{v}{(1+v)^2}\,dv                                         \\
&\ll
Z(\log X)^{2k}
\frac{\log\log X}{(\log X)^2}
=
o\!\left((\log X)^{2k-1}\right).
\end{align*}

As $|s|\le(\log X)^{-1+\varepsilon}$, replacing the condition
$\frac1{hZ}\le |z-s|\le \frac Zh$ by $\frac1{hZ}\le |z|\le \frac Zh$ changes only boundary intervals of total length $O(|s|)$.  
On these boundary intervals, the uniform RH estimate
\eqref{eq:full-range-bound-titch}, applied with $\eta=1$, gives
\begin{equation*}
\left|\Ldz(\sigma_*+iz)\right|
+
\left|\Ldz(\sigma_*-iz+is)\right|
\ll(\log T)^2.
\end{equation*}
In addition, \eqref{eq:L-kernel-bounds}, with $\ell=0$ and $s=0$,
gives
\begin{equation*}
\mathcal L(z,0)
\ll
\frac{X}{|\sigma_*+iz|\,|\sigma_*-iz|}
\ll
\frac{X}{z^2}.
\end{equation*}
As the boundary intervals have total length $O(|s|)$,
\eqref{eq:full-range-bound-titch} and
\eqref{eq:L-kernel-bounds}, applied with $\ell=0$ and $s=0$, show
that their contribution to the central-range double integral obtained
by expanding \eqref{eq:V-contour-reduction} is
\begin{equation*}
O\!\left(hZ^2(\log X)^{4k-2+2\varepsilon}\right)
=
o\!\left((\log X)^{2k-1}\right).
\end{equation*}
Thus the condition on $|z-s|$ may be replaced by the condition on
$|z|$ at the cost of an admissible error, and we obtain
\begin{align*}
&\frac{1}{H(2\pi)^2}
\iint_{\frac1{hZ}\le |z|, |z-s| \le \frac Zh}
\Bigl(-\Ldz(\sigma_*+iz)\Bigr)^k
\Bigl(-\Ldz(\sigma_*-iz+is)\Bigr)^k
\mathcal L(z,s)\,dz\,ds  \\
&=\frac{1}{H(2\pi)^2}
\int_{|s|\le(\log X)^{-1+\varepsilon}}
X^{is}
\int_{\frac1{hZ}\le |z|\le \frac Zh}
\mathcal L(z,0)
\Bigl(-\Ldz(\sigma_*+iz)\Bigr)^k
\Bigl(-\Ldz(\sigma_*-iz+is)\Bigr)^k
dz ds  \\
&\qquad \qquad \qquad \qquad\qquad\qquad \qquad\qquad\qquad\qquad\qquad\qquad\qquad\qquad\qquad\quad+  o\bigl((\log X)^{2k-1}\bigr).  
\end{align*}
We insert \eqref{eq:central-conj-main-expanded} through partial summation in the $z$-variable, with weight $\mathcal L(z,0)$. The
endpoint and derivative terms are controlled by the kernel bounds \eqref{eq:kernel-bounds-L},
so the conjectural error term contributes at most
\begin{equation*}
\ll
Z(\log X)^{2k-1+\xi}
\int_{|s|\le(\log X)^{-1+\varepsilon}}
\left(1+\min(1,|s|)\log X\right)^{-\theta_{k,c}}\,ds,
\end{equation*}
and, using the definition of $\theta_{k,c}$ in \eqref{eq:theta_kc}, we have uniformly in all cases
\begin{equation*}
\int_{|s|\le(\log X)^{-1+\varepsilon}}
(1+|s|\log X)^{-\theta_{k,c}}\,ds
\ll
\frac1{\log X}.
\end{equation*}
Thus the contribution of the conjectural error term is
$o((\log X)^{2k-1})$ provided $0<\xi<1$ for $k\geq3$ and $0<\xi< c$ otherwise. We now treat the conjectural main term in
\eqref{eq:central-conj-main-expanded}. By \eqref{eq:Wk-def}, this
term is
\begin{equation*}
\frac{Y}{k!^2}(2\kappa+is)^{-k^2}
\mathcal W_k(2\kappa+is;Y).
\end{equation*}
The two-centre decomposition leading to
\eqref{eq:J-factor-revised}, together with
\eqref{eq:Jtilde-expansion-revised}--\eqref{eq:Jtilde-final-bound-revised},
gives
\begin{equation}
\label{eq:central-main-scale-derivative}
\left|
\frac{\partial}{\partial\log Y}
\left(
(2\kappa+is)^{-k^2}\mathcal W_k(2\kappa+is;Y)
\right)
\right|
\ll
(\log X)^{2k-1},
\end{equation}
uniformly for
\begin{equation*}
\frac1{hZ}\leq Y\leq\frac Zh,
\qquad
|s|\leq(\log X)^{-1+\varepsilon}.
\end{equation*}
As $|\log(hY)|\leq\log Z$ on this interval,
\eqref{eq:central-main-scale-derivative} gives
\begin{align}
\label{eq:central-main-freezing}
&\left|
(2\kappa+is)^{-k^2}\mathcal W_k(2\kappa+is;Y)
-
(2\kappa+is)^{-k^2}\mathcal W_k(2\kappa+is;1/h)
\right| \\
&\qquad\ll
(\log Z)(\log X)^{2k-1}. \notag
\end{align} 
We now evaluate the remaining $z$-integral of $\mathcal L(z,0)$.  For fixed $x$, set
\begin{equation*}
        K_x(z)
        :=
        \frac{(x+H)^{\sigma_*+iz}-x^{\sigma_*+iz}}
             {\sigma_*+iz} =
        \int_x^{x+H} y^{\sigma_*+iz-1}\,dy.
\end{equation*}
Hence, when $s=0$, the second factor in $\mathcal L(z,0)$ is $K_x(-z)$,
and therefore
\begin{equation*}
        \mathcal L(z,0)
        =
        \frac1X
        \int_{\mathbb R}g\!\left(\frac{x}{X}\right)
        K_x(z)K_x(-z)\,dx  =
        \frac1X
        \int_{\mathbb R}g\!\left(\frac{x}{X}\right)
        |K_x(z)|^2\,dx,
\end{equation*}
because $K_x(-z)=\overline{K_x(z)}$ for real $z$. By Mellin--Plancherel \cite[Theorem~71, pp.~94--95]{kn:titchmarsh-fourier},
\begin{equation*}
        \frac1{2\pi}\int_{\mathbb R}|K_x(z)|^2\,dz
        =
        \int_x^{x+H}y^{2\sigma_*-1}\,dy = Hx^{2\sigma_* - 1} + O(H h x^{2\sigma_*-1})
\end{equation*}
and consequently, by setting $x = Xv$, we get
\begin{align*}
\frac1{2\pi}\int_{\mathbb R}\mathcal L(z,0)\,dz
&=
\frac1X
\int_{\mathbb R}g\!\left(\frac{x}{X}\right)
\int_x^{x+H}y^{2\sigma_*-1}\,dy\,dx \\
&\qquad\qquad\qquad=
H X^{2\sigma_*-1}
\int_0^\infty g(v)v^{2\sigma_*-1}\,dv +
O\!\left(HhX^{2\sigma_*-1}\right),
\end{align*}
since $x\asymp X$ and $H=o(X)$. As $X^{2\sigma_*-1}=X^{2/\log T}=O(1)$ and $v^{2\sigma_*-1}=O_g(1)$ on the support of $g$, the error above is $o(H)$. We now justify restricting the Mellin--Plancherel integral to the
central range. From the integral representation of $K_x(z)$ and the
fact that $x\asymp X$ on the support of $g(x/X)$, we have
\begin{equation*}
|K_x(z)|
\leq
\int_x^{x+H}y^{\sigma_*-1}\,dy
\ll
HX^{\sigma_*-1}.
\end{equation*}
The quotient representation of $K_x(z)$ also gives, for $|z|\geq1$,
\begin{equation*}
|K_x(z)|
\leq
\frac{(x+H)^{\sigma_*}+x^{\sigma_*}}
     {|\sigma_*+iz|}
\ll
\frac{X^{\sigma_*}}{|z|}.
\end{equation*}
Using the definition of $\mathcal L(z,0)$, it follows that
\begin{equation*}
\mathcal L(z,0)
\ll
H^2X^{2\sigma_*-2},
\end{equation*}
for all real $z$, while
\begin{equation*}
\mathcal L(z,0)
\ll
X^{2\sigma_*}|z|^{-2},
\end{equation*}
for $|z|\geq1$. Therefore,
\begin{align*}
\int_{|z|<1/(hZ)}\mathcal L(z,0)\,dz
&\ll
\frac{H^2X^{2\sigma_*-2}}{hZ}
=
\frac{HX^{2\kappa}}{Z},\\
&\int_{|z|>Z/h}\mathcal L(z,0)\,dz
\ll
X^{2\sigma_*}
\int_{Z/h}^{\infty}z^{-2}\,dz
\ll
\frac{hX^{2\sigma_*}}{Z}
=
\frac{HX^{2\kappa}}{Z}.
\end{align*}
Since $X^{2\kappa}=O(1)$ and $Z\to\infty$, both expressions are
$o(H)$. Thus the full Mellin--Plancherel integral may be restricted
to
\begin{equation*}
\frac1{hZ}\leq |z|\leq\frac Zh
\end{equation*}
with an error $o(H)$. It follows that
\begin{align*}
\frac{1}{2\pi H}
\int_{\frac1{hZ}\le |z|\le \frac Zh}
\mathcal L(z,0)\,dz
&=
X^{2\kappa}
\int_0^\infty g(v)v^{2\kappa}\,dv
+
o(1) =
X^{2\kappa}
\left(\int_0^\infty g(v)\,dv+o(1)\right)
\end{align*}
as we have $v^{2\kappa}=1+o(1)$ uniformly on
$\operatorname{supp}g$; this is due to $g$ being compactly supported in $(0,\infty)$ and $\kappa=1/\log T=o(1)$. The conjectural main term enters the Stieltjes integral through
\begin{equation*}
d\!\left[
\frac{Y}{k!^2}(2\kappa+is)^{-k^2}
\mathcal W_k(2\kappa+is;Y)
\right].
\end{equation*}
By the product rule,
\eqref{eq:central-main-scale-derivative} and
\eqref{eq:central-main-freezing} give
\begin{align}
\label{eq:central-main-density-freezing}
&\frac{d}{dY}\!\left[
\frac{Y}{k!^2}(2\kappa+is)^{-k^2}
\mathcal W_k(2\kappa+is;Y)
\right] \\
&\qquad\qquad\quad=
\frac{1}{k!^2}(2\kappa+is)^{-k^2}
\mathcal W_k(2\kappa+is;1/h)
+
O\!\left(
(1+\log Z)(\log X)^{2k-1}
\right).
\notag
\end{align}
The representation of $\mathcal L(z,0)$ above shows that it is
non-negative, while the Mellin--Plancherel calculation gives
\begin{equation*}
\int_{\frac1{hZ}\leq |z|\leq Z/h}
\mathcal L(z,0)\,dz
\ll H.
\end{equation*}
After including the outside factor $1/H$ and integrating over
$|s|\leq(\log X)^{-1+\varepsilon}$, the error term in
\eqref{eq:central-main-density-freezing} contributes
\begin{equation*}
\ll
(1+\log Z)(\log X)^{2k-2+\varepsilon}
=
o\!\left((\log X)^{2k-1}\right).
\end{equation*}
Thus the normalised conjectural main term may be evaluated at
$Y=1/h=X/H$. Using \eqref{eq:central-conj-main-expanded}, with $\frac{X}{H} = X^{1/\alpha}$ as the central value of the $z$-scale, we have
\begin{align}
\label{eq:central-main-before-s-line-move}
&\left(\int_0^\infty g(v)\,dv+o(1)\right)
\frac{1}{2\pi} \cdot \frac{1}{k!^2(2\pi i)^{2k}}
\int_{|s|\le(\log X)^{-1+\varepsilon}}
\left(X^{1/\alpha}\right)^{\alpha(2\kappa+is)}                         \notag\\
&\quad\times
\int_{\mathcal C_u^k}\int_{\mathcal C_v^k}
\left(X^{1/\alpha}\right)^{\frac12(\sum_{j=1}^k u_j-\sum_{\ell=1}^k v_\ell)}
\prod_{j,\ell=1}^k(v_\ell-u_j-2\kappa-is) 
\prod_{i<j}(u_i-u_j)^2
\prod_{i<j}(v_i-v_j)^2                                      \notag\\
&\quad\times
\biggl(
\sum_{j=1}^k\frac1{u_j}
+
\sum_{\ell=1}^k\frac1{v_\ell-2\kappa-is}
-\frac12\log X^{1/\alpha}
\biggr)^k 
\biggl(
\sum_{j=1}^k\frac1{u_j+2\kappa+is}
+
\sum_{\ell=1}^k\frac1{v_\ell}
+\frac12\log X^{1/\alpha}
\biggr)^k
 \notag\\
&\quad\times (2\kappa+is)^{-k^2}
\prod_{j=1}^k u_j^{-k}\,du_j
\prod_{\ell=1}^k v_\ell^{-k}\,dv_\ell\,ds \, +
o\bigl((\log X)^{2k-1}\bigr).
\end{align}
We make the change of variables $w=2\kappa+is$ where $dw=i\,ds$. The original $s$-interval is thereby mapped to the finite vertical
segment
\begin{equation*}
\Re w=2\kappa,
\qquad
|\Im w|\leq(\log X)^{-1+\varepsilon}.
\end{equation*}
We move this segment to $\Re w=1/\log X$ and then extend it to the
full vertical line. For this estimate, we apply Lemma~\ref{lem:P-to-Q-form} and
Lemma~\ref{lem:two-centre} termwise on the finite segment. The
factorisation \eqref{eq:J-factor-revised}, with $T$ there replaced
by $X^{1/\alpha}$, together with
\eqref{eq:Jtilde-zero-revised}--\eqref{eq:Jtilde-final-bound-revised},
shows that the terms with $a<r$ vanish and that the absolute value
of each remaining $(r,a)$-integrand satisfies
\begin{equation}
\label{eq:central-w-integrand-decay}
\ll_k
\begin{cases}
\bigl(\log X^{1/\alpha}\bigr)^k|w|^{-k},
    & r=0,\quad k\geq2,\\[2mm]
\bigl(\log X^{1/\alpha}\bigr)^{2k-2r}|w|^{-2r},
    & 1\leq r\leq k,
\end{cases}
\end{equation}
whenever $|\Im w|\geq(\log X)^{-1+\varepsilon}$. These estimates
hold uniformly in the strip joining the two vertical segments,
because $\Re w\asymp1/\log X$ there. For $k=1$ and $r=0$, the identities $\tilde J_0^{(0)}(w;X^{1/\alpha})=0$ and $\tilde J_0^{(1)}(w;X^{1/\alpha})=1$ show that the surviving integrand is $O(|w|^{-2})$.

It follows from \eqref{eq:central-w-integrand-decay} that the omitted
parts of the vertical line have total contribution
\begin{equation*}
\ll_k
\begin{cases}
\bigl(\log X^{1/\alpha}\bigr)^k
(\log X)^{(1-\varepsilon)(k-1)},
    & r=0,\quad k\geq2,\\[2mm]
\bigl(\log X^{1/\alpha}\bigr)^{2k-2r}
(\log X)^{(1-\varepsilon)(2r-1)},
    & 1\leq r\leq k,\\[2mm]
(\log X)^{1-\varepsilon},
    & r=0,\quad k=1,
\end{cases}
=
o\!\left((\log X)^{2k-1}\right).
\end{equation*}
The horizontal segments introduced when moving
$\Re w=2\kappa$ to $\Re w=1/\log X$ have length
$O(1/\log X)$. Applying
\eqref{eq:central-w-integrand-decay} at
$|\Im w|=(\log X)^{-1+\varepsilon}$ shows that their total
contribution is $o\!\left((\log X)^{2k-1}\right)$. Both vertical segments lie to the right of $w=0$, so no pole is
crossed. The finite segment may therefore be moved to
$\Re w=1/\log X$ and extended to the full vertical line at a cost
of $o((\log X)^{2k-1})$. After incorporating the Jacobian and
renaming $w$ as $s$, \eqref{eq:central-main-before-s-line-move}
becomes
\begin{align}
\label{eq:central-main-BC-form}
&\left(\int_0^\infty g(v)\,dv+o(1)\right)
\frac{1}{k!^2(2\pi i)^{2k+1}}
\int_{\Re s=\frac{1}{\log{X}}}
\int_{\mathcal{C}_{u}^k}
\int_{\mathcal{C}_{v}^k}
s^{-k^2} \left(X^{1/\alpha}\right)^{
\alpha s+\frac12(\sum_{j=1}^k u_j-\sum_{\ell=1}^k v_\ell)} \notag\\
&\quad\times
\biggl(
\sum_{j=1}^k\frac1{u_j}
+
\sum_{\ell=1}^k\frac1{v_\ell-s}
-\frac12\log X^{1/\alpha}
\biggr)^k
\biggl(
\sum_{j=1}^k\frac1{u_j+s}
+
\sum_{\ell=1}^k\frac1{v_\ell}
+\frac12\log X^{1/\alpha}
\biggr)^k                                    \notag\\
&\quad\times \prod_{j,\ell=1}^k(v_\ell-u_j-s)
\Delta(u)^2\Delta(v)^2 
\prod_{j=1}^k u_j^{-k}\,du_j
\prod_{\ell=1}^k v_\ell^{-k}\,dv_\ell\,ds
+
o \bigl((\log X)^{2k-1}\bigr).
\end{align}
Here $\mathcal{C}_{u}$ encloses $0$ and $-s$, while
$\mathcal{C}_{v}$ encloses $0$ and $s$, with no other singularities enclosed. We now apply Lemma~\ref{lem:P-to-Q-form} to \eqref{eq:central-main-BC-form}, which gives
\begin{align}
\label{eq:central-main-Q-form}
&\left(\int_0^\infty g(v)\,dv+o(1)\right)
\frac{(-1)^{k^2}}{k!^2(2\pi i)^{2k+1}}
\int_{\Re s=\frac{1}{\log{X}}}
s^{-k^2}
\int_{\mathcal{C}_{u}^k}
\int_{\mathcal{C}_{v}^k}
\left(X^{1/\alpha}\right)^{\alpha s+\sum_{j=1}^k u_j}
                                                        \notag\\
&\quad\times
\biggl(
\sum_{j=1}^k\frac1{u_j}
+
\sum_{\ell=1}^k\frac1{v_\ell-s}
-\log X^{1/\alpha}
\biggr)^k
\biggl(
\sum_{j=1}^k\frac1{u_j+s}
+
\sum_{\ell=1}^k\frac1{v_\ell}
\biggr)^k                                      \notag\\
&\quad\times \prod_{j,\ell=1}^k(u_j-v_\ell+s)  \Delta(u)^2\Delta(v)^2
\prod_{j=1}^k u_j^{-k}\,du_j
\prod_{\ell=1}^k v_\ell^{-k}\,dv_\ell\,ds
+
o \bigl((\log X)^{2k-1} \bigr).
\end{align}
In the above, we used $\prod_{j,\ell=1}^k(v_\ell-u_j-s)
        =
        (-1)^{k^2}
        \prod_{j,\ell=1}^k(u_j-v_\ell+s)$.
Next, Lemma~\ref{lem:two-centre} is applied; for each fixed $s$ on $\Re s=\frac{1}{\log{X}}$, take $\sigma=-s$, $L=\log X^{1/\alpha}$ and $\Phi_L(u_1,\ldots,u_k)=\left(X^{1/\alpha}\right)^{u_1+\cdots+u_k}$. The holomorphy and symmetry
hypotheses of Lemma~\ref{lem:two-centre} are immediate from this choice, and we obtain
\begin{align}
\label{eq:central-main-after-two-centre}
&\left(\int_0^\infty g(v)\,dv+o(1)\right)
\frac{1}{k!^2(2\pi i)^{2k+1}}
\sum_{r=0}^k\binom{k}{r}^2
\sum_{a=0}^k\binom{k}{a}
\left(\log X^{1/\alpha}\right)^{k-a}                 \notag\\
&\quad\quad\qquad\qquad\times
\int_{\Re s=1/\log X}
\left(X^{1/\alpha}\right)^{\alpha s}
s^{-a-k}
J_{\log X^{1/\alpha},r}^{(a)}(-s)\,ds
+
o\bigl((\log X)^{2k-1}\bigr),
\end{align}
where
\begin{align}
\label{eq:J-applied-form}
J_{\log X^{1/\alpha},r}^{(a)}(-s)
&=
\int_{c_j(r)}\int_{\hat c_\ell(r)}
\left(X^{1/\alpha}\right)^{-s(x_1+\cdots+x_k)}
\Delta(x)^2\Delta(y)^2
\prod_{j,\ell=1}^k(x_j-y_\ell-1)                     \notag\\
&\quad\times
\biggl(
\sum_{j=1}^k\frac1{x_j}
+
\sum_{\ell=1}^k\frac1{y_\ell+1}
\biggr)^a
\biggl(
\sum_{j=1}^k\frac1{x_j-1}
+
\sum_{\ell=1}^k\frac1{y_\ell}
\biggr)^k \prod_{j=1}^k x_j^{-k}\,dx_j
\prod_{\ell=1}^k y_\ell^{-k}\,dy_\ell .
\end{align}
The contours $c_j(r)$ and $\hat c_\ell(r)$ are as defined in \eqref{eq:prk}. We now evaluate the $s$-integral in
\eqref{eq:central-main-after-two-centre}. Substituting
\eqref{eq:J-applied-form} into this expression, the $s$-dependent
factor in the $(r,a)$-term is
\begin{equation*}
\left(X^{1/\alpha}\right)^{
s\left(\alpha-\sum_{j=1}^k x_j\right)}s^{-a-k}.
\end{equation*}
When the contour shift crosses $s=0$, the corresponding residue is
\begin{align}
\label{eq:central-elementary-s-residue}
&\operatorname*{Res}_{s=0}
\left[
\left(X^{1/\alpha}\right)^{
s\left(\alpha-\sum_{j=1}^k x_j\right)}
s^{-a-k}
\right] =
\frac{\left(\log X^{1/\alpha}\right)^{a+k-1}}
     {(a+k-1)!}
\left(\alpha-\sum_{j=1}^k x_j\right)^{a+k-1}.
\end{align}
Multiplication by the factor
$\left(\log X^{1/\alpha}\right)^{k-a}$ already present in
\eqref{eq:central-main-after-two-centre} gives the common scale
$\left(\log X^{1/\alpha}\right)^{2k-1}$. To identify the remaining coefficient, differentiate the
$\alpha$-dependent factor in \eqref{eq:prk}. The contours in \eqref{eq:prk} are independent of $\alpha$, and the
integrand depends polynomially on $\alpha$, so the differentiation
may be carried out under the contour integrals. All the remaining
factors in \eqref{eq:J-applied-form} agree with those in the
corresponding summand of \eqref{eq:prk}. Hence, after inserting
\eqref{eq:central-elementary-s-residue} into
\eqref{eq:central-main-after-two-centre} and summing over $a$, the
$r$-th residue contribution is
\begin{equation*}
\left(\log X^{1/\alpha}\right)^{2k-1}P'_{r,k}(\alpha).
\end{equation*}
The vanishing and degree bounds \eqref{eq:Jtilde-zero-revised} and
\eqref{eq:Jtilde-final-bound-revised} apply in the same way. Thus the contour
shift selection rule from \eqref{eq:z-residue-selection-finite} gives a
contribution only from $r<\alpha$, while the terms with $r>\alpha$ are
negligible. If $\alpha$ is an integer, the transition block $r=\alpha$ is
negligible by the same transition-block argument following \eqref{eq:not-log}. 
Applying \eqref{eq:central-elementary-s-residue} and summing over the
surviving values $r<\alpha$ gives
\begin{equation}
\label{eq:central-main-residue-evaluation}
        \frac{1}{k!^2}
        \sum_{0\le r<\alpha}
        \binom{k}{r}^2
        P'_{r,k}(\alpha)
        \left(\log X^{1/\alpha}\right)^{2k-1}
        +
        o \bigl((\log X)^{2k-1} \bigr).
\end{equation}
In particular, this proves
\begin{equation*}
        \frac1{HX}\int_{\mathbb R}g\!\left(\frac{x}{X}\right)
        |\mathcal P_k(x;H)|^2\,dx
        \ll
        (\log X)^{2k-1},
\end{equation*}
which in turn proves \eqref{eq:V-contour-reduction}. Hence, by \eqref{eq:back-to-reality}, we have
\begin{equation}
\label{eq:final-V}
\tilde V_{g,k}(X,H)
=
\frac1{HX}\int_{\mathbb R}g\!\left(\frac{x}{X}\right)
|\mathcal P_k(x;H)|^2\,dx
+
o \bigl((\log X)^{2k-1}\bigr).
\end{equation}
Therefore, we can now write the variance with a separated dependence on $g(v)$. Combining \eqref{eq:central-main-residue-evaluation} with the factor
$\int_0^\infty g(v)\,dv+o(1)$ in \eqref{eq:central-main-after-two-centre}, the
central contribution is
\begin{equation}
\label{eq:central-contribution-final}
\left(\int_0^\infty g(v)\,dv\right)
\frac{1}{k!^2}
\sum_{0\le r<\alpha}
\binom{k}{r}^{2}
P'_{r,k}(\alpha)
\left(\log X^{1/\alpha}\right)^{2k-1}
+
o \bigl((\log X)^{2k-1}\bigr).
\end{equation}
Due to the established negligibility of $E_1$, $E_2$ and
$E_3$, we obtain
\begin{align}
\label{eq:Pk-mean-square-asymptotic}
\frac1{HX}\int_{\mathbb R}g\!\left(\frac{x}{X}\right)
|\mathcal P_k(x;H)|^2\,dx
&=
\left(\int_0^\infty g(v)\,dv\right)
\frac{1}{k!^2}
\sum_{0\le r<\alpha}
\binom{k}{r}^{2}
P'_{r,k}(\alpha)
\left(\log X^{1/\alpha}\right)^{2k-1}        \notag\\
&\qquad\qquad\qquad\qquad\qquad\qquad \quad\quad+
o\bigl((\log X)^{2k-1}\bigr).
\end{align}
Combining the above with \eqref{eq:final-V} gives the desired asymptotic
for $\tilde V_{g,k}(X,H)$ and proves the claimed derivative relationship in \eqref{conj:derivative}. In particular, we have shown \eqref{eq:final-statement-to-prove} and, to then relate the variance to $\delta_k(\alpha)$, we can use \eqref{def:delta-derivative-Prk} to write
\begin{equation}
\lim_{X \rightarrow \infty} \frac{\mathcal{V}_{k}(X,X^{1-1/\alpha})}{(\log{X^{1/\alpha}})^{2k-1}} = \mathcal{D}_{k}(\alpha) = \delta_{k}(\alpha).
\end{equation}
As Theorem~\ref{mkalpha} requires stronger hypotheses for the parameter $Y(T)$, we may only write $\mathcal{D}_{k}(\alpha) = \mathcal{M}'_{k}(\alpha)$ under its assumptions.

\subsection{Proof of Theorem ~\ref{conj:arithmetic-variance}} 

Let $k\ge4$ and $\alpha>0$, and suppose that $\log X/\log q\to\alpha$ as
$q\to\infty$ through the primes. We assume GRH for Dirichlet $L$-functions
modulo $q$, and assume Conjecture~\ref{eq:ratios-arithmetic-corrected} with height $Y(q)=q^{\alpha/2 +\varepsilon}$ for some fixed $\varepsilon>0$. Choose $\eta>0$ such that $\alpha\eta<\varepsilon$ and let
$T=X^{1/2+\eta}$. Then $T=q^{\alpha/2+\alpha\eta+o(1)}$, and therefore
$2T\le Y(q)$ for all sufficiently large $q$. Let $\chi_0$ denote the principal character modulo $q$. For a Dirichlet character $\chi\bmod q$, define
\begin{equation}
\label{eq:ap-Ichi-def}
I_\chi(X)
:=
\sum_{n\le X}\Lambda^{*k}(n)\chi(n)
-
\mathds 1_{\chi=\chi_0}
\operatorname*{Res}_{s=1}
\left[
\left(-\frac{\zeta_q'}{\zeta_q}(s)\right)^k\frac{X^s}{s}
\right].
\end{equation}
For $\Re s>1$, we have
\begin{equation}
\label{eq:ap-Lseries-coefficients}
        \left(-\frac{L'}{L}(s,\chi)\right)^k
        =
        \sum_{n\ge1}\frac{\Lambda^{*k}(n)\chi(n)}{n^s}.
\end{equation}
As $L(s,\chi_0)=\zeta_q(s)$, character orthogonality gives
\begin{equation}
\label{eq:ap-character-variance}
        \mathcal W_k(X;q)
        =
        \frac1{\varphi(q)}
        \sum_{\chi\bmod q}
        \left|X^{-1/2}I_\chi(X)\right|^2,
\end{equation}
with $\mathcal W_k(X;q)$ defined in \eqref{def:W_k}. We now replace $X^{-1/2}I_\chi(X)$ by a truncated Perron integral. Applying the truncated Perron formula on $\Re s=1+\kappa$,
where $\kappa=1/\log q$, at height $T$ gives
\begin{align}
\label{eq:ap-perron-first-line}
X^{-1/2}I_\chi(X)
&=
\frac1{2\pi i}
\int_{1+\kappa-iT}^{1+\kappa+iT}
\left(-\frac{L'}{L}(s,\chi)\right)^k
\frac{X^{s-1/2}}{s}\,ds                                      \notag\\
&\qquad \quad \quad\quad
-
\mathds 1_{\chi=\chi_0}X^{-1/2}
\operatorname*{Res}_{s=1}
\left[
\left(-\frac{\zeta_q'}{\zeta_q}(s)\right)^k\frac{X^s}{s}
\right]
+
E_{\chi,\mathrm{tr}}(X;T).
\end{align}
Under GRH, uniformly in $\chi\bmod q$, we have
\begin{equation}
\label{eq:ap-GRH-logderivative-bound}
\frac{L'}{L}(\sigma+it,\chi)
\ll
\frac{\log(q(|t|+3))}{\sigma-\frac12},
\qquad
\frac12<\sigma\le2,\quad
|\sigma+it-1|\ge\kappa.
\end{equation}
For the principal character, the pole contribution is
$O(|\sigma+it-1|^{-1})=O(\log q)$ in this region and is therefore
absorbed by the bound in \eqref{eq:ap-GRH-logderivative-bound}. For sufficiently large $q$, both vertical contours and the horizontal
segments lie in the region
$\frac12<\Re s\le2$, $|s-1|\ge\kappa$.

The coefficient bound
$0\le\Lambda^{*k}(n)\le(\log n)^k$, which follows inductively from
$\sum_{d\mid n}\Lambda(d)=\log n$, allows us to apply the uniform
truncated Perron estimate. With our choices of $\kappa$ and $T$,
this gives, uniformly in $\chi\bmod q$,
\begin{equation}
\label{eq:ap-truncation-error}
E_{\chi,\mathrm{tr}}(X;T)
\ll_{k,\alpha,\eta}
X^{-1/2}(\log X)^k
+
\frac{X^{1/2+\kappa}}{T}
\left((\log q)\log(qT)\right)^k.
\end{equation}
The first term accounts for the possible endpoint and near-endpoint
contributions at the sharp cutoff $X$. We shift the integral in \eqref{eq:ap-perron-first-line} to
$\Re s=\frac12+\kappa$. By GRH, no zero of $L(s,\chi)$ is crossed. If
$\chi=\chi_0$, the pole at $s=1$ is crossed, and its residue is exactly the
residue subtracted in \eqref{eq:ap-perron-first-line}. On the horizontal
segments $s=\sigma\pm iT$, $\frac12+\kappa\le\sigma\le1+\kappa$,
\eqref{eq:ap-GRH-logderivative-bound} gives
\begin{equation}
\label{eq:ap-horizontal-log-bound}
        \frac{L'}{L}(\sigma\pm iT,\chi)
        \ll
        (\log q)\log(qT).
\end{equation}
Since $|s|\asymp T$ on these segments, the horizontal contribution after the
normalisation by $X^{-1/2}$ is
\begin{align}
\label{eq:ap-horizontal-error}
\ll
\int_{\frac12+\kappa}^{1+\kappa}
\left((\log q)\log(qT)\right)^k
\frac{X^{\sigma-1/2}}{T}\,d\sigma                                     \ll
\frac{X^{1/2+\kappa}}{T}
\left((\log q)\log(qT)\right)^k.
\end{align}
Because $T=X^{1/2+\eta}$ and $X^\kappa=O_\alpha(1)$, both
\eqref{eq:ap-truncation-error} and \eqref{eq:ap-horizontal-error} are
$O_{k,\alpha,\eta}(X^{-\delta})$ for some $\delta>0$. Writing $s=\frac12+z$,
we have
\begin{equation}
\label{eq:ap-perron-reduction}
        X^{-1/2}I_\chi(X)=I_{\chi,T}(X)+R_\chi(X;T),
        \qquad
        R_\chi(X;T)\ll_{k,\alpha,\eta}X^{-\delta},
\end{equation}
where
\begin{equation}
\label{eq:ap-truncated-perron}
I_{\chi,T}(X)
:=
\frac1{2\pi i}
\int_{\kappa-iT}^{\kappa+iT}
\left(-\frac{L'}{L}\left(\frac12+z,\chi\right)\right)^k
\frac{X^z}{\frac12+z}\,dz .
\end{equation}
We next remove $R_\chi(X;T)$ from the average. Put
\begin{equation}
\label{eq:ap-Iq-def}
        \mathcal I_q
        :=
        \frac1{\varphi(q)}
        \sum_{\chi\bmod q}|I_{\chi,T}(X)|^2.
\end{equation}
By \eqref{eq:ap-character-variance} and \eqref{eq:ap-perron-reduction},
\begin{equation}
\label{eq:ap-expanded-with-error}
        \mathcal W_k(X;q)
        =
        \frac1{\varphi(q)}
        \sum_{\chi\bmod q}|I_{\chi,T}(X)+R_\chi(X;T)|^2.
\end{equation}
Using $|A+B|^2=|A|^2+2\Re(A\overline B)+|B|^2$, we obtain
\begin{align}
\label{eq:ap-cauchy-difference}
\mathcal W_k(X;q)-\mathcal I_q
&=
\frac{2}{\varphi(q)}
\Re\sum_{\chi\bmod q}I_{\chi,T}(X)\overline{R_\chi(X;T)}
+
\frac1{\varphi(q)}
\sum_{\chi\bmod q}|R_\chi(X;T)|^2.
\end{align}
By the Cauchy--Schwarz inequality,
\begin{align}
\label{eq:ap-cauchy-bound}
\left|\mathcal W_k(X;q)-\mathcal I_q\right|
&\le
2\mathcal I_q^{1/2}
\left(
\frac1{\varphi(q)}
\sum_{\chi\bmod q}|R_\chi(X;T)|^2
\right)^{1/2}
+
\frac1{\varphi(q)}
\sum_{\chi\bmod q}|R_\chi(X;T)|^2 .
\end{align}
It remains only to give a crude bound for $\mathcal I_q$.
On the segment $z=\kappa+iv$, $|v|\le T$,
\eqref{eq:ap-GRH-logderivative-bound} gives
\begin{equation}
\label{eq:ap-line-log-bound}
        \frac{L'}{L}\left(\frac12+z,\chi\right)
        \ll
        (\log q)\log(qT)
        \ll_{\alpha,\eta}
        (\log q)^2.
\end{equation}
Also $X^\kappa=O_\alpha(1)$. Hence
\begin{align}
\label{eq:ap-Ichi-crude-bound}
I_{\chi,T}(X)
&\ll
(\log q)^{2k}
\int_{-T}^{T}\frac{dv}{\left|\frac12+\kappa+iv\right|}
\ll
(\log q)^{2k+1}.
\end{align}
Therefore $\mathcal I_q\ll(\log q)^{4k+2}$. Combining this bound with
\eqref{eq:ap-perron-reduction} in \eqref{eq:ap-cauchy-bound} gives
\begin{equation}
\label{eq:ap-character-variance-truncated}
        \mathcal W_k(X;q)
        =
        \mathcal I_q
        +
        o\left((\log q)^{2k-1}\right).
\end{equation}

Expanding the square in \eqref{eq:ap-Iq-def}, and writing the conjugate
integral on the line $\Re z_2=-\kappa$, gives
\begin{align}
\label{eq:ap-double-integral}
\mathcal W_k(X;q)
&=
\frac1{(2\pi i)^2}
\int_{\kappa-iT}^{\kappa+iT}
\int_{-\kappa-iT}^{-\kappa+iT}
\frac{X^{z_1-z_2}}
{\left(\frac12+z_1\right)\left(\frac12-z_2\right)}
\mathcal B_k(z_1,z_2;q)\,dz_2\,dz_1        \notag\\
&\quad+
o\left((\log q)^{2k-1}\right),
\end{align}
where
\begin{equation}
\label{eq:ap-Bk-def}
\mathcal B_k(z_1,z_2;q)
:=
\frac1{\varphi(q)}
\sum_{\chi\bmod q}
\left(\frac{L'}{L}\left(\frac12+z_1,\chi\right)\right)^k
\left(\frac{L'}{L}\left(\frac12-z_2,\overline\chi\right)\right)^k .
\end{equation}
On the contours in \eqref{eq:ap-double-integral}, with $z=z_1-z_2$, we have
\begin{equation}
\label{eq:ap-conjecture-region}
        \Re z=2\kappa=\frac2{\log q},
        \qquad
        |\Im z_2|\le T,
        \qquad
        |\Im z|\le2T.
\end{equation}
Thus Conjecture~\ref{eq:ratios-arithmetic-corrected} applies throughout
the integration region in \eqref{eq:ap-double-integral}, since
$2T\le Y(q)$ for all sufficiently large $q$.
Choose $0<\xi<1$. Since $k\ge4$, we have $\theta_{k,c}=k-1$.
Let $E_{\mathrm{conj}}$ denote the contribution of the conjectural
remainder to \eqref{eq:ap-double-integral}. Writing
$z_1=\kappa+i(v+y)$ and $z_2=-\kappa+iv$, and putting
$a=\frac12+\kappa$, we use
\begin{equation}
\int_{\mathbb R}
\frac{dv}{|a+i(v+y)|\,|a-iv|}
\ll
\frac{\log(2+|y|)}{1+|y|}.
\end{equation}
Here we have enlarged the inner integration interval to
$\mathbb R$ after taking absolute values.
Since $|X^{z_1-z_2}|=X^{2\kappa}=O_\alpha(1)$, it follows that
\begin{align}
\label{eq:ap-ratios-error-negligible}
|E_{\mathrm{conj}}|
&\ll
(\log q)^{2k-1+\xi}
\int_{-2T}^{2T}
\frac{\log(2+|y|)}{1+|y|}
\left(1+\min(|2\kappa+iy|,1)\log q\right)^{-(k-1)}
\,dy \notag\\
&\ll
(\log q)^{2k-1+\xi}
\left(
\frac1{\log q}
+\frac{(\log(2T))^2}{(\log q)^{k-1}}
\right) \notag\\
&\ll
(\log q)^{2k-2+\xi}
+(\log q)^{k+2+\xi}
=o\!\left((\log q)^{2k-1}\right).
\end{align}
The second bound follows by splitting the $y$-integral at
$|y|=1$. The final estimates use $\log T\asymp\log q$,
$k\ge4$, and $\xi<1$. We now insert the conjectural main term. Changing variables from $(z_1,z_2)$
to $(z,z_2)$, the inner $z_2$-integral is the same finite Perron kernel as in
\eqref{eq:perron-kernel-estimate-final}, with the shifted denominator
$\left(\frac12+z+z_2\right)\left(\frac12-z_2\right)$. Hence, uniformly for $z=2\kappa+iy$ with $|y|\le2T$,
\begin{equation}
\label{eq:ap-z2-kernel}
\frac1{2\pi i}
\int_{-\kappa+i\max(-T,-T-y)}^{-\kappa+i\min(T,T-y)}
\frac{dz_2}
{\left(\frac12+z+z_2\right)\left(\frac12-z_2\right)}
=
\frac1{1+z}
+
O\left(\frac{\log(qT)}{T}\right).
\end{equation}
The error term in \eqref{eq:ap-z2-kernel} gives
$o((\log q)^{2k-1})$ after insertion into the remaining $z$-integral. Therefore
\begin{align}
\label{eq:ap-single-z-main}
\mathcal W_k(X;q)
&=
\frac1{2\pi i}
\int_{2\kappa-2iT}^{2\kappa+2iT}
\frac{X^z}{1+z}\mathcal K_k(z;q)\,dz
+
o\left((\log q)^{2k-1}\right),
\end{align}
where
\begin{align}
\label{eq:ap-Kk-def}
\mathcal K_k(z;q)
&:=
\frac1{k!^2(2\pi i)^{2k}}
\int_{\mathcal C_u^k}\int_{\mathcal C_v^k}
q^{\frac12\left(\sum_{j=1}^k u_j-\sum_{\ell=1}^k v_\ell\right)}
\Delta(u)^2\Delta(v)^2                                      \notag\\
&\quad\times
\left(
\sum_{j=1}^k\frac1{u_j}
+
\sum_{\ell=1}^k\frac1{v_\ell-z}
-
\frac{\log q}{2}
\right)^k
\left(
\sum_{j=1}^k\frac1{u_j+z}
+
\sum_{\ell=1}^k\frac1{v_\ell}
+
\frac{\log q}{2}
\right)^k                                      \notag\\
&\quad\times
\prod_{j,\ell=1}^k(v_\ell-u_j-z)
z^{-k^2}
\prod_{j=1}^k u_j^{-k}\,du_j
\prod_{\ell=1}^k v_\ell^{-k}\,dv_\ell .
\end{align}
The two-centre estimates
\eqref{eq:Jtilde-zero-revised} and
\eqref{eq:Jtilde-final-bound-revised} give
\begin{equation}
|\mathcal K_k(2\kappa+iy;q)|
\ll_k(\log q)^{2k}(1+|y|\log q)^{-2}.
\end{equation}
Its absolute integral over $|y|\le2T$ is therefore
$O_k((\log q)^{2k-1})$. Since $X^{2\kappa}=O_\alpha(1)$
and $\log(qT)\asymp\log q$, the error in
\eqref{eq:ap-z2-kernel} contributes
\begin{equation}
O_{k,\alpha,\eta}\!\left(\frac{(\log q)^{2k}}{T}\right)
=o\!\left((\log q)^{2k-1}\right),
\end{equation}
as required in \eqref{eq:ap-single-z-main}. Make the changes of variables $z=w/\log q$ and $u_j \mapsto u_j/ \log{q}$, $v_\ell \mapsto v_\ell / \log{q}$. Put
\begin{equation}
L:=\log q,\qquad
H_k(w):=L^{-2k}\mathcal K_k(w/L;q).
\end{equation}
The scaled inner integral shows that $H_k$ is independent of $q$,
and the bound established above gives
\begin{equation}
|H_k(2+it)|\ll_k(1+|t|)^{-2}.
\end{equation}
Since $|1+w/L|^{-1}\le1$ on $\Re w=2$ and
$\log X/L\to\alpha$, dominated convergence yields
\begin{align}
\label{eq:ap-scaled-contour-identification}
\frac{\mathcal W_k(X;q)}{L^{2k-1}}
&=
\frac1{2\pi i}
\int_{2-2iTL}^{2+2iTL}
\frac{e^{w\log X/L}H_k(w)}{1+w/L}\,dw+o(1)
\notag\\
&\longrightarrow
\frac1{2\pi i}\int_{(2)}e^{\alpha w}H_k(w)\,dw.
\end{align}

For non-integer $\alpha$, Lemmas~\ref{lem:P-to-Q-form}
and~\ref{lem:two-centre}, followed by the residue calculation
in \eqref{eq:central-main-residue-evaluation}, evaluate this
limiting integral. The integral is continuous in $\alpha$
by the bound above, and the polynomial expression is continuous
at the positive integers since $P'_{m,k}(m)=0$ for $1\le m\le k$.
Consequently, for every $\alpha>0$,
\begin{equation}
\label{eq:ap-delta-definition}
\frac1{2\pi i}\int_{(2)}e^{\alpha w}H_k(w)\,dw
=
\frac1{k!^2}
\sum_{\substack{0\le r\le k\\r<\alpha}}
\binom{k}{r}^2P'_{r,k}(\alpha)
=
\delta_k(\alpha).
\end{equation}
Combining \eqref{eq:ap-scaled-contour-identification} and
\eqref{eq:ap-delta-definition}, we obtain
\begin{equation}
\label{eq:ap-W-final-asymptotic}
        \mathcal W_k(X;q)
        =
        (\log q)^{2k-1}\delta_k(\alpha)
        +
        o\left((\log q)^{2k-1}\right).
\end{equation}
Dividing by $(\log q)^{2k-1}$ and letting $q\to\infty$ through the primes gives the stated result. 

The restriction $k\ge4$ arises from the estimate for the
conjectural remainder. For $k=1,2,3$, the same argument applies if the decay exponent
$\theta_{k,c}$ in Conjecture~\ref{eq:ratios-arithmetic-corrected}
satisfies $\theta_{k,c}>2+\xi$. In particular, taking
$\theta_{k,c}=2+c$ with $0<\xi<c\le1$ extends the result to these
cases under the same GRH and height assumptions.

\section{Acknowledgments}
We would like to thank the American Institute of Mathematics for hosting several visits of the authors during the early stages of this project. We would also like to thank Sandro Bettin for answering questions related to this topic. This group was supported by NSF Grant DMS-185439. The third author was supported by the Heilbronn Institute for Mathematical Research. 

\end{document}